\documentclass[a4paper]{amsart}
\usepackage{amsmath,amssymb,amsfonts}
\usepackage{colordvi}
\usepackage[all]{xy}

\def\mpoint{\;.}
\def\mvirg{\;,}

\def\sou{\underline}

\def\mpn{\medskip\par\noindent}

\def\mmpn{\vskip 1em minus 1em\par\noindent}

\def\sp{\bigskip\par}

\def\smp{\smallskip\par}

\def\CB{{\mathcal B}}
\def\CC{{\mathcal C}}
\def\CD{{\mathcal D}}
\def\CE{{\mathcal E}}
\def\CF{{\mathcal F}}

\def\CI{{\mathcal I}}
\def\CL{{\mathcal L}}

\def\CP{{\mathcal P}}
\def\CQ{{\mathcal Q}}
\def\CR{{\mathcal R}}

\def\CY{{\mathcal Y}}

\def\Ker{\operatorname{Ker}\nolimits}

\def\Im{\operatorname{Im}\nolimits}
\def\Id{\operatorname{id}\nolimits}

\def\Mod{\operatorname{Mod}\nolimits}
\def\Hom{\operatorname{Hom}\nolimits}
\def\End{\operatorname{End}\nolimits}
\def\Aut{\operatorname{Aut}\nolimits}

\def\Tor{\operatorname{Tor}\nolimits}

\def\Irr{\operatorname{Irr}\nolimits}

\def\tot{\operatorname{tot}\nolimits}

\def\op{^{op}}
\def\dual{^{\scriptscriptstyle\natural}}

\def\ls#1#2{{\,^{#1}\!#2}}
\def\Iup{I^{\uparrow}}
\def\Idown{I_{\downarrow}}

\def\Z{\mathbb{Z}}
\def\N{\mathbb{N}}

\def\C{\mathbb{C}}

\def\F{\mathbb{F}}

\def\S{\mathbb{S}}

\newcommand{\sumb}[2]{\sum_{{\scriptstyle #1}\atop {\scriptstyle #2}}}

\def\bigoplusl{\bigoplus\limits}

\newcommand{\edge}[2]{\xymatrix{#1\ar@{->-}[r]&#2}}
\def\marc[#1]{\ar@{-}[#1]|(.4){\object@{<}}}
\def\mard[#1]{\ar@{-}[#1]|(.5){\object@{>}}}
\def\marb[#1]{\ar@{-}[#1]|{\object+{  }}}
\newcommand{\fleche}[2]{\xymatrix@C=4ex{*!U(0.2){#1\;}&*!U(0.5){\;#2}\marc[l]}}
\newcommand{\flecheb}[2]{\xymatrix@C=4ex{*!U(0.2){#1\;}&*!U(0.1){\;#2}\marc[l]}}

\def\pf{\par\bigskip\noindent{\bf Proof~: }}
\def\endpf{~\hfill\rlap{\hspace{-1ex}\raisebox{.5ex}{\framebox[1ex]{}}\sp}\bigskip\pagebreak[3]}

\makeatletter
\renewenvironment{enumerate}{\ifnum \@enumdepth >3 \@toodeep\else
       \advance\@enumdepth \@ne
       \edef\@enumctr{enum\romannumeral\the\@enumdepth}\list
       {\csname  label\@enumctr\endcsname}{\setlength{\topsep}{1ex}
 \setlength{\itemsep}{0 pt}\usecounter
         {\@enumctr}\def\makelabel##1{\hss\llap{##1}}}\fi}{\endlist}

\renewenvironment{itemize}{\ifnum \@itemdepth >3 \@toodeep\else
\advance\@itemdepth \@ne
\edef\@itemitem{labelitem\romannumeral\the\@itemdepth}
\list{\csname\@itemitem\endcsname}{\setlength{\topsep}{1ex}\setlength
{\itemsep}{0pt}\def\makelabel##1{\hss\llap{##1}}}\fi}
{\endlist}
\def\@seccntformat#1{\csname the#1\endcsname.\quad}

\def\section{\pagebreak[3]\setcounter{prop}{0}\setcounter{equation}{0}\@startsection{section}{1}{\z@}{4ex plus  6ex}{2ex}{\center\reset@font \large\bf}}

\makeatother

\def\theprop{\thesection.\arabic{prop}}

\renewenvironment{equation}{\refstepcounter{subsection}\refstepcounter
{prop}$$}{\leqno{\bf (\theprop)}$$}

\newenvironment{rem}[1]{\refstepcounter{subsection}\refstepcounter
{prop} \mpn{{\bf \thesection.\arabic{prop}.}\ \ \bf#1.}}{\smp}
\newenvironment{enonce}[1]{\pagebreak[3]\refstepcounter{prop}\mmpn
{{\bf  \thesection.\arabic{prop}.\ #1.}}\begin{it} }{\end{it}\smp}
\def\thesection{\arabic{section}}

\newcommand{\result}[1]{\begin{enonce}{#1}}
\newcommand{\fresult}{\end{enonce}}

\newcommand{\mbigvee}[1]{\mathop{\bigvee}_{#1}\limits}
\newcommand{\mbigveeb}[2]{\mathop{\bigvee}_{{\scriptstyle #1}\atop {\scriptstyle #2}}\limits}

\newcommand{\mbigsmash}[1]{\mathop{\bigwedge}_{#1}\limits}
\def\smash{\wedge}

\begin{document}

\title[Correspondence functors and lattices]
{Correspondence functors and lattices} 

\author{Serge Bouc}
\author{Jacques Th\'evenaz}
\date\today



\begin{abstract}
A correspondence functor is a functor from the category of finite sets and correspondences
to the category of $k$-modules, where $k$ is a commutative ring. 
A main tool for this study is the construction of a correspondence functor associated to any finite lattice~$T$.
We prove for instance that this functor is projective if and only if the lattice~$T$ is distributive.
Moreover, it has quotients which play a crucial role in the analysis of simple functors.
The special case of total orders yields some more specific and complete results.
\end{abstract}

\maketitle


\section{Introduction}
The present paper is the second in a series which develops the theory of correspondence functors,
namely functors from the category of finite sets and correspondences to the category of $k$-modules, where $k$ is a commutative ring.
In the first paper~\cite{BT2}, we showed that the category of finitely generated correspondence functors is artinian when $k$ is a field.
In representation theory, simple modules, or simple functors, are the most basic and important objects to understand.
We showed in~\cite{BT2} how to parametrize the simple correspondence functors $S_{E,R,V}$
by means of a finite set $E$, an order relation $R$ on~$E$, and a simple $k\Aut(E,R)$-module~$V$ (up to isomorphism).\par

The present paper establishes a connection between finite lattices and correspondence functors.
Associated to any finite lattice~$T$, we construct a correspondence functor $F_T$ (Section~\ref{Section-FT}).
This is the second indication of the importance of posets and lattices in our work and we describe the interplay between lattices and functors.
For instance, one of our first results asserts that the functor $F_T$ is projective if and only if the lattice~$T$ is distributive (Theorem~\ref{projective-distributive}).\par

The second main purpose of this paper is to introduce a {\em fundamental functor} $\S_{E,R}$ associated to any finite poset $(E,R)$. This is a precursor of each of the simple correspondence functors $S_{E,R,V}$ and it turns out that understanding $\S_{E,R}$ is the key for understanding those simple functors.
In particular, the fundamental functors $\S_{E,R}$ will play a crucial role for the determination of $\dim(S_{E,R,V}(X))$ for any finite set~$X$, which will appear in our next paper~\cite{BT3}.
Actually, the formula for this dimension involves a new invariant associated to lattices which will be introduced in~\cite{BT3} and which will give another important motivation for studying the link between finite lattices and correspondence functors.\par

The fundamental functors can be analyzed by using lattices. If $(E,R)$ is the subposet of irreducible elements in a finite lattice~$T$, then the functor $F_T$ has a fundamental functor as a quotient, which turns out to be $\S_{E,R\op}$ where $R\op$ denotes the opposite order relation (Theorem~\ref{surjection}).
The kernel of the morphism $F_T\to\S_{E,R\op}$ can be described by a system of linear equations.

We show that there is a duality between $F_T$ and $F_{T\op}$ over any commutative ring~$k$ (Theorem~\ref{F_T dual}).
Moreover, the fundamental functor $\S_{E,R}$ also appears as a subfunctor of~$F_{T\op}$ (Theorem~\ref{SERgamma}).
In Section~\ref{Section-embed-total}, some endomorphisms and idempotents of a lattice~$T$ are defined, associated with all possible quotients of~$T$ which are total orders. In Section~\ref{Section-total}, those idempotents are used to describe completely the functor $F_T$ in the special case where $T$ is totally ordered and they are also used to find all direct summands associated to total orders in a functor~$F_T$ corresponding to an arbitrary lattice~$T$.


\section{Correspondence functors} \label{Section-functors}

\bigskip
\noindent
In this introductory section, we recall the basic facts we need about correspondence functors (which also appear in~\cite{BT2}).
We denote by $\CC$ the category of finite sets and correspondences.
Its objects are the finite sets and the set $\CC(Y,X)$ of morphisms from $X$ to~$Y$ is the set of all correspondences from $X$ to~$Y$,
namely all subsets of $Y\times X$ (using a reverse notation which is convenient for left actions).
Given two correspondences $R\subseteq Z\times Y$ and $S\subseteq Y\times X$,
their composition $RS$ is defined by
$$RS:=\{ \, (z,x)\in Z\times X \,\mid\, \exists\;y\in Y \;\text{ such that } \; (z,y)\in R \,\text{ and }\, (y,x)\in S \,\} \mpoint$$
A correspondence from $X$ to~$X$ is also called a relation on~$X$.\par

Let $\Sigma_X$ be the symmetric group of all permutations of~$X$.
Associated with a permutation $\sigma\in\Sigma_X$, there is a relation on~$X$ which we write
$$\Delta_\sigma:=\{(\sigma(x),x)\in X\times X\mid x\in X\} \mpoint$$
In particular, $\Delta_X:=\Delta_{\Id}$ is the identity morphism of the object~$X$.
If $\sigma,\tau\in\Sigma_X$, then $\Delta_{\sigma\tau}=\Delta_\sigma\Delta_\tau$.
The symmetric group $\Sigma_X$ acts on relations by conjugation
and we write $\ls\sigma R=\Delta_\sigma R \Delta_{\sigma^{-1}}$.\par

For any commutative ring~$k$, we let $k\CC$ be the $k$-linearization of~$\CC$.
The objects are again the finite sets and $k\CC(Y,X)$ is the free $k$-module with basis $\CC(Y,X)$.
A correspondence functor is a $k$-linear functor from $k\CC$ to $k\text{-\!}\Mod$.
We let $\CF_k$ be the category of all correspondence functors (for some fixed commutative ring~$k$).
This category has the following feature~:

\result{Lemma} \label{vanishing-below} Let $E$ and $F$ be finite sets with $|E| \leq |F|$.
Let $M$ be a correspondence functor. If $M(F)=0$, then $M(E)=0$.
\fresult

\pf
Since $|E| \leq |F|$, there exists an injective map $i:E\hookrightarrow F$.
Let $i_*\subseteq F\times E$ denote the correspondence 
$$i_*=\big\{\big(i(e),e\big)\mid e\in E\big\}\mvirg$$
and $i^*\subseteq E\times F$ denote the correspondence 
$$i^*=\big\{\big(e,i(e)\big)\mid e\in E\big\}\mpoint$$
As $i$ is injective, one checks easily that $i^*i_*=\Delta_E$, that is, $i^*i_*=\Id_E$.
For any $m\in M(E)$, we have $m=i^*i_*\cdot m$. But $i_*\cdot m\in M(F)$, so $i_*\cdot m=0$.
Therefore $m=0$.
\endpf

We define a {\em minimal set} for a correspondence functor~$F$ to be a finite set~$X$ of minimal cardinality such that $F(X)\neq 0$.
For a nonzero functor, such a minimal set always exists and is unique up to bijection.

The first instances of correspondence functors are the representable functors $k\CC(-,E)$, where $E$ is a finite set, and the functors
$$L_{E,W}:=k\CC(-,E)\otimes_{k\CC(E,E)}W$$
where $W$ is a left $k\CC(E,E)$-module.
The proof of the following result is easy and is sketched in Lemma~2.3 of~\cite{BST} in the special case of biset functors for finite groups,
but it extends without change to representations of an arbitrary category.

\result{Lemma} \label{adjunction}
Let $\CF_k$ be the category of all correspondence functors and let $E$ be a finite set.
The functor
$$k\CC(E,E){-}\Mod \;\longrightarrow\; \CF_k  \;,\quad W\mapsto L_{E,W}$$
is left adjoint of the evaluation functor
$$\CF_k  \;\longrightarrow\; k\CC(E,E){-}\Mod \;,\quad F\mapsto F(E) \mpoint$$
\fresult

The functor $L_{E,W}$ has a subfunctor $J_{E,W}$ defined on any finite set~$X$ by
$$J_{E,W}(X):=\Big\{\sum_i\phi_i\otimes w_i \in L_{E,W}(X) \mid
\forall\psi\in k\CC(E,X), \sum_i (\psi\phi_i)\cdot w_i=0 \Big\} \mpoint$$

\result{Lemma} \label{JEW}
Let $E$ be a finite set and let $W$ be a $k\CC(E,E)$-module.
\begin{enumerate}
\item $J_{E,W}$ is the unique subfunctor of $L_{E,W}$ which is maximal
with respect to the condition that it vanishes at~$E$.
\item If $W$ is a simple $k\CC(E,E)$-module, then $J_{E,W}$ is the unique maximal subfunctor of $L_{E,W}$
and $L_{E,W}/J_{E,W}$ is a simple functor.
\end{enumerate}
\fresult

\pf
The result is a slight extension of the first lemma of~\cite{Bo1}.
The proof is also sketched in Lemma~2.3 of~\cite{BST} in the special case of biset functors for finite groups,
but it extends without change to the representation theory of an arbitrary category.
\endpf

Now we want to consider the functor $L_{E,W}/J_{E,W}$ for some specific choices of $k\CC(E,E)$-modules.
The algebra $k\CC(E,E)$ of all relations on~$E$ was studied in~\cite{BT1} and we need a few facts from this approach.
A relation~$R$ on~$E$ is called {\em essential\/} if it does not factor through a set of cardinality strictly smaller than~$|E|$.
The $k$-submodule generated by the set of inessential relations is a two-sided ideal
$$I_E:=\sum_{|Y|<|E|} k\CC(E,Y)k\CC(Y,E)$$
and the quotient
$$\CE_E:=k\CC(E,E)/I_E$$
is called the {\em essential algebra}.
A large part of its structure has been elucidated in~\cite{BT1}.
There is a quotient algebra $\CP_E=\CE_E/N$, where $N$ is a nilpotent two-sided ideal defined in~\cite{BT1}.
We call $\CP_E$ the {\em algebra of permuted orders},
because it has a $k$-basis consisting of all relations on~$E$ of the form
$\Delta_\sigma R$, where $\sigma$ runs through the symmetric group $\Sigma_E$ of all permutations of~$E$,
and $R$ is an order on~$E$. By an order, we always mean a partial order relation.
The product of two orders $R$ and $S$ in~$\CP_E$ is the transitive closure of $R\cup S$ if this closure is an order, and zero otherwise.
This describes completely the algebra structure of~$\CP_E$.\par

Among the $k\CC(E,E)$-modules, there is the {\em fundamental module} $\CP_E f_R$,
associated to any poset $(E,R)$, where $E$ is a finite set and $R$ denotes the order relation on~$E$ which defines the poset structure.
Here $f_R$ is a suitable idempotent in~$\CP_E$, depending on~$R$, and $\CP_E f_R$ is the left ideal generated by~$f_R$.
Actually, the fundamental module $\CP_E f_R$ only depends on the isomorphism type of the poset $(E,R)$, or in other words, for a fixed set~$E$, on the $\Sigma_E$-conjugacy class of~$R$. More explicitly,  if $\sigma\in\Sigma_E$,
then conjugation by~$\sigma$ induces an isomorphism of posets $(E,R)\cong (E,\ls\sigma R)$ and we also have an isomorphism of $\CP_E$-modules 
$\CP_E f_R\cong \CP_E f_{\ls\sigma R}$, because $f_{\ls\sigma R}=\ls\sigma f_R$ (see Lemma~7.1 in~\cite{BT1} for details).\par

The only thing we really need to know about the fundamental module $\CP_E f_R$ is its structure as a $k\CC(E,E)$-module.
This is described in the next result, which combines Corollary 7.3 and Proposition~8.5 of~\cite{BT1} (see also Proposition~4.5 of~\cite{BT2} and its use).

\result{Proposition} \label{fundamental-module}
Let $E$ be a finite set and $R$ an order on~$E$.
\begin{enumerate}
\item The fundamental module $\CP_E f_R$ is a left module for the algebra $\CP_E$, hence also a left module for the essential algebra~$\CE_E$ and for the algebra of relations~$k\CC(E,E)$.
\item $\CP_E f_R$ is a free $k$-module with a $k$-basis consisting of the elements $\Delta_\sigma f_R$,
where $\sigma$ runs through the group $\Sigma_E$ of all permutations of $E$.
\item $\CP_E f_R$ is a $(\CP_E,k\Aut(E,R))$-bimodule and the right action of $k\Aut(E,R)$ is free.
Explicitly, the right action of $\tau\in\Aut(E,R)$ maps the basis element $\Delta_\sigma f_R$ to the basis element $\Delta_{\sigma\tau} f_R$.
\item The action of the algebra of relations $k\CC(E,E)$ on the module $\CP_E f_R$
is given as follows. For any relation $Q\in \CC(E,E)$,
$$Q\cdot\Delta_\sigma f_R=\left\{\begin{array}{ll}
\Delta_{\tau\sigma}f_R&\hbox{if}\;\;\exists\tau\in\Sigma_E\;\hbox{such that}\;
\Delta_E\subseteq \Delta_{\tau^{-1}}Q\subseteq {\ls\sigma R},\\
0&\hbox{otherwise}\mvirg\end{array}\right.$$
where $\ls\sigma R=\big\{\big(\sigma(e),\sigma(f)\big)\mid (e,f)\in R\big\}$ (or equivalently $\ls\sigma R=\Delta_\sigma R\Delta_{\sigma^{-1}}$).
(Moreover, $\tau$ is unique in the first case.)
\end{enumerate}
\fresult

Using the $(\CP_E,k\Aut(E,R))$-bimodule structure on~$\CP_E f_R$, we define
$$T_{R,V}:=\CP_E f_R\otimes_{k\Aut(E,R)} V \mvirg$$
where $V$ is any $k\Aut(E,R)$-module.
Then $T_{R,V}$ is a left $\CP_E$-module, hence also a $k\CC(E,E)$-module since $\CP_E$ is a quotient of~$k\CC(E,E)$.
The left action of $k\CC(E,E)$ on~$T_{R,V}$ is induced from the action on~$\CP_Ef_R$ described in Proposition~\ref{fundamental-module} above.
Again, the module $T_{R,V}$ is invariant under $\Sigma_E$-conjugacy, that is, for a fixed set~$E$ and for $\sigma\in \Sigma_E$, we have an isomorphism of  $\CP_E$-modules $T_{\ls\sigma R,\ls\sigma V}\cong T_{R,V}$, where $\ls\sigma V$ denotes the conjugate module, namely a module for the group $\Aut(E,\ls\sigma R)=\sigma\Aut(E,R)\sigma^{-1}$ (see Theorem~8.1 in~\cite{BT1} for details).\par

The main thing we need to know about $T_{R,V}$ is the following result, which is part of Theorem~8.1 in~\cite{BT1}.

\result{Proposition} \label{simple-modules} Assume that $k$ is a field. Let $E$ be a finite set, $R$ an order on~$E$, and $V$ a simple $k\Aut(E,R)$-module.
Then $T_{R,V}$ is a simple $\CP_E$-module (hence also a simple $\CE_E$-module).
\fresult

Actually, Theorem~8.1 in~\cite{BT1} asserts that every simple $\CE_E$-module is isomorphic to some module $T_{R,V}$ and that, consequently,
the set of isomorphism classes of simple $\CE_E$-modules is parametrized by the set of conjugacy classes of pairs $(R,V)$ where $R$ is an order on~$E$ and $V$ is a simple $k\Aut(E,R)$-module.\par

Associated with the above $k\CC(E,E)$-modules, we can now define some specific correspondence functors.
Using the fundamental module $\CP_E f_R$, we define
$$\S_{E,R}:=L_{E,\CP_E f_R}/J_{E,\CP_E f_R}$$
and we call it the {\em fundamental functor} associated with the poset $(E,R)$.
Using the module $T_{R,V}$, we define
$$S_{E,R,V}:=L_{E,T_{R,V}}/J_{E,T_{R,V}} \mpoint$$
Note that $\S_{E,\ls\sigma R}\cong \S_{E,R}$ and $S_{E,\ls\sigma R,\ls\sigma V}\cong S_{E,R,V}$, for any permutation $\sigma\in\Sigma_E$.

\result{Proposition} \label{SER-SERV}
\begin{enumerate}
\item The set $E$ is a minimal set for $\S_{E,R}$ and $\S_{E,R}(E)\cong \CP_E f_R$ as left $k\CC(E,E)$-modules.
\item The set $E$ is a minimal set for $S_{E,R,V}$ and $S_{E,R,V}(E)\cong T_{R,V}$ as left $k\CC(E,E)$-modules.
\item If $k$ is a field and $V$ is a simple $k\Aut(E,R)$-module, then $S_{E,R,V}$ is a simple correspondence functor.
\end{enumerate}
\fresult

\pf
Let $W$ be any $\CP_E$-module and set $S=L_{E,W}/J_{E,W}$.
Suppose that $Y$ is a finite set such that $S(Y)\neq\{0\}$. Then $L_{E,W}(Y)\neq J_{E,W}(Y)$,
so there exists a correspondence $\phi\in \CC(Y,E)$ and $v\in W$ such that $\phi\otimes v\in L_{E,W}(Y)-J_{E,W}(Y)$.
By definition of~$J_{E,W}$, this means that there exists a correspondence $\psi\in \CC(E,Y)$
such that $\psi\phi\cdot v\neq 0$. Since $\CP_E$ is a quotient of the essential algebra~$\CE_E$,
it follows that $W$ is a module for $\CE_E=k\CC(E,E)/I_E$, so that the ideal $I_E$ acts by zero on~$W$.
Therefore $\psi\phi\notin I_E$. But $\psi\phi$ factorizes through~$Y$, so we must have $|Y|\geq |E|$.
Thus $E$ is a minimal set for~$S$.
In particular, this holds for $S=\S_{E,R}$ (taking $W=\CP_Ef_R$) and also for $S=S_{E,R,V}$ (taking $W=T_{R,V}$).\par

Since $J_{E,W}$ vanishes at~$E$ by Lemma~\ref{JEW}, evaluation at~$E$ yields
$$S(E)=L_{E,W}(E)/J_{E,W}(E)=L_{E,W}(E)=k\CC(E,E)\otimes_{k\CC(E,E)} W\cong W \mvirg$$
and therefore $\S_{E,R}(E)\cong \CP_E f_R$ and $S_{E,R,V}(E)\cong T_{R,V}$.\par

For the proof of (c), notice that $T_{R,V}$ is a simple $k\CC(E,E)$-module by Proposition~\ref{simple-modules} and therefore
$$S_{E,R,V}=L_{E,T_{R,V}}/J_{E,T_{R,V}}$$
is a simple correspondence functor by Lemma~\ref{JEW}.
\endpf

Although we do not need it here, let us mention that more is known about simple correspondence functors, when $k$ is a field.
It is proved in~\cite{BT2} that any simple functor has the form $S_{E,R,V}$ for some triple $(E,R,V)$ and that
the set of isomorphism classes of simple correspondence functors is parametrized by
the set of isomorphism classes of triples $(E,R,V)$ where $E$ is a finite set, $R$ is an order on~$E$, and $V$ is a simple $k\Aut(E,R)$-module.\par

We note that the fundamental functor $\S_{E,R}$ is a precursor of~$S_{E,R,V}$, in the sense of the following lemma.

\result{Lemma} \label{precursor}
Suppose that $V$ is a $k\Aut(E,R)$-module generated by a single element~$v$ (e.g. a simple module).
Consider the surjective morphism of correspondence functors
$$\Phi: L_{E,\CP_Ef_R} \longrightarrow L_{E,T_{R,V}}$$
induced by the surjective homomorphism of $\CP_E$-modules
$$\Phi_E:\CP_Ef_R \longrightarrow \CP_Ef_R \otimes_{k\Aut(E,R)} V =T_{R,V} \,, \qquad a\mapsto a\otimes v \mpoint$$
Then $\Phi$ induces a surjective morphism of correspondence functors
$$\S_{E,R} \longrightarrow S_{E,R,V} \mpoint$$
\fresult

\pf
Let $\varphi\otimes f_R \in k\CC(X,E)\otimes_{k\CC(E,E)} \CP_Ef_R = L_{E,\CP_Ef_R}(X)$.
If $\varphi\otimes f_R \in J_{E,\CP_Ef_R}(X)$, then for every $\psi\in\CC(E,X)$, we have $\psi \varphi\cdot f_R=0$.
Then $\Phi_X(\varphi\otimes f_R)=\varphi \otimes f_R\otimes v$ and we have
$$\psi \varphi \cdot (f_R\otimes v) = (\psi \varphi \cdot f_R)\otimes v=0 \mpoint$$
This shows that $\Phi_X(\varphi\otimes f_R) \in J_{E,T_{R,V}}(X)$, so that
$\Phi_X(J_{E,\CP_Ef_R}(X)) \subseteq J_{E,T_{R,V}}(X)$.
Therefore $\Phi$ induces a morphism of correspondence functors
$\S_{E,R} \longrightarrow S_{E,R,V}$ which remains surjective.
\endpf

By means of a very detailed analysis of $\S_{E,R}$ which will be carried out in~\cite{BT3}, we shall show that it is possible to recover $S_{E,R,V}$ from~$\S_{E,R}$ by simply tensoring with~$V$. Consequently, mastering $\S_{E,R}$ will be the key for obtaining information about the simple functors $S_{E,R,V}$ and this explains why the fundamental functors play a crucial role throughout our work.


\section{Posets and lattices} \label{Section-posets}

\bigskip
\noindent
In this section, we give some definitions, fix some notation, and prove some basic lemmas, which will be used throughout.

By an {\em order} $R$ on a finite set~$E$, we mean a partial order relation on~$E$.
In other words, $(E,R)$ is a finite poset.

\result{Notation and definitions} Let $(E,R)$ be a finite poset.
\begin{enumerate}
\item We write $\leq_R$ for the order relation, so that $(a,b)\in R$ if and only if $a\leq_R b$.
Moreover $a<_R b$ means that $a\leq_R b$ and $a \neq b$.
\item If $a,b\in E$ with $a\leq_R b$, we define intervals
$$\begin{array}{ll}
\qquad[a,b]_E:=\{x\in E\mid a\leq_R x\leq_R b\} \,,\qquad & {]a,b[_E}:=\{x\in E\mid a<_R x<_R b\} \,, \\
\qquad{[a,b[_E}:=\{x\in E\mid a\leq_R x<_R b\} \,,\qquad & {]a,b]_E}:=\{x\in E\mid A<_R x\leq_R b\} \,, \\
\qquad{[a,\cdot[_E}:=\{x\in E\mid a\leq_R x\} \,,\qquad & {]\cdot,b]_E}:=\{x\in E\mid x\leq_R b\} \mpoint
\end{array}$$
When the context is clear, we write ${[a,b]}$ instead of~${[a,b]_E}$.
\item A subset~$A$ of~$E$ is a {\em lower $R$-ideal}, or simply a {\em lower ideal},
if, whenever $a\in A$ and $x\leq_R a$, we have $x\in A$.
Similarly, a subset~$A$ of~$E$ is an {\em upper $R$-ideal}, or simply an {\em upper ideal},
if, whenever $a\in A$ and $a\leq_R x$, we have $x\in A$.
\item A {\em principal lower ideal}, or simply {\em principal ideal}, is a subset of the form ${]\cdot,a]_E}$, where $a\in E$.
A {\em principal upper ideal} is defined similarly.
\item The {\em opposite order relation} $R\op$ is defined by the property that $a\leq_{R\op}b$ if and only if $b\leq_R a$.
\end{enumerate}
\fresult

\result{Notation and definitions} Let $T$ be a finite lattice.
\begin{enumerate}
\item We write $\leq_T$, or sometimes simply $\leq$, for the order relation, $\vee$ for the join (least upper bound), $\wedge$ for the meet (greatest lower bound), $\hat{0}$ for the least element and $\hat{1}$ for the greatest element.
\item An element $e\in T$ is called {\em join-irreducible}, or simply {\em irreducible},
if, whenever $e=\mbigvee{a\in A} a$ for some subset $A$ of~$T$, then $e\in A$.
In case $A=\emptyset$, the join is~$\hat0$ and it follows that $\hat0$ is not irreducible.
An element $e\neq\hat0$ is irreducible if and only if the equality $e=s\vee t$, for $s,t\in T$, implies $e=s$ or $e=t$.
In other words, if $e\neq\hat0$, then $e$ is irreducible if and only if $[\hat0,e[$ has a unique maximal element.
\item Let $(E,R)$ be a subposet of~$T$. We say that it is a {\em full\/} subposet of~$T$ if for all $e,f\in E$ we have~:
$$e\leq_R f \Longleftrightarrow e\leq_T f \mpoint$$   

\end{enumerate}
\fresult

Note that if $(E,R)$ is the poset of irreducible elements in a finite lattice $T$, then $T$ is {\em generated\/} by~$E$ in the sense that any element $x\in T$ is a join of elements of~$E$.
To see this, define the height of $t\in T$ to be the maximal length of a chain in $[\hat0,t]_T$.
If $x$ is not irreducible and $x\neq\hat0$, then $x=t_1\vee t_2$ with $t_1$ and $t_2$ of smaller height than~$x$.
By induction on the height, both $t_1$ and $t_2$ are joins of elements of~$E$.
Therefore $x=t_1\vee t_2$ is also a join of elements of~$E$.\par

\result{Notation} \label{up-down}
Let $(E,R)$ be a finite poset.
\begin{enumerate}
\item Let $\Idown(E,R)$ denote the set of lower $R$-ideals of $E$.
Then $\Idown(E,R)$, ordered by inclusion of subsets, is a lattice~:
the join operation is union of subsets, and the meet operation is intersection.\par
\item Similarly, $\Iup(E,R)$ denotes the set of upper $R$-ideals of $E$, which is also a lattice.
If $R\op$ is the relation opposite to~$R$, then clearly $\Iup(E,R)=\Idown(E,R\op)$.
\end{enumerate}
\fresult

\begin{rem}{Remark} \label{preorder}
Let $R$ be a preorder on a finite set~$E$, that is, a relation which is reflexive and transitive.
There is an equivalence relation~$\sim$ associated with~$R$, defined by
$$x\sim y \;\Longleftrightarrow\;  (x,y)\in R \;\text{ and }\; (y,x)\in R \mpoint$$
Then $R$ induces an order relation~$\overline R$ on the quotient set $\overline E=E/{\sim}\,$ such that
$$(x,y)\in R \;\Longleftrightarrow\;  (\overline x,\overline y)\in \overline R \mvirg$$
where $\overline x$ denotes the equivalence class of~$x$ under~$\sim$.
It is easy to see that the quotient map $E\to \overline E$ induces an isomorphism of lattices
$\Idown(E,R) \cong \Idown(\overline E,\overline R)$.\par

Note that it is proved in Lemma~3.9 of~\cite{BT2} that the representable functors $k\CC(-,E)R$ and $k\CC(-,\overline E)\overline R$
are isomorphic, but actually we will view them in a new way in Proposition~\ref{C(-,E)R}.
These remarks show that, for our purposes, it is enough to consider orders rather than preorders,
and we shall do so in the rest of this paper, without loss of generality.
\end{rem}

\result{Lemma} \label{principal} Let $(E,R)$ be a finite poset.
\begin{enumerate}
\item The irreducible elements in the lattice~$\Idown(E,R)$ are the principal ideals $]\cdot,e]_E$, where $e\in E$.
Thus the poset $E$ is isomorphic to the poset of all irreducible elements in~$\Idown(E,R)$
by mapping $e\in E$ to the principal ideal $]\cdot,e]_E$.
\item $\Idown(E,R)$ is a distributive lattice.
\item If $T$ is a distributive lattice and $(E,R)$ is its subposet of irreducible elements, then $T$ is isomorphic to $\Idown(E,R)$.
\item For any finite lattice $T$ having $(E,R)$ as poset of irreducible elements, there is a join-preserving surjective map
$f: \Idown(E,R) \longrightarrow T$ which sends any lower ideal $A\in \Idown(E,R)$
to the join $\bigvee_{e\in A}e$ in~$T$.
\end{enumerate}
\fresult

\pf This is not difficult and well-known.
For details, see Theorem~3.4.1 and Proposition~3.4.2 in~\cite{St}.
\endpf

\result{Convention} \label{convention}
In the situation of Lemma~\ref{principal}, we shall identify $E$ with its image
via the map
$$E\longrightarrow\Idown(E,R) \,,\qquad  e \mapsto \;]\cdot,e]_E \mpoint$$
Thus we view $E$ as a full subposet of~$\Idown(E,R)$.
This abusive convention is a conceptual simplification and has many advantages for the rest of this paper.
\fresult

Given a poset $(E,R)$, the map
$$E\longrightarrow\Iup(E,R) \,,\qquad  e \mapsto \;[e,\cdot[_E$$
is order-reversing, so it is in fact $(E,R\op)$ which is identified with the poset of irreducible elements in~$\Iup(E,R)$.
Since $\Iup(E,R)=\Idown(E,R\op)$, this is actually just Convention~\ref{convention} applied to~$R\op$.\par

We now introduce a notation which will play an important role in our work
(and which was already used in the proof of Theorem~9.2 in~\cite{BT2}).

\result{Notation} \label{Gamma}
Let $T$ be a finite lattice and let $(E,R)$ be the full subposet of its irreducible elements.
For any finite set $X$ and any map $\varphi: X\to T$, we associate the correspondence
$$\Gamma_\varphi:= \{ (x,e)\in X\times E \mid e\leq_T \varphi(x) \} \subseteq X\times E \mpoint$$
In the special case where $T=\Idown(E,R)$ and in view of Convention~\ref{convention}, we obtain
$$\Gamma_\varphi= \{ (x,e)\in X\times E \mid e\in \varphi(x) \} \mpoint$$
\fresult

\result{Lemma} \label{rightR}
Let $T$ be a finite lattice and let $(E,R)$ be the full subposet of its irreducible elements.
\begin{enumerate}
\item For any map $\varphi: X\to T$, we have $\Gamma_\varphi R\op =\Gamma_\varphi$.
\item If $T=\Idown(E,R)$, then a correspondence $S\subseteq X\times E$ has the form $S=\Gamma_\varphi$
for some map $\varphi: X\to \Idown(E,R)$ if and only if $SR\op=S$.
\item If $T=\Iup(E,R)$, then a correspondence $S\subseteq X\times E$ has the form $S=\Gamma_\varphi$
for some map $\varphi: X\to \Iup(E,R)$ if and only if $SR=S$.
\end{enumerate}
\fresult

\pf
(a) Since $\Delta_E\subseteq R\op$, we always have
$\Gamma_\varphi=\Gamma_\varphi\Delta_E\subseteq\Gamma_\varphi R\op$.
Conversely, if $(x,f)\in \Gamma_\varphi R\op$, then there exists $e\in E$ such that
$(x,e)\in \Gamma_\varphi$ and $(e,f)\in R\op$, that is, $e\leq_T \varphi(x)$ and $f\leq_R e$.
But $f\leq_R e$ if and only if $f\leq_T e$, because $(E,R)$ is a full subposet of~$T$.
It follows that $f\leq_T \varphi(x)$, that is, $(x,f)\in \Gamma_\varphi$.
Thus $\Gamma_\varphi R\op\subseteq \Gamma_\varphi$ and equality follows.\mpn

(b) One direction follows from (a). For the other direction, let $S\in\CC(X,E)$ be such that $SR\op=S$,
or equivalently $S\in \CC(X,E)R\op$ (because $R\op$ is idempotent by reflexivity and transitivity).
Then the set
$$\phi(x)=\big\{e\in E\mid (x,e)\in S\big\}$$
is a lower $R$-ideal in $E$,
thus $\phi$ is a function $X\to \Idown(E,R)$. Clearly $\Gamma_\phi=S$.\mpn

(c) This follows from (b) applied to~$R\op$, because $\Iup(E,R)=\Idown(E,R\op)$
and $(E,R\op)$ is its poset of irreducible elements.
\endpf


\section{Functors associated to lattices} \label{Section-FT}

\bigskip
\noindent
A fundamental construction associates a correspondence functor $F_T$ to any finite lattice~$T$. This is one of our main tools for the analysis of correspondence functors. Throughout this section, $k$ is an arbitrary commutative ring.

\result{Definition} \label{FT}
Let $T$ be a finite lattice.
For a finite set $X$, we define $F_T(X)$ to be the free $k$-module with basis the set $T^X$ of all functions from $X$ to~$T$~:
$$F_T(X):=k(T^X)\mpoint$$
For two finite sets $X$ and $Y$ and a correspondence $R\subseteq Y\times X$, we define a map $F_T(R):F_T(X)\to F_T(Y)$ as follows~:
to a function $\varphi:X\to T$, we associate the function $F_T(R)(\varphi):Y\to T$, also simply denoted by $R\varphi$, defined by
$$(R\varphi)(y):=\mbigvee{\substack{x\in X \\ (y,x)\in R}} \varphi(x)\mvirg$$
with the usual rule that a join over the empty set is equal to $\hat{0}$. The map
$$F_T(R):F_T(X)\to F_T(Y)$$
is the unique $k$-linear extension of this construction.
More generally, for every element $\alpha=\sum_{R\in\mathcal{C}(Y,X)}\limits\alpha_RR \in k\mathcal{C}(Y,X)$,
where $\alpha_R\in k$, we set
$$F_T(\alpha)=\sum_{R\in\mathcal{C}(Y,X)}\limits\alpha_RF_T(R) \mpoint$$
\fresult

\result{Proposition} The assignment sending a finite set $X$ to the $k$-module $F_T(X)$ and a morphism $\alpha\in k\mathcal{C}(Y,X)$
to the $k$-linear map $F_T(\alpha):F_T(X)\to F_T(Y)$ is a correspondence functor.
\fresult

\pf First it is clear that if $X$ is a finite set and $\Delta_X\in\mathcal{C}(X,X)$ is the identity correspondence, then for any $\varphi:X\to T$ and any $y\in X$
$$(\Delta_X\varphi)(y)=\mbigvee{(y,x)\in \Delta_X}\varphi(x)=\varphi(y)\mvirg$$
hence $\Delta_X\varphi=\varphi$ and $F_T(\Delta_X)$ is the identity map of~$F_T(X)$.\par
Now if $X$, $Y$, and $Z$ are finite sets, if $R\in \mathcal{C}(Y,X)$ and $S\in\mathcal{C}(Z,Y)$, then for any $\varphi:X\to T$ and any $z\in Z$, we have
\begin{eqnarray*}
\big(S(R\varphi)\big)(z)&=&\mbigvee{(z,y)\in S}(R\varphi)(y)\\
&=&\mbigvee{(z,y)\in S}\mbigvee{(y,x)\in R}\varphi(x)\\
&=&\mbigvee{(z,x)\in S R}\varphi(x)\\
&=&(SR\varphi)(z)\mpoint
\end{eqnarray*}
By linearity, it follows that $F_T(\beta)\circ F_T(\alpha)=F_T(\beta\,\alpha)$, for any $\beta\in k\mathcal{C}(Z,Y)$ and any $\alpha\in k\mathcal{C}(Y,X)$.
\endpf

\begin{rem}{Remark} \label{join-semi-lattice}
The definition of~$F_T$ only uses the join operation in the lattice~$T$. It follows that the definition would work for a join semi-lattice, but it is actually well-known that a finite join semi-lattice has automatically a structure of lattice (the meet operation being uniquely determined from the sole join). This explains our choice of working with lattices. Such a choice will also be useful in Section~\ref{Section-duality-lattices} when we shall work with opposite lattices.
\end{rem} 

We now establish the link between the action of correspondences on functions $\varphi:X\to T$
(as in Definition~\ref{FT} above)
and the correspondences $\Gamma_\varphi$ defined in Notation~\ref{Gamma}.

\result{Lemma} \label{SGamma}
Let $T$ be a finite lattice and let $(E,R)$ be the full subposet of its irreducible elements.
Assume that $T$ is distributive, so that $T\cong \Idown(E,R)$.
Then, for any finite sets $X$, $Y$, any correspondence $S\in \CC(Y,X)$, and any function $\varphi:X\to T$, we have
$$\Gamma_{S\varphi\,}=S\Gamma_{\varphi} \mvirg$$
where $\Gamma_{\varphi}$ is defined in Notation~\ref{Gamma}.
\fresult

\pf
Let $y\in Y$ and $e\in E$. Then
$$\begin{array}{rl}
(y,e)\in \Gamma_{S\varphi} &\iff e\leq_T(S\varphi)(y) \iff e\leq_T \mbigvee{(y,x)\in S} \varphi(x) \\
&\iff e=e\,\smash\, \big(\mbigvee{(y,x)\in S} \varphi(x)\big) \mpoint
\end{array}$$
But, since $T$ is distributive, the latter equality is equivalent to $\;e=\mbigvee{(y,x)\in S} \big(e\smash \varphi(x)\big)$.
Now, since $e$ is irreducible, this is in turn equivalent to
$$\begin{array}{rl}
&\exists x\in X \,, \; (y,x)\in S \,\text{ and } \,e\smash \varphi(x)=e
\iff \exists x\in X \,, \;(y,x)\in S \,\text{ and } e\leq_T \varphi(x) \\
&\iff \exists x\in X \,, \;(y,x)\in S \,\text{ and } (x,e)\in \Gamma_{\varphi}
\iff (y,e)\in S \Gamma_{\varphi} \mpoint
\end{array}$$
This completes the proof.
\endpf

Now we can give another description of the correspondence functor associated to a distributive lattice.

\result{Proposition} \label{C(-,E)R} Let $(E,R)$ be a finite poset.
\begin{enumerate}
\item For any finite set $X$
$$\big\{\Gamma_\varphi \mid \varphi:X\to \Iup(E,R)\big\}=\big\{S\in \CC(X,E)\mid S R=S\big\}=\CC(X,E)R\mpoint$$
\item The correspondence functor $F_{\Iup(E,R)}$ is isomorphic to $k\CC(-,E)R$.
In particular $F_{\Iup(E,R)}$ is a projective object of~$\CF_k$.
\item The correspondence functor $F_{\Idown(E,R)}$ is isomorphic to $k\CC(-,E)R\op$.
In particular $F_{\Idown(E,R)}$ is a projective object of~$\CF_k$.
\end{enumerate}
\fresult

\pf (a) This is a restatement of Lemma~\ref{rightR}.\mpn

(b) The map
$$F_{\Iup(E,R)}(X) \longrightarrow k\CC(X,E)R \,,\qquad
\varphi\mapsto \Gamma_{\varphi}$$
is an isomorphism of correspondence functors, by (a) and Lemma~\ref{SGamma}.
Moreover $k\CC(-,E)$ is a projective functor by Yoneda's lemma and $k\CC(-,E)R$ is a direct summand of $k\CC(-,E)$ because $R$ is idempotent.
Therefore $k\CC(-,E)R$ is projective.\mpn

(c) follows from (b) and the obvious equality $\Idown(E,R)=\Iup(E,R\op)$.
\endpf

We now introduce a suitable category $\CL$ of lattices, as well as its $k$-linearization~$k\CL$.
Our aim is to show that the assignment $T\mapsto F_T$ becomes a $k$-linear functor from $k\CL$ to $\mathcal{F}_k$,
which will have the remarkable property of being full and faithful.

\result{Definition}\label{category L} Let $\CL$ and $k\CL$ denote the following categories~:
\begin{itemize}
\item The objects of $\CL$ and $k\CL$ are the finite lattices.
\item For any two lattices $T$ and $T'$, the set $\Hom_\CL(T,T')$ is the set 
of all maps ${f:T\to T'}$ which commute with joins, i.e. such that
$$f(\mbigvee{a\in A}a)=\mbigvee{a\in A}f(a)\mvirg$$
for any subset $A$ of~$T$.
\item For any two lattices $T$ and $T'$, the set $\Hom_{k\CL}(T,T')$ is the free $k$-module $k\Hom_\CL(T,T')$ with basis $\Hom_\CL(T,T')$.
\item The composition of morphisms in $\CL$ is the composition of maps.
\item The composition of morphisms in $k\CL$ is the $k$-bilinear extension of the composition in~$\CL$.
\end{itemize}
\fresult

It is easy to see that a morphism in $\CL$ is order-preserving, by considering the join $t_1\vee t_2$ in the case where $t_1\leq_T t_2$ in the lattice~$T$.

\begin{rem}{Remark} \label{not-lattice-morphism} The case $A=\emptyset$ in Definition~\ref{category L} shows that a morphism $f:T\to T'$ in $\CL$ always maps $\hat{0}\in T$ to $\hat{0}\in T'$.
Conversely, if $f:T\to T'$ satisfies $f(\hat{0})=\hat{0}$ and $f(a\vee b)=f(a)\vee f(b)$ for all $a,b\in T$,  then $f$ is a morphism in~$\CL$.\par

Morphisms in $\CL$ are morphisms of join semi-lattices (see Remark~\ref{join-semi-lattice}),
but they are generally not morphisms of lattices in the sense that they need not commute with the meet operation.
The choice of not using the meet operation turns out to be important for the next main theorem.
\end{rem}

\bigskip
For our next theorem, we need some notation. Let $f:T\to T'$ be a morphism in the category $\CL$.
For a finite set $X$, let $F_{f,X}:F_T(X)\to F_{T'}(X)$ be the $k$-linear map sending the function $\varphi:X\to T$ to the function $f\circ\varphi:X\to T'$. 

\result{Theorem} \label{fully faithful}
\begin{enumerate}
\item Let $f:T\to T'$ be a morphism in the category $\CL$. Then the collection of maps $F_{f,X}:F_T(X)\to F_{T'}(X)$, for all finite sets $X$, yields a natural transformation $F_f:F_T\to F_{T'}$ of correspondence functors.
\item The assignment sending a lattice $T$ to $F_T$, and a morphism $f:T\to T'$ in~$\CL$ to $F_f:F_T\to F_{T'}$, yields a functor $\CL\to \CF_k$. This functor extends uniquely to a $k$-linear functor
$$F_?:k\CL\longrightarrow \CF_k \mpoint$$
\item The functor $F_?$ is fully faithful.
\end{enumerate}
\fresult

\pf (a) Let $X$ and $Y$ be finite sets, let $\varphi: X\to T$ be a function, and let $U\in\CC(Y,X)$ be a correspondence.
Then
$F_{T'}(U) \big( F_{f,X}(\varphi) \big) =F_{T'}(U)(f\circ\varphi)=U(f\circ\varphi)$
and $F_{f,Y} \big( F_T(U)(\varphi)\big) = F_{f,Y} (U\varphi)=f\circ U\varphi$.
We show that they are equal by evaluating at any $y\in Y$~:
\begin{eqnarray*}
U(f\circ\varphi)(y)&=&\mbigvee{(y,x)\in U}(f\circ\varphi)(x)\\
&=&\mbigvee{(y,x)\in U}f\big(\varphi(x)\big)\\
&=&f\big(\mbigvee{(y,x)\in U}\varphi(x)\big)\\
&=&(f\circ U\varphi)(y)\mvirg
\end{eqnarray*}
hence $U(f\circ\varphi)=f\circ U\varphi$, which proves (a).\mpn

(b) It follows that the assignment $T\mapsto F_T$ is a functor $\CL\to\CF_k$. Since $k\CL$ is the $k$-linearization of $\CL$, this functor extends uniquely to a $k$-linear functor $F_?:k\CL\to\CF_k$.\mpn

(c) Let $S$ and $T$ be finite lattices, and $\Phi:F_{S}\to F_{T}$ be a morphism of functors. Thus, for any finite set $X$, we have a morphism of $k$-modules $\Phi_X:F_{S}(X)\to F_{T}(X)$ such that for any finite set $Y$ and any correspondence $R\subseteq (Y\times X)$, the diagram
$$\xymatrix{
F_{S}(X)\ar[r]^-{\Phi_X}\ar[d]_-{F_{S}(R)}&F_{T}(X)\ar[d]^-{F_{T}(R)}\\
F_{S}(Y)\ar[r]^-{\Phi_Y}&F_{T}(Y)
}
$$
is commutative. In other words, for any function $\alpha:X\to S$
\begin{equation}\label{fonctoriel}
R\Phi_X(\alpha)=\Phi_Y(R\alpha)\mpoint
\end{equation}
Taking $X=S$ and $\alpha=\Id_{S}$ in this relation, and setting
$$\varphi=\Phi_{S}(\Id_{S})=\sum_{\lambda:{S}\to{T}}u_\lambda\lambda\mvirg$$
where $u_\lambda\in k$, this gives
$$R\varphi=\Phi_Y(R\Id_{S})\mvirg$$
for any $Y$ and any $R\subseteq (Y\times S)$.\par

Given a function $\beta:Y\to S$ and taking $R=\Omega_\beta:=\{\big(y,\beta(y)\big)\mid y\in Y\}$,
one can check easily that $\Omega_\beta\Id_{S}=\beta$. It follows that
\begin{equation}\label{Phi}
\Phi_Y(\beta)=\Omega_\beta\,\varphi\mpoint
\end{equation}
Hence $\Phi$ is entirely determined by $\varphi$. 
Now Condition~(\ref{fonctoriel}) is fulfilled if and only if, for any finite sets $X$ and $Y$,
any correspondence $R\subseteq (Y\times X)$, and any function $\alpha:X\to S$, we have 
$$R\Omega_\alpha(\varphi)=\Omega_{R\alpha}(\varphi)\mpoint$$
In other words 
$$\sum_{\lambda}u_\lambda R\Omega_\alpha(\lambda)=\sum_{\lambda}u_\lambda\Omega_{R\alpha}(\lambda)\mpoint$$
Hence Condition~(\ref{fonctoriel}) is satisfied if and only if, for any finite sets $X$ and $Y$, any correspondence $R\subseteq (Y\times X)$, any function $\alpha:X\to S$, and any function $\psi:Y\to T$, we have
\begin{equation}\label{fonctoriel2}
\sum_{R\Omega_\alpha(\lambda)=\psi}u_\lambda=\sum_{\Omega_{R\alpha}(\lambda)=\psi}u_\lambda\mpoint
\end{equation}
But for $y\in Y$
\begin{eqnarray*}
R\Omega_\alpha(\lambda)(y)&=&\mbigvee{(y,s)\in R\Omega_\alpha}\lambda(s)\\
&=&\mbigvee{(y,x)\in R}\lambda\alpha(x)\mpoint
\end{eqnarray*}
On the other hand
\begin{eqnarray*}
\Omega_{R\alpha}(\lambda)(y)&=&\mbigvee{(y,s)\in \Omega_{R\alpha}}\lambda(s)\\
&=&\lambda\big(R\alpha(y)\big)\\
&=&\lambda\big(\mbigvee{(y,x)\in R}\alpha(x)\big)\mpoint
\end{eqnarray*}
Now take $X=S$ and $\alpha=\Id_{S}$ in~(\ref{fonctoriel2}). Then let $Y=\CB({S})$ be the set of subsets of~$S$ and let $R\subseteq (Y\times S)$ be the set of pairs $(A,s)$, where $A\subseteq {S}$ and $s\in A$.\par

Then for a given map $\lambda:{S}\to{T}$, let us {\em define} $\psi:Y\to{T}$ by $\psi=\Omega_{R\Id_S}(\lambda)$,
in other words
$$\forall A\subseteq {S},\;\psi(A)=\lambda(\mbigvee{s\in A}s)\mpoint$$
Suppose that there exists $\lambda':{S}\to{T}$ such that $\Omega_{R\Id_S}(\lambda')=\psi$. Then for $A\subseteq S$
$$\psi(A)=\lambda'(\mbigvee{s\in A}s)\mpoint$$
Taking $A=\{s\}$, it follows that $\lambda'=\lambda$.
Hence in~(\ref{fonctoriel2}) with our specific choices, the right hand side is simply equal to $u_\lambda$.\par

On the other hand the left hand side is equal to the sum of $u_{\lambda'}$, for all $\lambda'$ such that
$R\lambda'=\psi$, that is, satisfying
$$\forall A\subseteq {S},\; \psi(A)=\mbigvee{s\in A}\lambda'(s)\mpoint$$
Again, taking $A=\{s\}$, it follows that $\lambda'=\lambda$.
With our specific choices, the left hand side of~(\ref{fonctoriel2}) is then equal to $u_\lambda$
if and only if $R\lambda=\Omega_{R\Id_S}(\lambda)$,
that is, for any $A\subseteq S$,
$$\mbigvee{s\in A}\lambda(s)=\lambda(\mbigvee{s\in A}s)\mpoint$$
If this condition is not satisfied, then the left hand side of~(\ref{fonctoriel2}) is zero (empty sum).
In other words $u_\lambda=0$ if $\lambda$ is not a morphism in the category~$\CL$
(and this is where we see the relevance of the definition of morphisms in~$\CL$).\par

It follows that $\varphi=\sum\limits_{\lambda}u_\lambda\lambda$ is a morphism in~$k\CL$, from $S$ to~$T$.
We claim that the image of this morphism via the functor $F_?$ is equal to~$\Phi$
and this will prove that the functor $F_?:k\CL\to\CF_k$ is full.
To prove the claim, notice that, for any function $\beta:Y\to S$, we have
$$F_{\varphi,Y}(\beta)=\sum_{\lambda}u_\lambda F_{\lambda,Y}(\beta)
=\sum_{\lambda}u_\lambda (\lambda\circ\beta) = \sum_{\lambda}u_\lambda \Omega_\beta \, \lambda
=\Omega_\beta \, \varphi= \Phi_Y(\beta) \mvirg$$
using the equation~(\ref{Phi}). This proves the claim and completes the proof that $F_?$ is full.\par

It remains to show that the functor $F_?$ is faithful. So let $\varphi$ and $\psi$ be two linear combinations of morphisms $S\to T$ in~$\CL$, which induce the same morphism $\theta=F_\varphi=F_\psi:F_S\to F_T$. Evaluating this morphism at the set $S$ gives a map $\theta_S:F_S(S)\to F_T(S)$, and moreover
$$\theta_S(\Id_S)=F_{\varphi,S}(\Id_S)=\varphi\circ\Id_S=\varphi \in F_T(S)=k(T^S) \mpoint$$
For the same reason, $\theta_S(\Id_S)=F_{\psi,S}(\Id_S)=\psi$,
hence $\varphi=\psi$. This completes the proof of Theorem~\ref{fully faithful}.
\endpf

The connection between finite lattices and correspondence functors also has the following rather remarkable feature.

\result{Theorem} \label{projective-distributive} Let $T$ be a finite lattice. The functor $F_T$ is projective in $\CF_k$ if and only if $T$ is distributive.
\fresult

\pf Let $\CB({T})$ be the lattice of subsets of ${T}$. Let $\upsilon:\CB({T})\to{T}$ be the morphism in the category $\CL$ defined by
$$\forall A\subseteq {T},\;\upsilon(A)=\mbigvee{t\in A}t\mpoint$$
This morphism induces a morphism of functors $F_\upsilon:F_{\CB(\mathcal{T})}\to F_{T}$,  and $F_\upsilon$ is surjective~: indeed, if $X$ is a finite set and $\alpha:X\to {T}$ is a function,
and if we define $\hat{\alpha}:X\to \CB({T})$ by
$$\forall x\in X,\;\hat{\alpha}(x)=\{\alpha(x)\}\mvirg$$
then, for any $x\in X$
$$F_\upsilon(\hat{\alpha})(x)=(\upsilon\circ\hat{\alpha})(x)=\mbigvee{t\in \hat{\alpha}(x)}t=\alpha(x)\mvirg$$
thus $F_\upsilon(\hat{\alpha})=\alpha$, so $F_\upsilon$ is surjective.\par

Now if $F_T$ is projective, then the morphism $F_\upsilon$ splits and there exists a morphism $\Phi:F_T\to F_{\CB({T})}$  such that $F_\upsilon\circ \Phi$ is the identity morphism of $F_{T}$. It follows from Theorem~\ref{fully faithful} that $\Phi$ is of the form $\sum_{\sigma\in M}\limits u_\sigma F_\sigma$, where $M$ is a finite set of morphisms $\sigma:{T}\to \CB({T})$ in $\CL$, and $u_\sigma\in k$. Moreover $F_\upsilon\circ \Phi$ is then equal to $\sum_{\sigma\in M} \limits u_\sigma F_{\upsilon\circ\sigma}$, hence there exists at least one such $\sigma\in M$ such that $\upsilon\circ\sigma$ is equal to the identity of~${T}$. This means that
$$\forall t\in {T},\;t=\mbigvee{x\in\sigma(t)}x\mpoint$$
In particular $\sigma(t)\subseteq [\hat{0},t]_{T}$ for any $t\in{T}$. Then for $r,s\in{T}$
$$[\hat0,r\wedge s]_{T}=[\hat0,r]_{T}\cap [\hat0,s]_{T}\supseteq \sigma(r)\cap\sigma(s)\supseteq \sigma(r\wedge s)\mvirg$$
because $\sigma$ is order-preserving. It follows that
$$r\wedge s\;\geq\mbigvee{x\in \sigma(r)\cap\sigma(s)}x \; \geq \mbigvee{x\in \sigma(r\wedge s)}x=r\wedge s\mvirg$$
hence
$$r\wedge s=\mbigvee{x\in \sigma(r)\cap\sigma(s)}x\mpoint$$
Now, since $\sigma$ preserves joins, we obtain, for all $r,s,t\in{T}$,
\begin{eqnarray*}
t\wedge(r\vee s)&=&\mbigvee{x\in \sigma(t)\cap\sigma(r\vee s)}x\\
&=&\mbigvee{x\in\sigma(t)\cap\big(\sigma(r)\cup\sigma(s)\big)}x\\
&=&\mbigvee{x\in\big(\sigma(t)\cap\sigma(r)\big)\cup\big(\sigma(t)\cap\sigma(s)\big)}x\\
&=&\Big(\mbigvee{x\in\sigma(t)\cap\sigma(r)}x\Big)\vee \Big(\mbigvee{x\in\sigma(t)\cap\sigma(s)}x\Big)\\
&=&(t\wedge r)\vee(t\wedge s)\mpoint
\end{eqnarray*}
In other words the lattice ${T}$ is distributive.\par

Conversely, by Lemma~\ref{principal}, any finite distributive lattice $T$ is isomorphic to the lattice $\Idown(E,R)$ of lower ideals of a finite poset $(E,R)$. By Proposition~\ref{C(-,E)R}, the associated functor $F_T$ is projective in $\CF_k$. This completes the proof of Theorem~\ref{projective-distributive}.
\endpf


\section{Quotients of functors associated to lattices} \label{Section-quotients}

\bigskip
\noindent
We now introduce, for any finite lattice $T$, a subfunctor of $F_T$ naturally associated with the set of irreducible elements of~$T$.

\result{Notation} \label{H_T} Let $T$ be a finite lattice and let $(E,R)$ be the full subposet of its irreducible elements.
For a finite set $X$, let $H_T(X)$ denote the $k$-submodule of $F_T(X)=k(T^X)$ generated by all functions $\varphi:X\to T$ such that $E\nsubseteq \varphi(X)$.
\fresult

\result{Proposition} \label{F_T/H_T} Let $T$ be a finite lattice and let $(E,R)$ be the full subposet of its irreducible elements.
\begin{enumerate}
\item The assignment sending a finite set $X$ to $H_T(X)\subseteq F_T(X)$ is a subfunctor $H_T$ of $F_T$.
\item The evaluation $(F_T/H_T)(X)$ has a $k$-basis consisting of (the classes of) all functions $\varphi:X\to T$ such that $E\subseteq \varphi(X)$. 
\item The set~$E$ is a minimal set for the functor $F_T/H_T$ .
\end{enumerate}
\fresult

\pf (a) Let $X$ and $Y$ be finite sets, let $Q\in\CC(Y,X)$ be a correspondence, and let $\varphi:X\to T$ be a function.
Then
\begin{equation}\label{lemma H_T}(Q\varphi)(Y)\cap E\subseteq \varphi(X)\cap E\mpoint
\end{equation}
Indeed, if $e\in E$ and $e=(Q\varphi)(y)$, for $y\in Y$, then
$$e=\mbigvee{(y,x)\in Q}\varphi(x)\mpoint$$
As $e$ is irreducible in $T$, there exists $x\in X$ such that $(y,x)\in Q$ and $e=\varphi(x)$, and~(\ref{lemma H_T}) follows.\par

In particular, if $\varphi(X)\cap E$ is a proper subset of $E$, then $(Q\varphi)(Y)\cap E$ is a proper subset of $E$. Hence $H_T$ is a subfunctor of $F_T$.\mpn

(b) This follows from the definitions of $F_T$ and~$H_T$.\mpn

(c) If $|X|<|E|$, then $|\varphi(X)|<|E|$ and therefore there is no map $\varphi:X\to T$ such that $E\subseteq \varphi(X)$.
Consequently, the $k$-basis of~(b) is empty in that case, so that $(F_T/H_T)(X)=\{0\}$.
Now if $X=E$, then the $k$-basis of~(b) consists of all bijections $E\to E$ (followed by the inclusion map $E\to T$), so that $(F_T/H_T)(E) \cong k\Sigma_E$.
This shows that $E$ is a minimal set for~$F_T/H_T$.
\endpf

The quotient functor $F_T/H_T$ plays a important role in our work, in particular in Theorem~\ref{surjection} and for the description of the fundamental functors and the simple functors in~\cite{BT3}.
We now give another characterization of~$H_T(X)$.

\result{Proposition} \label{T=Iup(E,R)}
Let $T=\Iup(E,R)$ for a finite poset $(E,R)$ and let $X$ be a finite set.
\begin{enumerate}
\item Under the isomorphism $F_T\to k\CC({-},E)R$ of Proposition~\ref{C(-,E)R}, $H_T(X)$ is isomorphic to the $k$-submodule of $k\CC(X,E)R$ generated by the correspondences $S$ which have no retraction, that is, for which there is no $U\in\CC(E,X)$ such that $US=R$.
\item Under the isomorphism $F_T\to k\CC({-},E)R$ of Proposition~\ref{C(-,E)R}, the image of $F_T(X)/H_T(X)$ is a free $k$-module with basis consisting of all the correspondences $S\in \CC(X,E)R$ which have a retraction $U\in\CC(E,X)$.
\end{enumerate}
\fresult

\pf
By Proposition~\ref{C(-,E)R}, the functor $F_T$ is isomorphic to the functor $k\CC({-},E)R$ by sending, for a finite set $X$, a function $\varphi:X\to \Iup(E,R)$ to the correspondence $\Gamma_{\varphi}=\{(x,e)\in X\times E\mid e\in\varphi(x)\}$.\mpn

(a)
The set $E^{\uparrow}$ of irreducible elements of the lattice $\Iup(E,R)$ is the set of principal upper ideals 
$$[e,\cdot[_{R}=\big\{f\in E\mid (e,f)\in R\big\}\mvirg$$
for $e\in E$.
Let $\varphi:X\to \Iup(E,R)$ be such that $\varphi\notin H_T(X)$, that is, $\varphi(X)\supseteq E^{\uparrow}$. Then, for each $e\in E$, there exists $x_e\in X$ such that $\varphi(x_e)=[e,\cdot[_R$. Let $U\in\CC(E,X)$ be defined by
$$U=\big\{(e,x_e)\mid e\in E\big\}\subseteq E\times X\mpoint$$
Then for any $e\in E$
$$(U\varphi)(e)=\bigcup_{(e,x)\in U}\varphi(x)=\varphi(x_e)=[e,\cdot[_R\mpoint$$
By Lemma~\ref{SGamma}, it follows that
$$U\Gamma_{\varphi}=\Gamma_{U\varphi\,}=\big\{(e,f)\in E\times E\mid f\in[e,\cdot[_R\big\}=R\mvirg$$
so $\Gamma_\varphi$ has a retraction.\par

Conversely, let $S\in \CC(X,E)R$ be a correspondence such that there exists a correspondence $U\in\CC(X,E)$ with $U S=R$. Then $S=\Gamma_{\varphi}$, where $\varphi:X\to \Iup(E,R)$ is the function defined by $\varphi(x)=\big\{e\in E\mid (x,e)\in S\big\}$, for any $x\in X$. It follows that $US=
\Gamma_{U\varphi\,}=R$, or in other words
$$\forall e,f\in E,\;(e,f)\in R\iff \exists x\in X,\;(e,x)\in U,\;(x,e)\in S\mpoint$$
As $\Delta_E\subseteq R$, for any $e\in E$, there exists $x_e\in X$ such that $(e,x_e)\in U$ and $(x_e,e)\in S$. Moreover if $(x_e,f)\in S$, then $(e,f)\in R$, and conversely, if $(e,f)\in R$, then $(x_e,f)\in SR=S$.
In other words, $f\in \varphi(x_e)$ if and only if $(e,f)\in R$.
It follows that $\varphi(x_e)=[e,\cdot[_R$, hence $\varphi(X)\supseteq E^{\uparrow}$.
This proves that $\varphi\notin H_T(X)$.\mpn

(b) This follows from (a).
\endpf

\begin{rem}{Remark}
In the special case when $R=\Delta_E$ is the equality relation, then $\CC(X,E)R=\CC(X,E)$ and a retraction of $S\in\CC(X,E)$ is a correspondence $U\in\CC(E,X)$ such that $US=\Id_E$ (a retraction in the usual sense).
Moreover, if $S\in \CC(X,E)$ has a retraction, then $S$ is a monomorphism in the category~$\CC$.
It can be shown conversely that any monomorphism in the category~$\CC$ has a retraction.
Thus in this case, the evaluation $F_T(X)/H_T(X)$ of the quotient functor $F_T/H_T$
has a $k$-basis consisting of all the monomorphisms in~$\CC(X,E)$.
\end{rem}

\bigbreak
In order to deal with quotients of the functor $F_T$, we need information on morphisms starting from~$F_T$.
We first need a lemma.

\result{Lemma} \label{surgamma}  
Let $T$ be a finite lattice, let $(E,R)$ be the full subposet of its irreducible elements, and
let $\iota:E\to T$ denote the inclusion map.
\begin{enumerate}
\item If $\varphi:X\to T$ is a function, then $\Gamma_\varphi\,\iota=\varphi$ and $\Gamma_\varphi R\op=\Gamma_\varphi$,
where $\Gamma_{\varphi}$ is defined in Notation~\ref{Gamma}.
\item $\Gamma_\iota=R\op$.
\item The correspondence functor $F_T$ is generated by $\iota\in F_T(E)$.
\end{enumerate}
\fresult

\pf
(a) By definition, the map $\Gamma_\varphi\,\iota:X\to T$ satisfies
$$\forall x\in X,\;\;(\Gamma_\varphi\, \iota)(x)=\mbigvee{(x,e)\in \Gamma_\varphi}\iota(e)
=\mbigvee{e\leq_T\varphi(x)}e\;=\varphi(x)\mvirg$$
as any element $t$ of $T$ is equal to the join of the irreducible elements of $T$ smaller than~$t$.
Thus we have $\Gamma_\varphi\,\iota=\varphi$.
The equality $\Gamma_\varphi R\op=\Gamma_\varphi$ was proved in Lemma~\ref{rightR}.\mpn

(b) $\Gamma_\iota=\{(x,e)\in E\times E\mid e\leq_T \iota(x)\}=\{(x,e)\in E\times E\mid e\leq_R x\}=R\op$.\mpn

(c) For every function $\varphi:X\to T$, we have $\varphi=\Gamma_\varphi\iota$ by~(a).
Therefore $F_T$ is generated by $\iota\in F_T(E)$.
\endpf

\result{Proposition} \label{HomFT}  
Let $T$ be a finite lattice, let $(E,R)$ be the full subposet of its irreducible elements, and
let $\iota:E\to T$ denote the inclusion map. Let $M$ be a correspondence functor.
\begin{enumerate}
\item The $k$-linear map
$$\Hom_{\CF_k}(F_T, M) \longrightarrow M(E) \,,\quad \Phi \mapsto \Phi_E(\iota)$$
is injective. Its image is contained in the $k$-submodule
$$R\op M(E)=\{ m\in M(E) \mid R\op m=m \} \mpoint$$
\item If $T$ is distributive, so that $T\cong \Idown(E,R)$, then the image of the above map is equal to $R\op M(E)$,
so that $\Hom_{\CF_k}(F_T, M) \cong R\op M(E)$ as $k$-modules.
\end{enumerate}
\fresult

\pf
(a) By Lemma~\ref{surgamma}, for any $\Phi\in\Hom_{\CF_k}(F_T, M)$ and any map $\varphi: X\to T$, we have
$$\Phi_X(\varphi)=\Phi_X(\Gamma_\varphi\,\iota)=\Gamma_\varphi\,\Phi_E(\iota) \mpoint$$
This shows that $\Phi$ is entirely determined by $\Phi_E(\iota)$,
proving the injectivity of the map $\Phi\mapsto \Phi_E(\iota)$.\par

Moreover, $\Phi_E(\iota)=\Gamma_\iota\Phi_E(\iota)=R\op\Phi_E(\iota)$, because $\Gamma_\iota=R\op$ by Lemma~\ref{surgamma}.
Therefore $\Phi_E(\iota)$ is contained in~$R\op M(E)$.\mpn

(b) Since $T$ is distributive, we have
$$\Gamma_{Q\varphi} = Q \Gamma_\varphi$$
by Lemma~\ref{SGamma}.
Now given $m\in R\op M(E)$, we can define $\Phi:F_T\to M$ by setting
$$\Phi_X(\varphi)=\Gamma_\varphi m \;,\quad \forall \varphi:X\to T \mpoint$$
This is indeed a natural transformation of functors since
$$\Phi_Y(Q\varphi)=\Gamma_{Q\varphi} m=Q \Gamma_\varphi m=Q\Phi_X(\varphi)$$
for any correspondence $Q\subseteq Y\times X$.
Moreover,
$$\Phi_E(\iota)=\Gamma_\iota m=R\op m=m \mvirg$$
because $m\in R\op M(E)$ by assumption and $R\op$ is idempotent.
Thus $m$ is indeed in the image of the map $\Phi\mapsto \Phi_E(\iota)$.
\endpf

When $k$ is a field, we wish to give some information on simple functors $S_{F,Q,V}$ appearing as quotients of~$F_T$.
We prove a more general result over an arbitrary commutative ring~$k$,
involving the not necessarily simple functors $S_{F,Q,V}$ introduced in Section~\ref{Section-functors}.

\result{Theorem} \label{simple-quotients}   
Let $T$ be a finite lattice and let $(E,R)$ be the full subposet of its irreducible elements.
Let $(F,Q)$ be a poset and let $V$ be a $k\Aut(F,Q)$-module generated by a single element (e.g. a simple module).
\begin{enumerate}
\item  If $S_{F,Q,V}$ is isomorphic to a quotient of~$F_T$ , then $|F|\leq |E|$.
\item Assume that $F=E$. If $S_{E,Q,V}$ is isomorphic to a quotient of~$F_T$, then there exists a permutation $\sigma\in\Sigma_E$ such that $R\op\subseteq \ls\sigma Q$.
\item Assume that $F=E$ and that $T$ is distributive, so that $T\cong \Idown(E,R)$.
Then $S_{E,Q,V}$ is isomorphic to a quotient of~$F_T$ if and only if there exists a permutation $\sigma\in\Sigma_E$ such that $R\op\subseteq \ls\sigma Q$.
\end{enumerate}
\fresult

\pf
(a) If $S_{F,Q,V}$ is isomorphic to a quotient of~$F_T$, then $\Hom_{\CF_k}(F_T, S_{F,Q,V})\neq\{0\}$,
so we have $S_{F,Q,V}(E)\neq\{0\}$ by Proposition~\ref{HomFT}.
But $F$ is a minimal set for $S_{F,Q,V}$ by Proposition~\ref{SER-SERV}, so $|F|\leq |E|$.\mpn

(b) If $S_{E,Q,V}$ is isomorphic to a quotient of~$F_T$, there exists a nonzero morphism $\Phi: F_T\to S_{E,Q,V}$.
By Proposition~\ref{HomFT}, $\Phi_E(\iota)=m\neq 0 \in R\op S_{E,Q,V}(E)$. By Proposition~\ref{SER-SERV}, we know that
$$S_{E,Q,V}(E)\cong T_{Q,V}=\CP_E f_Q\otimes_{k\Aut(E,Q)} V$$
and $\CP_E f_Q$ is a free $k$-module with basis $\{\Delta_\sigma f_Q \mid \sigma\in\Sigma_E\}$, by Proposition~\ref{fundamental-module}.
Thus we can write
$$m=\sum_{\sigma\in\Sigma_E}\lambda_\sigma\Delta_\sigma f_Q\otimes v \qquad (\lambda_\sigma\in k) \mpoint$$
Since $m\in R\op S_{E,Q,V}(E)$, we have $R\op m=m$ and so there exists $\sigma\in\Sigma_E$ such that
$R\op \Delta_\sigma f_Q\neq0$. Hence $R\op \Delta_\sigma f_Q\Delta_\sigma^{-1}\neq0$, that is,
$R\op\subseteq \ls\sigma Q$, by Theorem~6.2 of~\cite{BT1}.\mpn

(c) One implication follows from (b). Assume now that there exists a permutation $\sigma\in\Sigma_E$ such that $R\op\subseteq \ls\sigma Q$.
We first note that $S_{E,Q,V}$ is generated by $f_Q\otimes v\in S_{E,Q,V}(E)\cong \CP_E f_Q\otimes_{k\Aut(E,Q)} V$, where $v$ is a generator of~$V$.
This follows from the definition of $S_{E,Q,V}$ as a quotient of $L_{E,T_{Q,V}}$ and the fact that any functor $L_{E,W}$ is generated by $L_{E,W}(E)=W$ by definition.\par

It is easy to see that $S_{E,Q,V}\cong S_{E,\ls\sigma Q,\ls\sigma V}$ for any $\sigma\in\Sigma_E$,
because $T_{Q,V}\cong T_{\ls\sigma Q,\ls\sigma V}$ by construction (see Theorem~8.1 in~\cite{BT1} for more details).
Since $R\op$ is contained in a conjugate of~$Q$, we can assume that $R\op\subseteq Q$.
This is equivalent to $R\op f_Q=f_Q$, by Theorem~6.2 of~\cite{BT1}.\par

Thus $m=f_Q\otimes v\in S_{E,Q,V}(E)$ is invariant under left multiplication by~$R\op$. By Proposition~\ref{HomFT} and the assumption that $T$ is distributive, there exists a morphism $\Phi: F_T\to S_{E,Q,V}$ such that $\Phi_E(\iota)=f_Q\otimes v$. Since this is nonzero and generates $S_{E,Q,V}$, this functor is isomorphic to a quotient of~$F_T$.
\endpf

The similar question of finding fundamental functors appearing as quotients of~$F_T$ will be considered later in Theorem~\ref{fundamental-quotient}.


\section{The fundamental functor associated to a poset} \label{Section-fundamental}

\bigskip
\noindent
The fundamental functor $\S_{E,R}$ associated to a poset $(E,R)$ was introduced in Section~\ref{Section-functors}.
One of our important goals is to give a precise description of its evaluations and use it to deduce a precise description of the evaluations of simple functors, but this will be fully achieved only in~\cite{BT3}.
We prepare the ground by proving several main results about~$\S_{E,R}$.
Recall from Proposition~\ref{SER-SERV} that $E$ is a minimal set for~$\S_{E,R}$ and that $\S_{E,R}(E)$ is isomorphic to the fundamental module $\CP_E f_R$, which is described in Proposition~\ref{fundamental-module}.\par

Since $\S_{E,R}=L_{E,\CP_Ef_R}/J_{E,\CP_Ef_R}$ by definition, it is important to know when an element of~$L_{E,\CP_Ef_R}(X)$ belongs to~$J_{E,\CP_Ef_R}(X)$, where $X$ is a finite set. For this analysis, we note that an element of~$L_{E,\CP_Ef_R}(X)$ is written
$$\sum_{S\in\CC(X,E)} \lambda_S S\otimes f_R
\in k\CC(X,E) \otimes_{k\CC(E,E)} \CP_Ef_R = L_{E,\CP_Ef_R}(X) \mvirg$$
where $\lambda_S\in k$ for every $S$, because the tensor product is over $k\CC(E,E)$ and $\CP_E$ is a quotient algebra of~$k\CC(E,E)$.
But since $Rf_R=f_R$ (by Proposition~\ref{fundamental-module}) and the tensor product is over $k\CC(E,E)$, we can replace $S$ by $SR$ and obtain a sum running only over $S\in\CC(X,E)R$.

\result{Lemma} \label{belong-to-J}
Consider an element
$$\sum_{S\in\CC(X,E)R} \lambda_S S\otimes f_R \in L_{E,\CP_Ef_R}(X) \mpoint$$
This element belongs to~$J_{E,\CP_Ef_R}(X)$ if and only if
$$\forall \, U\in R\CC(E,X) \,, \quad \sum_{\substack{S\in\CC(X,E)R \\ US=R}} \lambda_S=0 \mpoint$$
\fresult

\pf
By the definition of~$J_{E,\CP_Ef_R}(X)$, we have
$$\begin{array}{rl}
\displaystyle \sum_{S\in\CC(X,E)R} \lambda_S S\otimes f_R &\in \; J_{E,\CP_Ef_R}(X) \;\iff\\
&
\forall \, U\in \CC(E,X) \,, \quad \displaystyle\sum_{S\in\CC(X,E)R} \lambda_S US\cdot f_R =0 \mpoint
\end{array}$$
By Proposition~\ref{fundamental-module}, the action of the relation $US\in \CC(E,E)$ on~$f_R$ is given by
$$US\cdot f_R=\left\{\begin{array}{ll}
\Delta_{\tau}f_R&\hbox{if}\;\;\exists\tau\in\Sigma_E\;\hbox{such that}\;
\Delta_E\subseteq \Delta_{\tau^{-1}}US\subseteq R,\\
0&\hbox{otherwise}\mpoint\end{array}\right.$$
We claim that
$$\Delta_E\subseteq \Delta_{\tau^{-1}}US\subseteq R \;\iff \; US=\Delta_\tau R \mpoint$$
If the left hand side holds, then multiply on the right by~$R$ and use the fact that $SR=S$ and $R^2=R$ (by transitivity and reflexivity of~$R$) to obtain $\Delta_{\tau^{-1}}US= R$, hence $US= \Delta_\tau R$.
Conversely, if the right hand side holds, then $\Delta_{\tau^{-1}}US= R$, hence
$$R\Delta_{\tau^{-1}}US= R^2=R$$
by transitivity and reflexivity of~$R$. In particular, by reflexivity again,
$$\Delta_E\subseteq R\Delta_{\tau^{-1}}US$$
so that, for any $(a,a)\in\Delta_E$, there exists $b\in E$ with $(a,b)\in R$
and $(b,a)\in\Delta_{\tau^{-1}}US=R$.
By antisymmetry of~$R$, it follows that $b=a$ and therefore
$(a,a)\in\Delta_{\tau^{-1}}US$,
so that $\Delta_E\subseteq \Delta_{\tau^{-1}}US$.
This shows that the left hand side holds, proving the claim.\par

It follows that our given element belongs to $J_{E,\CP_Ef_R}(X)$ if and only if
\begin{equation} \label{equation-for-all-U}
\forall \, U\in \CC(E,X) \,, \quad
\sum_{\substack{S\in\CC(X,E)R \\ \tau\in\Sigma_E \\ US=\Delta_\tau R}} \lambda_S \Delta_\tau f_R =0 \mpoint
\end{equation}
But by Proposition~\ref{fundamental-module}, $\CP_Ef_R$ is a free $k$-module with basis $\{\Delta_\tau f_R\mid \tau\in\Sigma_E\}$. Therefore (\ref{equation-for-all-U}) is equivalent to
$$\forall \, U\in \CC(E,X) \,, \quad \forall \, \tau\in\Sigma_E \,,\quad
\sum_{\substack{S\in\CC(X,E)R \\  US=\Delta_\tau R}} \lambda_S =0 \mpoint$$
Replacing $U$ by $\Delta_{\tau^{-1}}U$, this is equivalent to
\begin{equation} \label{equation-without-tau}
\forall \, U\in \CC(E,X) \,, \quad
\sum_{\substack{S\in\CC(X,E)R \\ US= R}} \lambda_S  =0 \mpoint
\end{equation}
Now we claim that
$$US=R \Longleftrightarrow RUS=R \mpoint$$
If the left hand side holds, then multiply on the left by~$R$ and use the fact that $R^2=R$. 
Conversely, if the right hand side holds, then $US=\Delta_EUS\subseteq RUS=R$.
Moreover, $\Delta_E\subseteq R=RUS$, so that, for any $(a,a)\in\Delta_E$, there exists $b\in E$ with $(a,b)\in R$
and $(b,a)\in US \subseteq R$.
By antisymmetry of~$R$, it follows that $b=a$ and therefore
$(a,a)\in US$,
so that $\Delta_E\subseteq US$. Multiplying on the right by $R$ and using the fact that $SR=S$, we obtain
$R\subseteq US$. This shows that $R=US$, as required, proving the claim.\par

It now follows that in the family of equations~(\ref{equation-without-tau}), the equation for $U$ is the same as the equation for~$RU$. Therefore we can assume that $U\in R\CC(E,X)$ and (\ref{equation-without-tau}) is equivalent to
$$\forall \, U\in R\CC(E,X) \,, \quad
\sum_{\substack{S\in\CC(X,E)R \\ US= R}} \lambda_S  =0 \mvirg
$$
proving the lemma.
\endpf

Before stating the main theorem, we need another lemma.

\result{Lemma} \label{W-phi}
Let $T$ be a finite lattice, let $(E,R)$ be the full subposet of its irreducible elements, and let $\iota:E\to T$ denote the inclusion map.
For any correspondence $W\in \CC(E,X)$ and any map $\varphi:X\to T$, the following two statements are equivalent~:
$$W\Gamma_\varphi =R\op \iff W\varphi=\Id_E \mvirg$$
where $\Gamma_\varphi\in\CC(X,E)$ is defined in Notation~\ref{Gamma}.
\fresult

\pf
Let $e,f\in E$. Since $(E,R)$ is a full subposet of~$T$, the relation $(f,e)\in R\op$ is equivalent to $e\leq f$,
where we write $\leq$ instead of $\leq_T$ for simplicity.
Suppose that $W\Gamma_\varphi =R\op$. Then $e\leq f$ if and only if
$$\exists\, x\in X \,\text{ with } \, (f,x)\in W \,, \; (x,e)\in \Gamma_\varphi \iff
\exists\, x\in X \,\text{ with } \, (f,x)\in W \,,\; e\leq \varphi(x)\mpoint$$
Now we can write
$$\varphi(x)=\bigvee_{\substack{e\in E \\ e\leq \varphi(x)}}e$$
because any element of~$T$ is a join of irreducible elements. It follows that
$$(W\varphi)(f)=\bigvee_{\substack {x\in X \\ (f,x)\in W}} \varphi(x)
=\bigvee_{\substack {x\in X \\ (f,x)\in W}} \bigvee_{\substack{e\in E \\ e\leq \varphi(x)}} e=\bigvee_{e\leq f}e \mvirg$$
by the equivalence above. Thus $(W\varphi)(f)=f$, so that $W\varphi=\Id_E$.\par

Conversely, suppose that $W\varphi=\Id_E$. Since $\varphi=\Gamma_\varphi\iota$ by Lemma~\ref{surgamma},
we have $W\Gamma_\varphi\iota=\Id_E$ and we obtain, for every $e\in E$,
$$e=(W\Gamma_\varphi\iota)(e)=\bigvee_{\substack {f\in E \\ (e,f)\in W\Gamma_\varphi}} \iota(f)
=\bigvee_{\substack {f\in E \\ (e,f)\in W\Gamma_\varphi}} f \mpoint$$
Since $e$ is irreducible, it follows that $e=f$ for some~$f\in E$ such that $(e,f)\in W\Gamma_\varphi$, and so $(e,e)\in W\Gamma_\varphi$, showing that $\Delta_E\subseteq W\Gamma_\varphi$.
On the other hand, we also obtain $f\leq e$ for every $f\in E$ such that $(e,f)\in W\Gamma_\varphi$, showing that $W\Gamma_\varphi \subseteq R\op$.
Therefore $\Delta_E\subseteq W\Gamma_\varphi\subseteq R\op$. Multiplying on the right by~$R\op$, we deduce that $W\Gamma_\varphi= R\op$, thanks to the fact that $\Gamma_\varphi R\op=\Gamma_\varphi$ by Lemma~\ref{rightR}.
\endpf

The following theorem establishes the link between the functor $F_T$ associated with a lattice~$T$ and the fundamental correspondence functors.

\result{Theorem} \label{surjection} Let $T$ be a finite lattice, let $(E,R)$ be the full subposet of its irreducible elements, and let $\iota:E\to T$ denote the inclusion map.
\begin{enumerate}
\item There exists a unique morphism of correspondence functors
$$\Theta_T:F_T\longrightarrow \S_{E,R\op}$$
mapping $\iota\in F_T(E)$ to $f_{R\op}\in\S_{E,R\op}(E)\cong \CP_Ef_{R\op}$.
Moreover, $\Theta_T$ is surjective.
\item The subfunctor $H_T$ of $F_T$ is contained in the kernel of~$\Theta_T$. In particular, $\Theta_T$ induces a surjective morphism
$\overline\Theta_T:F_T/H_T\longrightarrow \S_{E,R\op}$.
\end{enumerate}
\fresult

\pf (a) For any finite set $X$ and any function $\varphi: X\to T$, we let
$$\widetilde\Theta_{T,X}(\varphi)=\Gamma_\varphi \otimes f_{R\op} \in k\CC(X,E) \otimes_{k\CC(E,E)} f_{R\op} =L_{E,\CP_Ef_{R\op}}(X) \mvirg$$
where $\Gamma_\varphi\in\CC(X,E)$ is defined in Notation~\ref{Gamma}.
This extends to a $k$-linear map
$$\widetilde\Theta_{T,X} : F_T(X) \longrightarrow L_{E,\CP_Ef_{R\op}}(X)$$
which we compose with the canonical surjection
$$\pi_X:L_{E,\CP_Ef_{R\op}}(X) \longrightarrow L_{E,\CP_Ef_{R\op}}(X)/J_{E,\CP_Ef_{R\op}}(X)=\S_{E,R\op}(X)$$
to obtain a $k$-linear map
$$\Theta_{T,X} : F_T(X) \longrightarrow \S_{E,R\op}(X) \mpoint$$
The family of maps $\widetilde\Theta_{T,X}$ is not a morphism of functors,
but we are going to prove that it becomes so after composition with~$\pi:L_{E,\CP_Ef_{R\op}}\to \S_{E,R\op}$.
We have to show that, for any correspondence $V\in\CC(Y,X)$ and any map $\varphi: X\to T$, we have
$$\Theta_{T,Y}(V\varphi) = V\Theta_{T,X}(\varphi) \,,\quad\text{that is,} \quad
\widetilde\Theta_{T,Y}(V\varphi) - V\widetilde\Theta_{T,X}(\varphi) \in J_{E,\CP_Ef_{R\op}}(Y) \mpoint$$
In other words, we need to prove that
$$
(\Gamma_{V\varphi} - V\Gamma_\varphi) \otimes f_{R\op} \in J_{E,\CP_Ef_{R\op}}(Y) \mpoint$$
Since $\Gamma_{V\varphi} R\op=\Gamma_{V\varphi}$ and $\Gamma_\varphi R\op=\Gamma_\varphi$ by Lemma~\ref{rightR},
we can apply Lemma~\ref{belong-to-J} to the element $(\Gamma_{V\varphi} - V\Gamma_\varphi) \otimes f_{R\op}$ (with $R\op$ instead of~$R$).
Since we have a difference of two terms, we only have to show that the additional conditions in the equations occur simultaneously, that is,
\begin{equation} \label{V-Gamma}
\forall \, U\in R\op\CC(E,Y) \,, \quad U\Gamma_{V\varphi} =R\op \iff UV\Gamma_\varphi =R\op \mpoint
\end{equation}
We emphasize that the equation $V\Gamma_\varphi =\Gamma_{V\varphi}$ does not hold in general (it holds if the lattice~$T$ is distributive by Lemma~\ref{V-Gamma}), but it will become correct after left multiplication by $U\in R\op\CC(E,Y)$.
Using Lemma~\ref{W-phi}, we have to show that, for all $U\in R\op\CC(E,Y)$,
$$U(V\varphi) =\Id_E \iff (UV)\varphi =\Id_E \mpoint$$
But this is obvious in view of the action of correspondences on~$\varphi\in F_T(X)$.
Therefore $\Theta_T$ is a morphism of correspondence functors, as was to be shown.\par

By Lemma~\ref{surgamma}, the image of $\iota\in F_T(E)$ is
$$\Theta_{T,E}(\iota)=\Gamma_\iota\otimes f_{R\op}=R\op\otimes f_{R\op}=\Delta_E\otimes R\op f_{R\op}=\Delta_E\otimes f_{R\op} \mvirg$$
because the tensor product is over $k\CC(E,E)$
and because $R\op f_{R\op}= f_{R\op}$ by the action of relations on~$ f_{R\op}$ (see Proposition~\ref{fundamental-module}).
Now $\Delta_E\otimes f_{R\op}$ corresponds to $f_{R\op}$ under the isomorphism
$\S_{E,R\op}(E)\cong \CP_Ef_{R\op}$. Therefore $\Theta_T$ maps $\iota$ to~$f_{R\op}$.\par

Since $R\op  f_{R\op}= f_{R\op}$, any generator $S\otimes f_{R\op}$ of $L_{E,\CP_e f_{R\op}}(X)$, where $S\in \CC(X,E)$, can be written
$$S\otimes f_{R\op}= S\otimes R\op f_{R\op}=S R\op\otimes f_{R\op}=\Gamma_\varphi\otimes f_{R\op}$$
for some $\varphi:X\to T$, in view of Lemma~\ref{rightR}. This shows that $\Theta_T$ is surjective.
\mpn

(b) A generator of $H_T(X)$ is a map $\varphi: X\to T$ such that $E\not\subseteq \varphi(X)$.
We need to prove that $\varphi\in\Ker\Theta_T$, or in other words $\Gamma_\varphi\otimes f_{R\op} \in J_{E,\CP_e f_{R\op}}(X)$.
In order to prove this, we apply Lemma~\ref{belong-to-J} and we note that we have here a sum with a single term,
so we have to check that the additional condition $U\Gamma_\varphi=R\op$ never holds, where $U\in R\op\CC(E,X)$.
By Lemma~\ref{W-phi}, $U\Gamma_\varphi\neq R\op$ is equivalent to $U\varphi \neq\Id_E$.
But $U\varphi\in H_T(E)$ because $H_T$ is a subfunctor and obviously $\Id_E\notin H_T(E)$, so $U\varphi \neq\Id_E$, as required.
\endpf

\begin{rem}{Remark} \label{omega}
By Lemma~\ref{principal}, any lattice $T$ having $(E,R)$ as its full subposet of irreducible elements is isomorphic to a quotient of the distributive lattice $\Idown(E,R)$ via the map
$$\pi:\Idown(E,R) \longrightarrow T \,, \qquad A \mapsto \bigvee_{e\in A} e$$
where the join is taken in the lattice~$T$.
Since $\Idown(E,R)$ also has $(E,R)$ as its full subposet of irreducible elements, there is a surjective morphism
$$\Theta_{\Idown(E,R)} : F_{\Idown(E,R)} \longrightarrow \S_{E,R\op}$$
which is easily seen to correspond, via the isomorphism $F_{\Idown(E,R)}\cong k\CC(-,E)R\op$ of Proposition~\ref{C(-,E)R}, to the morphism
$$\omega: k\CC(-,E)R\op\longrightarrow \S_{E,R\op} \,, \qquad R\op \mapsto f_{R\op} \mpoint$$
Moreover, it is also elementary to check that the composite
$$\xymatrix{
k\CC(-,E)R\op \cong F_{\Idown(E,R)} \ar[r]^-{F_\pi} & F_T \ar[r]^-{\Theta_T} & \S_{E,R\op}
}$$
is equal to~$\omega$.
\end{rem}

\bigskip
For our next result, we need a well-known result of algebraic $K$-theory.

\result{Lemma} \label{K-theory} Let $L$ and $L'$ be finitely generated free $k$-modules of the same rank, where $k$ is a commutative ring. Then any surjective homomorphism of $k$-modules $f:L\to L'$ is an isomorphism.
\fresult

\pf 
The homomorphism $f$ splits, because $L'$ is free.
Let $M=\Ker(f)$ and let $r$ be the rank of~$L$, so that we have $k^r\cong k^r\oplus M$.
By localization at a prime ideal~$\wp$ of~$k$, we obtain $k_\wp^r\cong k_\wp^r\oplus M_\wp$.
The projective module $M_\wp$ is free because $k_\wp$ is local.
If $m$ is the rank of~$M_\wp$, we obtain $r=r+m$ by reducing modulo the maximal ideal $\wp k_\wp$.
Therefore $m=0$, that is, $M_\wp=\{0\}$. This holds for every prime ideal~$\wp$, hence $M=\{0\}$.
\endpf

We now show that the relationship between $F_T$ and $\S_{E,R\op}$ is very strong, in the sense that $\S_{E,R\op}$ is the only fundamental functor appearing as a quotient of~$F_T/H_T$, where $H_T$ is defined in Notation~\ref{H_T}. Recall that $E$ is a minimal set for the functor~$F_T/H_T$, by Proposition~\ref{F_T/H_T}.

\result{Theorem} \label{fundamental-quotient} Let $T$ be a finite lattice and let $(E,R)$ be the full subposet of its irreducible elements.
Let $(D,Q)$ be a finite poset such that there exists a surjective morphism $\Phi:F_T/H_T \to \S_{D,Q}$.
\begin{enumerate}
\item $|D|=|E|$ (so we can assume that $D=E$).
\item Assuming that $D=E$, there exists $\sigma\in\Sigma_E$ such that $\ls\sigma Q=R\op$. In other words, the posets $(E,Q)$ and $(E,R\op)$ are isomorphic.
\item The kernel of $\overline\Theta_T: F_T/H_T \to \S_{E,R\op}$ vanishes at~$E$ and
$$\overline\Theta_{T,E}: F_T(E)/H_T(E) \longrightarrow \S_{E,R\op}(E) \cong \CP_Ef_{R\op} $$
is an isomorphism.
\end{enumerate}
\fresult

\pf
(a) By Proposition~\ref{HomFT}, any morphism $F_T/H_T\to \S_{D,Q}$ is determined by an element of~$\S_{D,Q}(E)$. Since our given morphism $\Phi$ is surjective, it is nonzero, and
therefore $\S_{D,Q}(E) \neq\{0\}$. Since $D$ is a minimal set for $\S_{D,Q}$, we have $|D|\leq|E|$.
Now there is a surjective morphism
$$\Phi_D:F_T(D)/H_T(D) \to S_{D,Q}(D)\neq\{0\} \mvirg$$
so $F_T/H_T$ does not vanish at~$D$.
Since $E$ is a minimal set for $F_T/H_T$ by Proposition~\ref{F_T/H_T}, we have $|E|\leq |D|$.\mpn

(b) We prove that $\Phi_E:F_T(E)/H_T(E) \to S_{E,Q}(E)$ is an isomorphism.
By Proposition~\ref{H_T}, $F_T(E)/H_T(E)$ is a free $k$-module with basis the set of all bijections $E\to E$ (followed by the inclusion map $E\to T$).
By Proposition~\ref{fundamental-module}, $S_{E,Q}(E) \cong \CP_E f_Q$ is a free $k$-module with a basis consisting of the elements $\Delta_\sigma f_Q$,
where $\sigma$ runs through the group $\Sigma_E$ of all permutations of~$E$.
Thus the two modules under consideration are finitely generated free $k$-modules with the same rank 
and we know that the map $\Phi_E$ is surjective. Therefore $\Phi_E$ is an isomorphism by Lemma~\ref{K-theory}.\par

The same argument applies to the surjective morphism $\overline\Theta_{T,E}: F_T(E)/H_T(E) \to \S_{E,R\op}(E)$, which is therefore also an isomorphism.
It follows that there are isomorphisms of $k\CC(E,E)$-modules
$$\CP_E f_Q \cong S_{E,Q}(E) \cong F_T(E)/H_T(E) \cong \S_{E,R\op}(E) \cong \CP_E f_{R\op} \mpoint$$
This isomorphism maps $f_Q$ to $a f_{R\op}$ for some $a\in \CP_E$.
Therefore it maps $f_Q=Qf_Q$ to $a f_{R\op}=Qa f_{R\op}$ and in particular $Qa f_{R\op}\neq 0$.
By Proposition~\ref{fundamental-module}, this is possible only if there exists $\sigma,\tau\in\Sigma_E$ such that
$$\Delta_E\subseteq \Delta_{\tau^{-1}} Q\subseteq \ls \sigma (R\op) \mpoint$$
In particular, $\Delta_\tau\subseteq Q$, forcing $\tau=\Id$ because $Q$ is an order relation.
Thus $Q\subseteq \ls \sigma (R\op)$.
Swapping the role of $Q$ and~$R\op$, the same argument shows that $R\op$ is contained in a conjugate of~$Q$.
Therefore $Q$ and~$R\op$ are conjugate, as was to be shown.\mpn

(c) We have just seen that $\overline\Theta_{T,E}: F_T(E)/H_T(E) \longrightarrow \S_{E,R\op}(E)$
is an isomorphism. Therefore $\Ker\Theta_{T,E}=\{0\}$. In other words, the subfunctor $\Ker\Theta_T$ vanishes at~$E$.
\endpf


\section{The kernel of $\Theta_T$} \label{Section-kertheta}

\bigskip
\noindent
Theorem~\ref{surjection} shows that $\S_{E,R\op}(X)$ is isomorphic to a quotient of~$F_T$ and we want to understand the kernel.
We do this in the following result in terms of a system of equations.
The solution of this system of equations is quite hard and will only be obtained in~\cite{BT3},
when we will compute the dimension of the evaluations of fundamental functors and simple functors.

\result{Theorem} \label{kertheta} Let $T$ be a finite lattice, let $(E,R)$ be the full subposet of its irreducible elements, and let $X$ be a finite set.
The kernel of the map
$$\Theta_{T,X}:F_T(X)\to \S_{E,R\op}(X)$$
is equal to the set of linear combinations $\sum_{\varphi:X\to T}\limits\lambda_\varphi\varphi$, where $\lambda_\varphi\in k$, such that for any map $\psi:X\to \Iup(E,R)$
$$\sum_{\substack{\varphi:X\to T \\ \Gamma_\psi\op\Gamma_\varphi = R\op}} \lambda_\varphi=0\mpoint$$
Recall from Notation~\ref{Gamma} that $\Gamma_{\varphi}=\{(x,e)\in X\times E\mid e\leq\varphi(x)\}\subseteq \CC(X,E)$ and
$\Gamma_{\psi}=\{(x,e)\in X\times E\mid e\in\psi(x)\}\subseteq \CC(X,E)$.
\fresult

\pf The image of $\varphi:X\to T$ under the map $\Theta_{T,X}$ is equal to the class of $\Gamma_\varphi\otimes f_{R\op}\in L_{E,\CP_Ef_{R\op}}(X)$ in the quotient
$$\S_{E,R\op}(X) = L_{E,\CP_Ef_{R\op}}(X) / J_{E,\CP_Ef_{R\op}}(X) \mpoint$$
Therefore a linear combination $u=\sum_{\varphi:X\to T}\limits\lambda_\varphi\varphi$ lies in $\Ker\Theta_{T,X}$ if and only if
$\sum_{\varphi:X\to T}\limits\lambda_\varphi\Gamma_\varphi\otimes f_{R\op}$ belongs to $J_{E,\CP_Ef_{R\op}}(X)$.
We apply Lemma~\ref{belong-to-J}, using the fact that $\Gamma_\varphi R\op= \Gamma_\varphi$.
It follows that $u\in\Ker\Theta_{T,X}$ if and only if
$$\forall \, U\in R\op\CC(E,X) \,, \quad \sum_{\substack{\varphi:X\to T \\ U\Gamma_\varphi=R\op}} \lambda_\varphi=0 \mpoint$$
Now $U\in R\op\CC(E,X)$ if and only if $U\op\in\CC(X,E)R$ if and only if $U=\Gamma_\psi$ for some map $\psi:X\to \Iup(E,R)$
(by Lemma~\ref{rightR}). Thus the condition becomes
$$ \forall \,\psi:X\to \Iup(E,R) \,, \qquad
\sum_{\substack{\varphi:X\to T \\ \Gamma_\psi\op\Gamma_\varphi = R\op}} \lambda_\varphi=0 \mvirg$$
as was to be shown.
\endpf
 
The condition $\Gamma_\psi\op\Gamma_\varphi = R\op$ which appears in the system of equations in Theorem~\ref{kertheta} is the key for the description of the fundamental functor $\S_{E,R\op}$, and consequently for understanding the simple functors $S_{E,R,V}$.
We need to characterize this condition in various useful ways.\par

We first introduce the following notation.

\result{Notation} 
Let $T$ be a finite lattice, let $(E,R)$ be the full subposet of its irreducible elements, let
$\psi:X\to \Iup(E,R)$ be any map, and let $\varphi:X\to T$ be any map.
We define the function $\smash\psi:X\to T$ by 
$$\forall x\in X,\;\;\smash\psi(x)=\mbigsmash{e\in\psi(x)}e \mvirg$$
where $\bigwedge$ is the meet in the lattice~$T$.\par

Moreover, the notation $\varphi\leq\smash\psi$ means that $\varphi(x)\leq_T\smash\psi(x)$ for all $x\in X$.
\fresult

\medskip
We can now state the various characterizations of the condition which we need.

\result{Theorem} \label{Gammapsigamma}
Let $T$ be a finite lattice, let $(E,R)$ be the full subposet of its irreducible elements, and let $X$ be a finite set.
Let $\varphi:X \to T$ be a map and let
$$\Gamma_\varphi=\{(x,e)\in X\times E\mid e\leq_T \varphi(x)\}$$
be the associated correspondence.
Let $\psi:X\to \Iup(E,R)$ be a map and let
$$\Gamma_\psi\op=\big\{(e,x)\in E\times X\mid e\in\psi(x)\big\}$$
be the associated correspondence.
The following conditions are equivalent.
\begin{enumerate}
\item $\Gamma_\psi\op\varphi =\iota$.
\item $\Gamma_\psi\op\Gamma_\varphi \iota=\iota$.
\item $\Delta_E\subseteq\Gamma_\psi\op\Gamma_\varphi\subseteq R\op$.
\item $\Gamma_\psi\op\Gamma_\varphi= R\op$.
\item $\varphi\leq\smash\psi$ and
$\forall e\in E,\;\exists x\in X$ such that $\varphi(x)=e$ and $\psi(x)=[e,\cdot[_E$.
\item $\forall t\in T,\;\psi\big(\varphi^{-1}(t)\big)\subseteq[t,\cdot[_T\cap E$ and
$\;\forall e\in E,\;\psi\big(\varphi^{-1}(e)\big)=[e,\cdot[_E$.
\end{enumerate}
\fresult

\pf
(a) $\Leftrightarrow$ (b). By Lemma~\ref{surgamma}, we have $\varphi=\Gamma_\varphi \iota$.\mpn

(b) $\Leftrightarrow$ (c). We prove more generally that $Q\iota =\iota$ if and only if $\Delta_E\subseteq Q\subseteq R\op$,
for any relation $Q\in \CC(E,E)$.
Suppose first that $\Delta_E\subseteq Q\subseteq R\op$. Then for any $e\in E$, 
$$(Q\iota)(e)=\mbigvee{(e,f)\in Q}\iota(f)=\mbigvee{(e,f)\in Q}f\mpoint$$
Since $(e,e)\in Q$, it follows that $(Q\iota)(e)\geq_T e$.
On the other hand since $Q\subseteq R\op$, if $(e,f)\in Q$, then $f\leq_R e$, hence $f\leq_T e$.
Thus $(Q\iota)(e)\leq_T e$. It follows that $(Q\iota)(e)=e$ for any $e\in E$, i.e. $Q\iota=\iota$.\par

Conversely, if $Q\iota=\iota$, then 
$$\forall e\in E,\;\;\mbigvee{(e,f)\in Q}f=e\mpoint$$
As $e$ is irreducible, it follows that $(e,e)\in Q$ for any $e\in E$, i.e. $\Delta_E\subseteq Q$.
Moreover if $(e,f)\in Q$, then $f\leq_T e$, hence $f\leq_R e$ and so $Q\subseteq R\op$. This proves the claim and completes the proof of~(b).\mpn

(c) $\Leftrightarrow$ (d). If (c) holds, multiply on the right by~$R\op$ and use the equality
$\Gamma_\varphi R\op=\Gamma_\varphi$ of Lemma~\ref{rightR} to obtain~(d). On the other hand, it is clear that (d) implies~(c).\mpn

(d) $\Rightarrow$ (e). 
Suppose that $\Gamma_\psi\op\Gamma_\varphi= R\op$ and let $x\in X$. Then for all $f\leq_T \varphi(x)$ and for all $e\in \psi(x)$, we have $(e,x)\in \Gamma_\psi\op$ and $(x,f)\in \Gamma_\varphi$,
hence $(e,f)\in R\op$, that is, $f\leq_R e$, hence $f\leq_T e$. Therefore 
$\varphi(x)=\mbigvee{f\leq_T \varphi(x)}f \leq_T e$, whenever $e\in\psi(x)$. Thus 
$$\forall x\in X,\;\;\varphi(x)\leq_T \mbigsmash{e\in\psi(x)}e=\smash\psi(x)\mvirg$$
that is, $\varphi\leq\smash\psi$. This shows that the first property in~(e) holds.\par

Since $(e,e)\in R\op$, there exists $x_e\in X$ such that $e\leq_T\varphi(x_e)$ and
$e\in\psi(x_e)$. Then for all $f\leq_T \varphi(x_e)$, we have
$(e,x_e)\in \Gamma_\psi\op$ and $(x_e,f)\in \Gamma_\varphi$,
hence $(e,f)\in R\op$, that is, $f\leq_R e$, or in other words $f\leq_T e$. Thus again
$\varphi(x_e)=\mbigvee{f\leq_T \varphi(x_e)}f \leq_T e$, hence $\varphi(x_e)=e$.
Moreover, if $g\in E$ with $g\in\psi(x_e)$, then $(g,x_e)\in \Gamma_\psi\op$ and $(x_e,e)\in \Gamma_\varphi$, hence $(g,e)\in R\op$, that is, $e\leq_T g$. Therefore $\psi(x_e)\subseteq [e,\cdot[_E$.
But we also have $[e,\cdot[_E\subseteq \psi(x_e)$, as $e\in\psi(x_e)$ and $\psi(x_e)$ is an upper ideal of~$E$. Thus $\psi(x_e)= [e,\cdot[_E$. This shows that the second property in~(e) holds.\mpn

(e) $\Rightarrow$ (d). For any $e\in E$, there exists $x_e\in X$ such that $\varphi(x_e)=e$ and $\psi(x_e)=[e,\cdot[_E$.
If now $(f,e)\in R\op$, then $e\leq_R f$, hence $f\in\psi(x_e)$. Since we also have $e\leq_R \varphi(x_e)$, we obtain $(f,x_e)\in \Gamma_\psi\op$ and $(x_e,e)\in \Gamma_\varphi$.
Thus $R\op\subseteq \Gamma_\psi\op\Gamma_\varphi$.\par

Moreover if $(f,e)\in \Gamma_\psi\op\Gamma_\varphi$, then there exists $x\in X$ such that $f\in\psi(x)$ and $e\leq_T\varphi(x)$. Since $\varphi\leq\smash\psi$, we have $\varphi(x)\leq_T \mbigsmash{f\in\psi(x)}f$. It follows that $e\leq_T f$, hence
$e\leq_R f$, that is, $(f,e)\in R\op$. Thus $\Gamma_\psi\op\Gamma_\varphi\subseteq R\op$.
Therefore we obtain $\Gamma_\psi\op\Gamma_\varphi=R\op$.\mpn

(e) $\Leftrightarrow$ (f). 
We are going to slightly abuse notation by setting, for any subset $Y$ of $X$,
$\psi(Y)=\bigcup\limits_{x\in Y}\psi(x)$.
Taking $t=\varphi(x)$, the first condition in~(e) is equivalent to
$$\forall t\in T, \; e\in \psi\big(\varphi^{-1}(t)\big) \Longrightarrow t\leq_T e \mvirg$$
which in turn is equivalent to
$$\forall t\in T, \; \psi\big(\varphi^{-1}(t)\big)\subseteq[t,\cdot[_T\cap E \mpoint$$
In particular $\psi\big(\varphi^{-1}(e)\big)\subseteq [e,\cdot[_E$ for all $e\in E$ because $[e,\cdot[_T\cap E=[e,\cdot[_E$.
But the second condition in~(e) says that $e$ must belong to $\psi\big(\varphi^{-1}(e)\big)$, so we get 
$\psi\big(\varphi^{-1}(e)\big)=[e,\cdot[_E$.
This shows that the second condition in~(e) is equivalent to
$$\forall e\in E,\;\;\psi\big(\varphi^{-1}(e)\big)=[e,\cdot[_E\mpoint$$
This completes the proof of Theorem~\ref{Gammapsigamma}.\endpf

Condition (d) will play an important role in the proof of Theorem~\ref{SERgamma} below, while conditions (e) and~(f) will be main tools used in~\cite{BT3}.


\section{Duality and opposite lattices} \label{Section-duality-lattices}

\bigskip
\noindent
In this section, we prove a perfect duality between the functor associated to a lattice~$T$ and the functor associated to the opposite lattice~$T\op$.
We work with an arbitrary commutative base ring~$k$.\par

Let $F$ be a correspondence functor over~$k$.
The {\em dual} $F\dual$ of~$F$ is the correspondence functor defined on a finite set $X$ by
$$F\dual(X):=\Hom_k\big(F(X),k\big)\mpoint$$
If $Y$ is a finite set and $R\subseteq Y\times X$, then the map $F\dual(R): F\dual(X)\to F\dual(Y)$ is defined by
$$\forall \alpha\in F\dual(X),\;\;F\dual(R)(\alpha):=\alpha\circ F(R\op)\mpoint$$

Recall that $\CL$ denotes the category of finite lattices and $k\CL$ its $k$-linearization (Definition~\ref{category L}).
For any finite lattice $T=(T,\vee,\smash)$, denote by $T\op=(T,\smash,\vee)$ the {\em opposite} lattice, i.e. the set $T$ ordered with the opposite partial order.
For simplicity throughout this section, we write $\leq$ for $\leq_T$ and $\leq\op$ for $\leq_{T\op}$.

\result{Lemma} \label{Lop} The assignment $T\mapsto T\op$ extends to an isomorphism $\CL\to \CL\op$,
and to a $k$-linear isomorphism $k\CL\to k\CL\op$.
\fresult

\pf Let $f:T_1\to T_2$ be a morphism in the category $\CL$. For any $t\in T_2$, let $f\op(t)$ denote the join in $T_1$ of all the elements $x$ such that $f(x)\leq t$, i.e.
\begin{equation}\label{f-op}
f\op(t)=\mbigvee{f(x)\leq t}x\mpoint
\end{equation}
Then $f\big(f\op(t)\big)=\mbigvee{f(x)\leq t}f(x)\leq t$, so $f\op(t)$ is actually the greatest element of $f^{-1}\big([\hat{0},t]_{T_2}\big)$, i.e. $f^{-1}\big([\hat{0},t]_{T_2}\big)=[\hat{0},f\op(t)]_{T_1}$. In other words,
\begin{equation}\label{adjoint lattices}
\forall t_1\in T_1,\;\forall t_2\in T_2,\;\;f(t_1)\leq t_2 \iff t_1\leq f\op(t_2)\mvirg
\end{equation}
that is, the pair $(f,f\op)$ is an adjoint pair of functors between the posets $T_1$ and~$T_2$, viewed as categories. In those terms, saying that $f$ is a morphism in $\CL$ is equivalent to saying that $f$ commutes with colimits in $T_1$ and~$T_2$. Hence $f\op$ commutes with limits, that is, $f\op$ commutes with the meet operation, i.e.~it is a morphism of lattices $T_2\op\to T_1\op$.\par

In more elementary terms, for any subset $A\subseteq T_2$,
\begin{eqnarray*}
[\hat{0},f\op\big(\mbigsmash{t\in A}t\big)]_{T_1}&=&f^{-1}\big([\hat{0},\mbigsmash{t\in A}t]_{T_2}\big)\\
&=&f^{-1}\big(\bigcap_{t\in A}[\hat{0},t]_{T_2}\big)\\
&=&\bigcap_{t\in A}f^{-1}\big([\hat{0},t]_{T_2}\big)\\
&=&\bigcap_{t\in A}[\hat{0},f\op(t)]_{T_1}\\
&=&[\hat{0},\mbigsmash{t\in A}f\op(t)]_{T_1}\mpoint
\end{eqnarray*}
It follows that $f\op\big(\mbigsmash{t\in A}t\big)=\mbigsmash{t\in A}f\op(t)$, i.e. $f\op$ is a morphism $T_2\op\to T_1\op$ in~$\CL$.\par

Now denoting by $\leq\op$ the opposite order relations on both $T_1$ and $T_2$, Equation~\ref{adjoint lattices} reads
$$\forall t_2\in T_2,\;\forall t_1\in T_1,\;\;f\op(t_2)\leq\op t_1 \iff t_2\leq\op f(t_1)\mvirg$$
which shows that the same construction applied to the morphism $f\op:T_2\op\to T_1\op$ yields $(f\op)\op=f$.
This proves that the map $f\mapsto f\op$ is a bijection from $\Hom_\CL(T_1,T_2)$ to $\Hom_\CL(T_2\op,T_1\op)$.\par

Now if $f:T_1\to T_2$ and $g:T_2\to T_3$ are morphisms in $\CL$, the adjunction~(\ref{adjoint lattices}) easily implies that $(gf)\op=f\op g\op$.
It is clear moreover that $(\Id_T)\op=\Id_{T\op}$ for any finite lattice $T$.
Hence the assignment $T\mapsto T\op$ and $f\mapsto f\op$ is an isomorphism $\CL\to \CL\op$, which extends linearly to an isomorphism $k\CL\to k\CL\op$.
\endpf

\result{Definition} \label{bilinear}
Let $T$ be a finite lattice and let $X$ be a finite set. For two functions $\varphi:X\to T$ and $\psi:X\to T\op$, set
$$(\varphi,\psi)_X:=\left\{ \begin{array}{cl}
1& \text{ if } \; \varphi\leq \psi, \;\text{ i.e.~if } \; \varphi(x)\leq_T \psi(x), \;\forall x\in X, \\
0& \text{ otherwise.}\end{array} \right.$$
This definition extends uniquely to a $k$-bilinear form 
$$ ({-},{-})_X:F_T(X)\times F_{T\op}(X)\longrightarrow k\mpoint $$
This bilinear form induces a $k$-linear map $\Psi_{T,X}:F_{T\op}(X)\to(F_T)\dual(X)$ defined by
$\Psi_{T,X}(\psi)(\varphi)=(\varphi,\psi)_X$.
\fresult

\medskip
We need some notation.

\result{Notation} \label{notation-star}
Let $T$ be a finite lattice, $X$ and $Y$ finite sets, $Q\subseteq Y\times X$ a correspondence,
and $\psi:X\to T\op$ a map.
We denote by $Q\star\psi$ the action of the correspondence~$Q$ on~$\psi$.
In other words, $Q\star\psi$ is the map $F_{T\op}(Q)(\psi):Y\to T\op$. Recall that it is defined by
$$\forall y\in Y,\;\;(Q\star\psi)(y)=\mbigsmash{(y,x)\in Q}\psi(x)\mvirg$$
because the join in $T\op$ is the meet in~$T$.
\fresult

\result{Lemma}\label{pairing}
\begin{enumerate}
\item With the notation~\ref{notation-star}, the family of bilinear forms in Definition~\ref{bilinear} satisfy
$$(\varphi,Q\star\psi)_Y=(Q\op\varphi,\psi)_X \mpoint$$
\item The family of maps $\Psi_{T,X}:F_{T\op}(X)\to (F_T)\dual(X)$ form a morphism of correspondence functors $\Psi_T:F_{T\op}\to (F_T)\dual$.
\end{enumerate}
\fresult

\pf
(a) We have
\begin{eqnarray*}
\varphi\leq Q\star\psi&\iff&\forall y\in Y,\;\varphi(y)\leq_T Q\star\psi(y)\\
&\iff&\forall y\in Y,\;\varphi(y)\leq_T \mbigsmash{(y,x)\in Q}\psi(x)\\
&\iff&\forall (y,x)\in Q,\;\varphi(y)\leq_T \psi(x)\\
&\iff&\forall x\in X,\;\mbigvee{(x,y)\in Q\op}\varphi(y)\leq_T \psi(x)\\
&\iff&Q\op\varphi\leq \psi\mpoint
\end{eqnarray*}

(b) The equation in part (a) also reads
$$\Psi_{T,X}(\psi)(Q\op\varphi)=\Psi_{T,Y}(Q\star\psi)(\varphi)\mvirg$$
that is $Q\Psi_{T,X}(\psi)=\Psi_{T,Y}(Q\star\psi)$. 
\endpf

\begin{rem}{Remark} \label{same-forms} Let $T=\Idown(E,R\op)$ be the distributive lattice corresponding to a poset $(E,R\op)$.
Then $T\op=\Idown(E,R\op)\op$ is isomorphic, via complementation, to the lattice $\Idown(E,R)$. Using the isomorphisms of Proposition~\ref{C(-,E)R}
$$F_T=F_{\Idown(E,R\op)}\cong k\CC(-,E)R \mvirg \qquad F_{T\op}=F_{\Idown(E,R)}\cong k\CC(-,E)R\op \mvirg$$
we can transport the bilinear forms $(-,-)_X$ defined in~(\ref{bilinear}) and obtain a pairing
$$k\CC(-,E)R \times k\CC(-,E)R\op \longrightarrow k \mpoint$$
It is easy to check, using complementation, that this pairing coincides with the one obtained in~Remark~10.5 of~\cite{BT2}.
\end{rem}

\bigskip

\result{Notation} \label{varphi-star}
Let $T$ be a finite lattice, $X$ a finite set, and $\varphi:X\to T$ a map.
We denote by $\varphi^\star$ the element of $F_{T\op}(X)$ defined by
$$\varphi^\star:=\sum_{\substack{\rho:X\to T \\ \rho\leq \varphi}} \mu(\rho,\varphi)\rho^\circ\mvirg$$
where $\rho^\circ$ is the function~$\rho$, viewed as a map $X\to T\op$, and where $\mu(\rho,\varphi)$ is the M\"obius function of the poset of maps from $X$ to $T$, for which $\rho\leq \varphi$ if and only if $\rho(x)\leq \varphi(x)$ in $T$ for any $x\in X$.
Recall that $\mu(\rho,\varphi)$ can be computed as follows~:
$$\mu(\rho,\varphi)=\prod_{x\in X}\mu_T\big(\rho(x),\varphi(x)\big)\mvirg$$
where $\mu_T$ is the M\"obius function of the poset~$T$.
\fresult

\medskip
Now we can prove that we have a perfect duality.

\result{Theorem} \label{F_T dual} Let $T$ be a finite lattice.
\begin{enumerate}
\item Let $X$ be a finite set. The bilinear form (\ref{bilinear}) is nondegenerate, in the strong sense,
namely it induces an isomorphism
$$\Psi_{T,X}:F_{T\op}(X)\longrightarrow (F_T)\dual(X) \mpoint$$
More precisely, $\{\varphi^\star \mid \varphi:X\to T\}$ is the dual basis, in $F_{T\op}(X)$, of the $k$-basis of functions $X\to T$, in~$F_T(X)$.
\item $\Psi_T:F_{T\op}\to (F_T)\dual$ is an isomorphism of correspondence functors.
\item The functor $T\mapsto F_{T\op}$ and the functor $T\mapsto (F_T)\dual$ are naturally isomorphic functors from $k\CL$ to~$\CF_k\op$. More precisely, the family of isomorphisms $\Psi_T$, for finite lattices~$T$, form a natural transformation $\Psi$ between the functor $T\mapsto F_{T\op}$ and the functor $T\mapsto (F_T)\dual$.
\end{enumerate}
\fresult

\pf (a) The set $\{\rho^\circ \mid \rho^\circ:X\to T\op\}$ is a $k$-basis of the free $k$-module $F_{T\op}(X)$.
It follows that $\{\varphi^\star \mid \varphi:X\to T\}$ is also a $k$-basis of $F_{T\op}(X)$,
because the integral matrix of M\"obius coefficients $\mu(\rho,\phi)$ is unitriangular, hence invertible over~$\Z$.
Actually its inverse is the adjacency matrix of the order relation $\rho\leq \varphi$ on the set of maps $X\to T$.\par

Now, for any two functions $\varphi, \lambda:X\to T$,
\begin{equation} \label{dual basis}
(\lambda,\varphi^\star)_X=\sumb{\rho:X\to T}{\rho\leq \varphi}\mu(\rho,\varphi)(\lambda,\rho^\circ)_X
=\sumb{\rho:X\to T}{\lambda\leq\rho\leq \varphi}\mu(\rho,\varphi)
=\delta_{\lambda,\varphi}\mvirg
\end{equation}
where $\delta_{\lambda,\varphi}$ is the Kronecker symbol (the last equality coming from the definition of the M\"obius function). This shows that $\{\varphi^\star \mid \varphi:X\to T\}$ is the dual basis, in $F_{T\op}(X)$, of the $k$-basis of functions $X\to T$, in~$F_T(X)$.\mpn

(b) This follows immediately from~(a).
Another way of seeing this is to build an explicit inverse $\Phi_T$ of~$\Psi_T$.
For each finite set $X$, we define a linear map $\Phi_{T,X}:(F_T)\dual(X)\to F_{T\op}(X)$ by setting
$$\forall \alpha\in (F_T)\dual(X),\;\;\Phi_{T,X}(\alpha)=\sum_{\varphi:X\to T}\alpha(\varphi)\varphi^\star\mpoint$$
Then, for any function $\lambda:X\to T$,
$$\big(\Psi_{T,X}\Phi_{T,X}(\alpha)\big)(\lambda)=\big(\lambda,\Phi_{T,X}(\alpha)\big)_X=\sum_{\varphi:X\to T}\alpha(\varphi)(\lambda,\varphi^\star)_X=\alpha(\lambda)\mvirg$$
so $\Psi_{T,X}\Phi_{T,X}$ is the identity map of $(F_T)\dual(X)$.
In particular, $\Psi_{T,X}$ is surjective.\par

On the other hand, $\Psi_{T,X}$ is injective, because if $\Psi_{T,X}(\beta)=0$, then we write
$\beta=\sum\limits_{\varphi:X\to T}a_\varphi\varphi^\star$, where $a_\varphi\in k$, and
then for all~$\lambda:X\to T$, we get
$$0=\Psi_{T,X}\big(\sum_{\varphi:X\to T}a_\varphi\varphi^\star\big)(\lambda)
=(\lambda,\sum_{\varphi:X\to T}a_\varphi\varphi^\star)_X
=\sum_{\varphi:X\to T}a_\varphi (\lambda,\varphi^\star)_X=a_\lambda\mvirg$$
so that $\beta=0$.
Therefore $\Psi_{T,X}$ is an isomorphism and $\Phi_{T,X}$ is its inverse.\mpn

(c) Let $T'$ be another finite lattice, and let $\Psi_{T'}:F_{{T'}\op}\to (F_{T'})\dual$ be the corresponding morphism.
Let moreover $f:T\to T'$ be a morphism in~$\CL$. We claim that for any finite set~$X$, the square
$$\xymatrix{
F_{T\op}(X)\ar[r]^-{\Psi_{T,X}}&(F_T)\dual(X)\\
F_{{T'}\op}(X)\ar[r]_-{\Psi_{T',X}}\ar[u]^-{F_{f\op}}&(F_{T'})\dual(X)\ar[u]_-{(F_f)\dual}
}
$$
is commutative: indeed, for any functions $\psi:X\to {T'}\op$ and $\varphi:X\to T$,
$$\big((F_f)\dual\Psi_{T',X}(\psi)\big)(\varphi)=\Psi_{T',X}(\psi)(f\circ\varphi)=(f\circ\varphi,\psi)_X\mvirg$$
whereas
$$\big(\Psi_{T,X}F_{f\op}(\psi)\big)(\varphi)=\big(\varphi,F_{f\op}(\psi)\big)_X=(\varphi,f\op\circ\psi)_X\mpoint$$
Now by~\ref{adjoint lattices}, we have that
\begin{eqnarray*}
f\circ\varphi\leq\psi&\iff& \forall x\in X,\;\;f\big(\varphi(x)\big)\leq \psi(x)\\
&\iff&\forall x\in X,\;\;\varphi(x)\leq f\op\big(\psi(x)\big)\\
&\iff&\varphi\leq f\op\circ\psi\mvirg
\end{eqnarray*}
which proves our claim. This shows that the isomorphisms $\Psi_T$, for finite lattices~$T$, form a natural transformation $\Psi$ of the functor $T\mapsto F_{T\op}$ to the functor $T\mapsto (F_T)\dual$ from $k\CL$ to~$\CF_k\op$.
This completes the proof of Theorem~\ref{F_T dual}.
\endpf

\result{Corollary} \label{self-injective} Let $k$ be a self-injective ring.
Then for any distributive lattice~$T$, the functor $F_T$ is projective and injective in~$\CF_k$.
\fresult

\pf Since $T$ is distributive, the functor $F_T$ is projective by Theorem~\ref{projective-distributive}, without further assumption on~$k$.\par

If $k$ is self-injective, the functor sending a $k$-module $A$ to its $k$-dual $\Hom_k(A,k)$ is exact.
It follows that the functor $M\mapsto M\dual$ is an exact contravariant endofunctor of the category~$\CF_k$,
where $M\dual$ denotes the dual correspondence functor.\par

Let $\alpha:M\to N$ be an injective morphism in $\CF_k$, and let $\lambda:M\to F_T$ be any morphism.
Then $\alpha\dual:N\dual\to M\dual$ is surjective, and we have the following diagram with exact row in $\CF_k$
$$\xymatrix{
&(F_T)\dual\ar[d]^-{\lambda\dual}&\\
N\dual\ar[r]^-{\alpha\dual}&M\dual\ar[r]&0
}
$$
Now $(F_T)\dual\cong F_{T\op}$ by Theorem~\ref{F_T dual}, and $T\op$ is distributive. Hence $F_{T\op}$ is projective in $\CF_k$, and there exists a morphism $\beta:(F_T)\dual\to N\dual$ such that $\alpha\dual\circ\beta=\lambda\dual$. Dualizing once again the previous diagram yields the commutative diagram
$$\xymatrix{
&(F_T)\dual{}\dual&F_T\ar[l]_-{\eta_{F_T}}\\
N\dual{}\dual\ar[ur]^-{\beta\dual}&M\dual{}\dual\ar[u]_-{\lambda\dual{}\dual}\ar[l]^-{\alpha\dual{}\dual}
&\\
N\ar[u]^-{\eta_N}&M\ar[u]_-{\eta_M}\ar[l]^-{\alpha}\ar@/_/[uur]_-{\lambda}
}
$$
where for any functor $M$, we denote by $\eta_M$ the canonical morphism from $M$ to $M\dual{}\dual$. Now $\eta_{F_T}$ is an isomorphism, because for any finite set $X$, the module $F_T(X)$ is a finitely generated free $k$-module. Let $\varepsilon:N\to F_T$ be defined by $\varepsilon=\eta
_{F_T}^{-1}\circ\beta\dual\circ\eta_N$. Then
$$\varepsilon\circ\alpha=\eta_{F_T}^{-1}\circ\beta\dual\circ\eta_N\circ\alpha
=\eta_{F_T}^{-1}\circ\lambda\dual{}\dual\circ\eta_M=\eta_{F_T}^{-1}\circ\eta_{F_T}\circ\lambda=\lambda\mpoint$$
Thus for any injective morphism $\alpha:M\to N$ and any morphism $\lambda:M\to F_T$, there exists a morphism $\varepsilon:N\to F_T$ such that $\varepsilon\circ \alpha=\lambda$. Hence $F_T$ is injective in~$\CF_k$.
\endpf


\section{Duality for fundamental functors and simple functors} \label{Section-duality-simple}

\bigskip
\noindent
By Theorem~\ref{surjection}, any fundamental functor $\S_{E,R}$ is isomorphic to a quotient of some functor associated to a lattice.
One of the main purposes of this section is to use duality to realize $\S_{E,R}$ as a subfunctor of some functor associated to another lattice.
We also determine what is the dual of a simple functor.\par

We will study the subfunctor generated by a specific element of~$F_{T\op}(E)$ which will be defined below.
We need some more notation.

\result{Notation} \label{r(t)} Let $T$ be a finite lattice. If $t\in T$, let $r(t)$ denote the join of all the elements of $T$ strictly smaller than~$t$, i.e.
$$r(t):=\mbigvee{s<t}s$$
\fresult

Thus $r(t)=t$ if and only if $t$ is not irreducible.
If $t$ is irreducible, then $r(t)$ is the unique maximal element of~$[\hat0,t[$.

\result{Notation} \label{def-eta-gamma} Let $T$ be a finite lattice and let $(E,R)$ be the full subposet of its irreducible elements. If $A\subseteq E$, let $\eta_A:E\to T$ be the map defined by
$$\forall e\in E,\;\;\eta_A(e):=\left\{\begin{array}{ll}r(e)&\hbox{if}\;e\in A\\e&\hbox{if}\;e\notin A\end{array}\right.\mpoint$$
Moreover, let $\gamma_T$ denote the element of $F_{T\op}(E)$ defined by
$$\gamma_T:=\sum_{A\subseteq E}(-1)^{|A|}\eta_A^\circ\mvirg$$
where $\eta_A^\circ$ denotes the function $\eta_A$, viewed as a map $E\to T\op$. 
\fresult

We now show that this element $\gamma_T$ has another characterization.
Recall that we use a star $\star$, as in Notation~\ref{notation-star}, for the action of a correspondence on evaluations of $F_{T\op}$.

\result{Lemma} \label{gamma-iota} Let $T$ be a finite lattice, let $(E,R)$ be the full subposet of its irreducible elements, and let $\iota:E\to T$ be the inclusion map.
\begin{enumerate}
\item The element $\gamma_T$ is equal to $\iota^*$ (using Notation~\ref{varphi-star}).
\item $R\star\gamma_T=\gamma_T$.
\end{enumerate}
\fresult

\pf (a) By definition, $\iota^\star=\sum\limits_{\rho\leq\iota}\mu(\rho,\iota)\rho^\circ$,
where $\rho^\circ$ denotes the function~$\rho$, viewed as a map $E\to T\op$,
and where $\mu$ is the M\"obius function of the poset of functions from $E$ to~$T$ (see Notation~\ref{varphi-star}).
Furthermore
$$\mu(\rho,\iota)=\prod_{e\in E}\mu_T\big(\rho(e),\iota(e)\big)\mvirg$$
where $\mu_T$ is the M\"obius function of the poset~$T$.
Now $\mu_T\big(\rho(e),\iota(e)\big)=\mu_T\big(\rho(e),e\big)$ is equal to 0 if $\rho(e)<r(e)$, because in that case the interval $]\rho(e),e[_T$ has a greatest element~$r(e)$. Moreover $\mu_T\big(\rho(e),e\big)$ is equal to -1 if $\rho(e)=r(e)$, and to +1 if $\rho(e)=e$. It follows that the only maps $\rho$ appearing in the sum above are of the form $\rho=\eta_A$ for some subset $A\subseteq E$ and $\mu(\eta_A,\iota)=(-1)^{|A|}$. Therefore
$$\iota^\star=\sum_{A\subseteq E} (-1)^{|A|} \eta_A^\circ=\gamma_T \mpoint$$

(b) For any $A\subseteq E$ and any $e\in E$,
$$(R*\eta_A)(e)=\mbigsmash{(e,e')\in R}\eta_A(e')=\mbigsmash{e\leq e'}\eta_A(e')=\eta_A(e)\mvirg$$
since $e<e'$ implies $\eta_A(e)\leq e\leq r(e')\leq\eta_A(e')$. Therefore $R\star \gamma_T=\gamma_T$.
\endpf

\bigskip   
Our aim is to show that the subfunctor $\langle\gamma_T\rangle$ of $F_{T\op}$ generated by $\gamma_T$ is isomorphic to the fundamental correspondence functor $\S_{E,R}$. We first show that $\langle\gamma_T\rangle$ is independent of the choice of~$T$.

\result{Lemma} \label{iso gamma_T}   
Let $f:T\to T'$ be a morphism in $\CL$ and let $(E,R)$, respectively $(E',R')$, be the full subposet of irreducible elements in~$T$, respectively~$T'$.
Suppose that the restriction of $f$ to $E$ is an isomorphism of posets $f_|:(E,R)\stackrel{\cong}{\longrightarrow}(E',R')$. 
\begin{enumerate}
\item The map $f$ is surjective and $fr(e)=rf(e)$ for any $e\in E$.
\item The map $f\op:{T'}\op\to T\op$ restricts to a bijection $f\op_|:E' \stackrel{\cong}{\longrightarrow}E$, which is inverse to~$f_|$. Moreover $f\op r(e')=rf\op(e')$ for any $e'\in E'$.
\item $f\op:{T'}\op\to T\op$ induces an injective morphism $F_{f\op}:F_{{T'}\op}\to F_{T\op}$ and an isomorphism $\langle\gamma_{T'}\rangle\cong\langle\gamma_T\rangle$.
\end{enumerate}
\fresult

\pf (a) Since any element of $T'$ is a join of irreducible elements, which are in the image of~$f$, and since $f$ commutes with joins, the map $f$ is surjective.\par

Let $e\in E$. By assumption $f(e)\in  E'$. The condition $r(e)<e$ implies $f\big(r(e)\big)\leq f(e)$. Moreover $r(e)=\mbigveeb{e_1\in E}{e_1<e}e_1$, hence $fr(e)=\mbigveeb{e_1\in E}{e_1<e}f(e_1)$. Thus if $fr(e)=f(e)$, then there exists $e_1<e$ such that $f(e_1)=f(e)$, contradicting the assumption on~$f$. It follows that $fr(e)\leq rf(e)$.\par

Now $rf(e)=\mbigveeb{e'\in E'}{e'<f(e)}e'$, and each $e'\in E'$ with $e'<f(e)$ can be written $e'=f(e_1)$, for $e_1\in E$ with $e_1<e$. It follows that $rf(e)\leq\mbigveeb{e_1\in E}{e_1<e}f(e_1)=fr(e)$. Thus $rf(e)=fr(e)$, as was to be shown.\mpn

(b) Recall from Equation~\ref{f-op} that $f\op$ is defined by $f\op(t')=\mbigvee{f(t)\leq t'}t$. Let $e'\in E'$. Then there exists $e\in E$ such that $f(e)=e'$. Let $t\in T$ be such that $f(t)\leq e'$ and write $t=\mbigveeb{e_1\in E}{e_1\leq t}e_1$. For each $e_1\in  E$ with $e_1\leq t$, we have $f(e_1)\leq f(t)\leq e'=f(e)$, hence 
$e_1\leq e$, and $t\leq e$. It follows that $f\op(e')=\mbigvee{f(t)\leq e'}t=e$, so $f\op_|$ is a bijection $E'\to E$, inverse to~$f_|$. This proves the first statement in~(b).\par

Now let $e\in E$, and set $e'=f(e)\in E'$. First we have $r(e')\leq e'$, thus $f\op r(e')\leq f\op(e')=e$. If $f\op r(e')=e$, then $\mbigvee{f(t)\leq r(e')}t=e$, hence 
$$f(e)\leq r(e')<e'=f(e)\mvirg$$ 
a contradiction. Thus $f\op r(e')\leq r(e)=rf\op(e')$. But we also have 
$$rf\op(e')=r(e)=\mbigveeb{e_1\in E}{e_1<e}e_1=\mbigveeb{e_1\in E}{f(e_1)<f(e)}e_1\leq\mbigveeb{t\in T}{f(t)\leq rf(e)}t=f\op rf(e)=f\op r(e')\mvirg$$
so $f\op r(e')=rf\op(e')$, which proves the second statement in~(b).\mpn

(c) Since $f$ is surjective by~(a), so is the morphism $F_f:F_T\to F_{T'}$. By duality and Theorem~\ref{F_T dual}, the morphism $F_{f\op}:F_{{T'}\op}\to F_{T\op}$ can be identified with the dual of $F_f$ and is therefore injective.
This proves the first statement in~(c).\par

Now for any $B\subseteq  E'$, consider the map $\eta_B^\circ: E' \to {T'}\op$. Then for any $e'\in  E'$
$$f\op \eta_B^\circ(e')=\left\{\begin{array}{ll}f\op(e')&\hbox{if}\;e'\notin B,\\ f\op r(e')=rf\op(e')&\hbox{if}\;e'\in B.\end{array}\right.$$
Hence $f\op\circ\eta_B^\circ=\eta_{f\op_|(B)}^\circ\circ f\op_|$, and therefore $f\op\circ\gamma_{T'}=\gamma_{T}\circ f\op_|$. It follows that
$$F_{f\op}(\gamma_{T'}\circ f_|)=f\op\circ\gamma_{T'}\circ f_|=\gamma_{T}\circ f\op_|\circ f_|=\gamma_{T} \mpoint$$
    
Therefore the injective morphism $F_{f\op}$ maps the subfunctor $\langle\,\gamma_{T'}\circ f_| \,\rangle$ isomorphically to the subfunctor $\langle\gamma_T\rangle$. But since $f_|:E\to E'$ is a bijection, the subfunctor $\langle\,\gamma_{T'}\circ f_| \,\rangle$ of $F_{{T'}\op}$ is equal to the subfunctor $\langle\gamma_{T'}\rangle$. This proves the second statement in~(c).
\endpf

Recall that we use a star $\star$, as in Notation~\ref{notation-star}, for the action of a correspondence on evaluations of $F_{T\op}$. We now come to our main result.

\result{Theorem} \label{SERgamma} Let $T$ be a finite lattice and let $(E,R)$ be the full subposet of its irreducible elements.
\begin{enumerate}
\item The subfunctor $\langle\gamma_T\rangle$ of $F_{T\op}$ generated by~$\gamma_T$ is isomorphic to $\S_{E,R}$.
\item In other words, for any finite set $X$, the module $\S_{E,R}(X)$ is isomorphic to the $k$-submodule of $F_{T\op}(X)$ generated by the elements $S\star\gamma_T$, for ${S\subseteq X\times E}$.
\end{enumerate}
\fresult

\pf
We first show that it suffices to prove the result in the case when $T$ is the lattice $\Idown(E,R)$.
For any other lattice $T'$ with the same poset $(E,R)$ of irreducible elements,
the inclusion $E\subseteq T'$ extends to a unique surjective map
$$f:T=\Idown(E,R)\longrightarrow T'$$
in the category~$\CL$ which induces the identity on $E$ (see Lemma~\ref{principal}).
Then $\langle\gamma_{T'}\rangle$ is isomorphic to $\langle\gamma_T\rangle$ by Lemma~\ref{iso gamma_T}, so we now assume that $T=\Idown(E,R)$, which is a distributive lattice.\par

We now apply Proposition~\ref{HomFT} to the element $\gamma_T \in F_{T\op}(E)$,
using the fact that $\gamma_T\in R\star F_{T\op}(E)$,
because $R\star\gamma_T=\gamma_T$ by Lemma~\ref{gamma-iota}.
We deduce that $\gamma_T \in F_{T\op}(E)$ determines a unique morphism 
$$\xi:F_{\Idown(E,R\op)} \to F_{T\op}$$
such that $\xi_E(j)=\gamma_T$, where $j:E\to \Idown(E,R\op)$ is the inclusion map.\par

Then for any finite set~$X$ and any function $\psi:X\to \Idown(E,R\op)$,
we can compute the map $\xi_X: F_{\Idown(E,R\op)}(X)\to F_{T\op}(X)$ as follows~:
$$\xi_X(\psi)=\xi_X(\Gamma_\psi j)=\Gamma_\psi\star\xi_E(j)=\Gamma_\psi\star \gamma_T \mvirg$$
using the equality $\psi=\Gamma_\psi j$ of Lemma~\ref{surgamma}.
In particular, the image of~$\xi$ is the subfunctor $\langle\gamma_T\rangle$ generated by~$\gamma_T$.

By Theorem~\ref{surjection} and the fact that the lattice $\Idown(E,R\op)$ has $(E,R\op)$ as its full subset of irreducible elements,
there is a surjective morphism
$$\Theta_{\Idown(E,R\op)} : F_{\Idown(E,R\op)} \longrightarrow \S_{E,R}$$
and its kernel is described in Theorem~\ref{kertheta}.
We want to prove that, for any finite set $X$, the kernel of the surjection
$$\xi_X: F_{\Idown(E,R\op)}(X)\longrightarrow \langle \gamma_T\rangle(X)$$
is equal to the kernel of the surjection
$$\Theta_{\Idown(E,R\op),X}: F_{\Idown(E,R\op)}(X) \longrightarrow \S_{E,R}(X) \mvirg$$
from which the isomorphism $\langle\gamma_T\rangle\cong \S_{E,R}$ will follow.\par

The kernel of the surjection $\xi_X$ is the set of all linear combinations
$$u=\sum_{\psi:X\to \Idown(E,R\op)}\limits\lambda_\psi \psi \mvirg$$
where $\lambda_\psi\in k$, such that 
$$\sum_{\psi:X\to\Idown(E,R\op)}\limits\lambda_\psi \Gamma_\psi\star\gamma_T=0\mpoint$$
Equivalently, using the nondegeneracy of the bilinear form $(-,-)_X$ of~(\ref{bilinear}), proved in Theorem~\ref{F_T dual},
$$\forall \,\varphi: X\to T,\;\;(\varphi,\sum_{\psi}\lambda_\psi \Gamma_\psi\star\gamma_T)_X =0 \mpoint$$
By Equation~\ref{pairing}, this is in turn equivalent to
$$\forall \,\varphi: X\to T,\;\;\sum_{\psi}\lambda_\psi (\Gamma_\psi\op\varphi,\gamma_T)_E=0\mpoint$$
Now $\gamma_T=\iota^\star$ by Lemma~\ref{gamma-iota} and we use $\iota^\star$ instead.
By Equation~\ref{dual basis}, we have
$$(\Gamma_\psi\op\varphi,\gamma_T)_E=(\Gamma_\psi\op\varphi,\iota^\star)_E=\delta_{\Gamma_\psi\op\varphi,\iota}$$
and therefore we obtain the condition
$$\forall \, \varphi: X\to T,\;\;
\sum_{\substack{\psi:X\to\Idown(E,R\op) \\ \Gamma_\psi\op\varphi=\iota}}\lambda_\psi=0\mpoint$$
Since we have assumed that $T=\Idown(E,R)$ and since $\Idown(E,R\op)=\Iup(E,R)$,
we have maps $\varphi:X\to \Idown(E,R)$ and $\psi:X\to\Iup(E,R)$
and we know from Theorem~\ref{Gammapsigamma} that the property
$\Gamma_\psi\op\varphi=\iota$ is equivalent to $\Gamma_\psi\op\Gamma_\varphi= R\op$.\par

It follows that $u=\sum_{\psi:X\to \Iup(E,R)}\limits\lambda_\psi \psi$ is in the kernel of $\xi_X$ if and only if
$$\forall \, \varphi:X\to \Idown(E,R),\;\;
\sum_{\substack{\psi:X\to \Idown(E,R\op) \\ \Gamma_\psi\op\Gamma_\varphi= R\op}} \lambda_\psi=0\mpoint$$
But the condition $\Gamma_\psi\op\Gamma_\varphi= R\op$ is in turn is equivalent to $\Gamma_\varphi\op\Gamma_\psi= R$, by passing to the opposite.
Moreover $\Idown(E,R)=\Iup(E,R\op)$ and $\Iup(E,R)=\Idown(E,R\op)$.
Therefore $u=\sum_{\psi:X\to\Idown(E,R\op)}\limits\lambda_\psi \psi$ is in the kernel of $\xi_{X}$ if and only if 
$$\forall \varphi:X\to \Iup(E,R\op),\;\;
\sum_{\substack{\psi:X\to \Idown(E,R\op) \\ \Gamma_\varphi\op\Gamma_\psi= R}} \lambda_\psi=0 \mpoint$$
By Theorem~\ref{kertheta}, this is equivalent to requiring that $u\in\Ker\Theta_{\Idown(E,R\op),X}$.
It follows that $\Ker\xi_X=\Ker\Theta_{\Idown(E,R\op),X}$.
Consequently, the images of $\xi_X$ and $\Theta_{\Idown(E,R\op),X}$ are isomorphic, that is,
$\langle\gamma_T\rangle\cong\S_{E,R}$.
This completes the proof of Theorem~\ref{SERgamma}.
\endpf

Since we now know that the subfunctor $\langle\gamma_T\rangle$ of $F_{T\op}$ is isomorphic to $\S_{E,R}$, we use again duality to obtain more.

\result{Theorem} \label{fundamental-dual}
Let $T$ be a finite lattice and let $(E,R)$ be the full subposet of its irreducible elements.
We consider orthogonal $k$-submodules with respect to the family of $k$-bilinear forms $(-,-)_X$ defined in~(\ref{bilinear}).
\begin{enumerate}
\item $\langle\gamma_T\rangle^\perp=\Ker\Theta_T$, where $\Theta_T:F_T\to \S_{E,R\op}$ is the morphism of Theorem~\ref{surjection}.
\item $F_T/\langle\gamma_T\rangle^\perp \cong \S_{E,R\op}$.
\item $\langle\gamma_T\rangle^{\perp\perp}\cong \S_{E,R\op}\dual$.
\item There is a canonical injective morphism $\alpha_{E,R}:\S_{E,R}\longrightarrow\S_{E,R\op}\dual$.
\end{enumerate}
\fresult

\pf (a) Let $\sum\limits_{\varphi:X\to T} \lambda_\varphi\varphi \in F_T(X)$, where $X$ is a finite set. Then
$$
\begin{array}{rl}
&\sum\limits_{\varphi} \lambda_\varphi\varphi \in \langle\gamma_T\rangle(X)^\perp
\iff
\big( \sum\limits_{\varphi} \lambda_\varphi\varphi, Q\star \gamma_T \big)_X=0 \quad\forall Q\in\CC(X,E) \\
&\iff
\big( \sum\limits_{\varphi} \lambda_\varphi\varphi, Q\star \gamma_T \big)_X=0 \quad\forall Q\in\CC(X,E)R
\qquad (\text{because } R\star\gamma_T=\gamma_T)  \\
&\iff
\sum\limits_{\varphi} \lambda_\varphi \big(Q\op\varphi, \gamma_T \big)_E=0 \quad\forall Q\in\CC(X,E)R
\qquad (\text{by \ref{pairing}})  \\
&\iff
\sumb{\varphi}{Q\op\varphi=\iota}\limits \lambda_\varphi =0 \quad\forall Q\in\CC(X,E)R
\qquad (\text{by \ref{dual basis} and Lemma~\ref{gamma-iota}})  \\
&\iff
\sumb{\varphi}{\Gamma_\psi\op\varphi=\iota}\limits \lambda_\varphi =0 \quad\forall \psi:X\to \Iup(E,R)
\qquad (\text{by Proposition~\ref{C(-,E)R}})  \\
&\iff
\sumb{\varphi}{\Gamma_\psi\op\Gamma_\varphi=R\op}\limits \lambda_\varphi =0 \quad\forall \psi:X\to \Iup(E,R)
\qquad (\text{by Theorem~\ref{Gammapsigamma}})  \\
&\iff
\sum\limits_{\varphi:X\to T} \lambda_\varphi\varphi \in \Ker\Theta_{T,X}
\quad (\text{by Theorem~\ref{kertheta}}) 
\end{array}$$
Therefore $\langle\gamma_T\rangle(X)^\perp=\Ker\Theta_{T,X}$.\mpn

(b) This follows immediately from (a) and Theorem~\ref{surjection}.\mpn

(c) This follows immediately from (b) and duality.\mpn

(d) There is an obvious inclusion $\langle\gamma_T\rangle \subseteq \langle\gamma_T\rangle^{\perp\perp}$.
Now we have $\langle\gamma_T\rangle\cong \S_{E,R}$ by Theorem~\ref{SERgamma}
and $\langle\gamma_T\rangle^{\perp\perp}\cong \S_{E,R\op}\dual$ by~(c).
Thus we obtain a canonical injective morphism $\S_{E,R}\longrightarrow\S_{E,R\op}\dual$.
\endpf

\begin{rem}{Remark} We will prove in a future paper that $\alpha_{E,R}:\S_{E,R}\longrightarrow\S_{E,R\op}\dual$ is actually an isomorphism.
This is easy to prove if $k$ is a field, because the inclusion $\langle\gamma_T\rangle \subseteq \langle\gamma_T\rangle^{\perp\perp}$ must be an equality since the pairing~(\ref{bilinear}) is nondegenerate, by Theorem~\ref{F_T dual}.
\end{rem}

\bigskip

We end this section with a description of the dual of a simple functor.
We assume that $k$ is a field and we let $S_{E,R,V}$ be the simple correspondence functor (over~$k$) parametrized by~$(E,R,V)$.
Part~(d) of Theorem~\ref{fundamental-dual} suggests that the index~$R$ must become $R\op$ after applying duality.
We now show that this is indeed the case. 

\result{Theorem} \label{simple-dual}  
Let $k$ be a field. The dual $S_{E,R,V}\dual$ of the simple functor $S_{E,R,V}$ is isomorphic to $S_{E,R\op,V\dual}$, where $V\dual$ denotes the ordinary dual of the $k\Aut(E,R)$-module~$V$.
\fresult

\pf
For simplicity, write $\CR_E=k\CC(E,E)$ for the algebra of all relations on~$E$
and recall from Section~\ref{Section-functors} that $\CP_E$ is a quotient algebra of~$\CR_E$.
The evaluation $S_{E,R,V}(E)$ is the $\CR_E$-module
$$T_{R,V}=\CP_Ef_R\otimes_{k\Aut(E,R)}V \mpoint$$
(which is actually a simple $\CR_E$-module by Proposition~\ref{simple-modules}).\par

Clearly the dual $S_{E,R,V}\dual$ is again a simple functor and its minimal set is $E$ again.
Moreover, by evaluation at~$E$, we find that
$$S_{E,R,V}\dual(E)\cong S_{E,R,V}(E)\dual = T_{R,V}\dual \mpoint$$
Here the action of a relation $Q\in \CR_E$ on a $\CR_E$-module $W\dual$ is defined by
$$(Q\cdot \alpha)(w) = \alpha(Q\op \cdot w)\mvirg \;\forall \alpha\in W\dual \mvirg \; \forall w\in W \mpoint$$
We are going to define a nondegenerate pairing
$$\langle -,-\rangle :T_{R,V} \times T_{R\op,V\dual} \longrightarrow k \mvirg$$
satisfying $\langle Q\cdot x \,,\,y\rangle = \langle x\,,\,Q\op\cdot y\rangle$ for all $x\in T_{R,V}$, $y\in T_{R\op,V\dual}$, and $Q\in \CR_E$.
This will induce an isomorphism of $\CR_E$-modules
$$T_{R,V}\dual \cong T_{R\op,V\dual} \mpoint$$
But a simple functor with minimal set~$E$ is completely determined by its evaluation at~$E$,
because it is generated by this minimal nonzero evaluation (by simplicity). 
Since $S_{E,R,V}\dual$ and $S_{E,R\op,V\dual}$ have both an evaluation at~$E$ isomorphic to~$T_{R\op,V\dual}$, it follows that
$$S_{E,R,V}\dual \cong S_{E,R\op,V\dual} \mvirg$$
as required.\par

Now we construct the required pairing.
By Proposition~\ref{fundamental-module}, $\CP_Ef_R$ has a $k$-basis $\{\Delta_\sigma f_R \mid \sigma\in \Sigma_E\}$,
where $\Sigma_E$ is the group of all permutations of~$E$.
Moreover, it is a free right $k\Aut(E,R)$-module and it follows that we can write
$$T_{R,V}=\CP_Ef_R\otimes_{k\Aut(E,R)}V
=\bigoplusl_{\sigma\in[\Sigma_E/\Aut(E,R)]} \Delta_\sigma f_R \otimes V\mvirg$$
where $[\Sigma_E/\Aut(E,R)]$ denotes a set of representatives of the left cosets of $\Aut(E,R)$ in~$\Sigma_E$.
Noticing that $\Aut(E,R\op)=\Aut(E,R)$, we have a similar decomposition
$$T_{R\op,V\dual}=\CP_Ef_{R\op}\otimes_{k\Aut(E,R)}V\dual
=\bigoplusl_{\varepsilon\in[\Sigma_E/\Aut(E,R)]} \Delta_\varepsilon f_{R\op} \otimes V\dual\mpoint$$
We define the pairing
$$\langle -,-\rangle :T_{R,V} \times T_{R\op,V\dual} \longrightarrow k \mvirg \qquad
\langle \Delta_\sigma f_R\otimes v\,,\,\Delta_\varepsilon f_{R\op} \otimes \alpha\rangle
:= \delta_{\sigma,\varepsilon}\, \alpha(v) \mvirg$$
where $\sigma, \varepsilon \in [\Sigma_E/\Aut(E,R)]$, $v\in V$, $\alpha\in V\dual$.\par

By choosing dual bases of $V$ and $V\dual$, we easily find dual bases of $T_{R,V}$ and $T_{R\op,V\dual}$, and it follows that this pairing is nondegenerate. We are left with the proof of the required property of this pairing, namely 
$$\langle Q\cdot x \,,\,y\rangle = \langle x\,,\,Q\op\cdot y\rangle\qquad \text{for all } \;x\in T_{R,V},\; y\in T_{R\op,V\dual}, \; Q\in \CC(E,E) \mpoint $$
By Proposition~\ref{fundamental-module}, the action of $Q$ on~$\CP_E f_R$ is given by:
$$Q\cdot \Delta_\sigma f_R=\left\{\begin{array}{ll}
\Delta_{\tau\sigma}f_R&\hbox{if}\;\;\exists\tau\in\Sigma_E\;\hbox{such that}\;
\Delta_E\subseteq \Delta_{\tau^{-1}}Q\subseteq {\ls\sigma R},\\
0&\hbox{otherwise}\mpoint \end{array}\right.$$
It follows that
$$\langle Q\cdot\Delta_\sigma f_R\otimes v\,,\,\Delta_\varepsilon f_{R\op} \otimes \alpha\rangle
=\left\{ \begin{array}{cl} \alpha(v)  \quad &\text{ if } \; \varepsilon=\tau\sigma \;\text{ and }\;
\Delta_E\subseteq \Delta_{\tau^{-1}}Q\subseteq {\ls\sigma R} \mvirg \\
0\quad &\text{ otherwise} \mpoint \end{array}\right.$$
On the other hand
$$\langle \Delta_\sigma f_R\otimes v\,,\,Q\op\cdot\Delta_\varepsilon f_{R\op} \otimes \alpha\rangle
=\left\{ \begin{array}{cl} \alpha(v)  \quad &\text{ if } \; \rho\varepsilon=\sigma \;\text{ and }\;
\Delta_E\subseteq \Delta_{\rho^{-1}}Q\op\subseteq  \ls\varepsilon R\op \mvirg \\
0\quad &\text{ otherwise} \mpoint \end{array}\right.$$
We put $\rho=\tau^{-1}$ and we assume that $\varepsilon=\tau\sigma$ (that is, $\rho\varepsilon=\sigma$).
We then obtain equivalent conditions:
$$\begin{array}{rcl}
\Delta_E\subseteq \Delta_{\tau}Q\op\subseteq  \ls\varepsilon R\op
&\iff & \Delta_E\subseteq Q\Delta_{\tau^{-1}}\subseteq  \ls\varepsilon R \qquad\text{(taking the opposite)} \\
&\iff & \Delta_E\subseteq Q\Delta_{\tau^{-1}}\subseteq  \ls{\tau\sigma} R \qquad\text{(because $\varepsilon=\tau\sigma$)} \\
&\iff & \Delta_E\subseteq \Delta_{\tau^{-1}}Q \subseteq  \ls{\sigma} R \qquad\text{(conjugating by $\Delta_{\tau^{-1}}$)}
\end{array}$$
Therefore, we obtain the required equality
$$\langle Q\cdot\Delta_\sigma f_R\otimes v\,,\,\Delta_\varepsilon f_{R\op} \otimes \alpha\rangle
=\langle \Delta_\sigma f_R\otimes v\,,\,Q\op\cdot\Delta_\varepsilon f_{R\op} \otimes \alpha\rangle \mvirg$$
from which it follows that we have an isomorphism of $\CR_E$-modules
$T_{R,V}\dual \cong T_{R\op,V\dual}$.
This completes the proof.
\endpf


\section{Embeddings and idempotents corresponding to total orders} \label{Section-embed-total}

\bigskip
\noindent
In this section, we construct morphisms of lattices in the category~$k\CL$ between a finite lattice~$T$ and a totally ordered lattice.
This will be used in Section~\ref{Section-total} to obtain detailed information about correspondence functors associated to a total order.\par

For a (non negative) integer $n\in\N$, we denote by $\sou{n}$ the set $\{0,1,\ldots,n\}$, linearly ordered by $0<1<\ldots<n$.
Then $\sou{n}$ is a distributive lattice, with least element 0 and greatest element~$n$.
Moreover $x\vee y=\sup(x,y)$ and $x\smash y=\inf(x,y)$, for any $x,y\in\sou{n}$.
We denote by $[n]=\{1,\ldots,n\}$ the set of irreducible elements of $\sou{n}$, viewed as a full subposet of~$\sou{n}$.\par

Let $\pi:T\to \sou n$ be a surjective morphism of finite lattices.
For every $h\in \sou n$, let
$$b_h=\sup\big(\pi^{-1}(h)\big) = \bigvee_{\pi(t)=h} t = \bigvee_{\pi(t)\leq h} t \mpoint$$
In other words, $b_h=\pi\op(h)$ with the notation of Lemma~\ref{Lop}.
Then we have $\widehat 0 \leq b_0 < b_1 < \ldots < b_{n-1} < b_n=\widehat 1$.
If $h\in\Irr(\sou n)=[n]=\{1,\ldots,n\}$, the whole interval $]b_{h-1},b_h]_T$ is mapped to~$h$ under~$\pi$,
while $b_{h-1}$ is mapped to~$h-1$.
The morphism $\pi$ determines uniquely the totally ordered $n$-tuple $B=(b_0,b_1,\ldots,b_{n-1})$.\par

Conversely, if $B=(b_0,b_1,\ldots,b_{n-1})$ is a strictly increasing $n$-tuple in~$T-\{\widehat1\}$ and if $b_n=\widehat1$
(so that $\widehat 0 \leq b_0 < b_1 < \ldots < b_{n-1} < b_n=\widehat 1$),
then $B$ determines a unique surjective morphism $\pi^B:T\to \sou n$
by setting $\pi^B(t)=h$ if $t\leq b_h$ and $t\not\leq b_{h-1}$
(hence in particular $\pi^B(\widehat 0)=0$).\par

For every $h\in [n]$, choose $a_h\in [b_{h-1},b_h]_T$, and define the $n$-tuple
$$A=(a_1,a_2,\ldots,a_n) \mpoint$$
Set also $a_0=\widehat 0$. Notice that $A$ is totally ordered~: $\widehat 0=a_0\leq a_1\leq \ldots \leq a_n\leq \widehat1$, with possible equalities.
Define
$$j_A^B:\sou n \longrightarrow T \,, \qquad h\mapsto a_h \mpoint$$
It is easy to see that $j_A^B$ is a morphism of lattices (because if $e<f$ in~$\sou n$,
then $e\leq f-1$, hence $a_e\leq b_e\leq b_{f-1}\leq a_f$).
Note that $j_A$ is a section of~$\pi^B$ if $a_h\in ]b_{h-1},b_h]_T$ for every $h\in [n]$, but not if $a_h=b_{h-1}$ for some~$h$.

Write
$$\mu(B,A)=\prod_{h\in [n]} \mu_T(b_{h-1},a_h) \mvirg$$
where $\mu_T(-,-)$ denotes the M\"obius function of the lattice~$T$.
For simplicity, we write $\mu=\mu_T$ throughout this section and the next.
Allowing the $n$-tuple $A$ to vary (i.e. $a_h$ varies in $[b_{h-1},b_h]_T$ for each~$h$), define
$$j^B=(-1)^n\sum_{A}\mu(B,A) j_A^B \in k\CL(\sou n,T) \mvirg$$
or in other words
$$j^B=(-1)^n \prod_{h=1}^n \Big(\sum_{a_h\in [b_{h-1},b_h]_T} \mu(b_{h-1},a_h)\Big) \,j_{(a_1,a_2,\ldots,a_n)}^B \mpoint$$
 
\result{Proposition} \label{vanish-on-H}
Let $T$ be a finite lattice and let $B$ be a strictly increasing $n$-tuple in $T-\{\widehat 1\}$.
Let $j^B:\sou n\to T$ be the corresponding morphism, constructed as above.
\begin{enumerate}
\item For any finite set $X$ and any function $\varphi:X\to \sou n$ such that $[n]\not\subseteq \varphi(X)$,
we have $j^B\varphi=0$.
\item $j^B$ induces $F_{j^B}:F_{\sou n}\to F_T$ vanishing on $H_{\sou n}$, hence induces in turn a morphism
$$\overline F_{j^B}: F_{\sou n}/H_{\sou n} \longrightarrow F_T \mpoint$$
\end{enumerate}
\fresult

\pf (a) We have
$$j^B\varphi=(-1)^n\sum_A\mu(B,A) j_A^B\varphi
=\sum_{\psi:X\to T}(-1)^n \Big(\sum_{\substack{A \\ j_A^B\varphi=\psi}}\mu(B,A)\Big) \psi \mpoint$$
Let $g\in [n]$ be such that $g\notin\varphi(X)$.
Then we can modify each $n$-tuple~$A$ by changing freely the image $j_A^B(g)=a_g\in [b_{g-1},b_g]_T$ without changing the equality $j_A^B\varphi=\psi$.
This is because $\varphi(x)\neq g$ for all $x\in X$.
We set $A'=(a_1, \ldots, \widehat a_g,\ldots,a_n)$ and $B'=(b_0, \ldots, \widehat b_{g-1},\ldots,b_{n-1})$ (where $\;\widehat{}\;$ denotes omission) and we let
$$j_{A'}^{B'}:\sou n-\{g\}  \longrightarrow T\,, \qquad h\mapsto a_h \mpoint$$
We obtain
$$\sum_{\substack{A \\ j_A^B\varphi=\psi}}\mu(B,A)
= \sum_{\substack{A' \\ j_{A'}^{B'}\varphi=\psi}}\mu(B',A') \sum_{a_g\in [b_{g-1},b_g]_T} \mu(b_{g-1},a_g) \mpoint$$
But the inner sum on the right is zero by definition of the M\"obius function (and because $b_{g-1}<b_g$).
Therefore the coefficient of every $\psi$ is zero, hence $j^B\varphi=0$.\mpn

(b) This follows from (a).\endpf

For any subset $Y$ of~$[n]$, define
$$\rho_Y: \sou n \longrightarrow \sou n \,, \qquad \rho_Y(h)=
\begin{cases}{0} & {\text { if }\, h=0 \mvirg} \\
{h} & {\text { if }\, h\in Y \mvirg} \\
{h-1} & {\text { if }\, h\notin Y \mpoint} \end{cases}$$

\result{Proposition} \label{ker-j}
Let $T$ be a finite lattice and let $B$ be a strictly increasing $n$-tuple in $T-\{\widehat 1\}$.
Let $\pi^B:T\to \sou n$ be the corresponding surjective morphism and $j^B:\sou n\to T$, constructed as above.
Let $p:F_{\sou n} \to F_{\sou n}/H_{\sou n} $ be the canonical map.
\begin{enumerate}
\item The composite of $\overline F_{j^B}: F_{\sou n}/H_{\sou n} \to F_T$ and $p\, F_{\pi^B}: F_T\to F_{\sou n}/H_{\sou n}$ is the identity morphism of~$F_{\sou n}/H_{\sou n}$.
\item $\overline F_{j^B}:F_{\sou n}/H_{\sou n} \to F_T$ is injective and embeds $F_{\sou n}/H_{\sou n}$ as a direct summand of~$F_T$.
\item $F_{j^B}F_{\pi^B}$ is an idempotent endomorphism of~$F_T$ whose image is isomorphic to~$F_{\sou n}/H_{\sou n}$.
\item $j^B\pi^B$ is an idempotent endomorphism of~$T$.
\item $\pi^B j^B=(-1)^n \displaystyle \sum_{\emptyset \subseteq Y \subseteq [n]} (-1)^{|Y|} \rho_Y$.
\end{enumerate}
\fresult

\pf
(a) It is clear that
$$\pi^B j_A^B(h)=\begin{cases} {0}& {\text { if }\, h=0 \mvirg} \\
{h} & {\text { if }\, j_A^B(h)\in \, ]b_{h-1},b_h]_T \mvirg} \\
{h-1} & {\text { if }\, j_A^B(h)=b_{h-1} \mpoint} \end{cases}$$
Therefore
$$\pi^B j^B=(-1)^n \sum_{\emptyset \subseteq Y \subseteq [n]}
\Big( \sum_{\substack{A \; \\ \pi^Bj_A^B=\rho_Y}}\mu(B,A) \Big) \rho_Y \mpoint$$
If $h\in Y$, then $j_A^B(h)$ runs freely over $]b_{h-1},b_h]_T$. If $h\notin Y$, then $j_A^B(h)=b_{h-1}$ is fixed.
It follows that
$$\sum_{\substack{A \; \\ \pi^Bj_A^B=\rho_Y}}\mu(B,A) = \prod_{h\in Y} \Big( \sum_{a_h\in \,]b_{h-1},b_h]_T} \mu(b_{h-1},a_h) \Big)
=\prod_{h\in Y} (-1) = (-1)^{|Y|} \mvirg$$
using the fact that
$$0=\sum_{a_h\in \,[b_{h-1},b_h]_T} \mu(b_{h-1},a_h) =1+ \sum_{a_h\in \,]b_{h-1},b_h]_T} \mu(b_{h-1},a_h) \mpoint$$
This shows that
$$\pi^B j^B=(-1)^n \sum_{\emptyset \subseteq Y \subseteq [n]} (-1)^{|Y|} \rho_Y \mpoint$$

Suppose now that $Y$ is a proper subset of~$[n]$ and let $g\in [n]$ be maximal such that $g\notin Y$.
If $h>g$, then $\rho_Y(h)=h>g$, while if $h\leq g$, then $\rho_Y(h)\leq g-1$. Therefore $g\notin \rho_Y(\sou n)$.
This shows that $[n]\not\subseteq\rho_Y(\sou n)$, for any proper subset $Y$ of~$[n]$.\par

Now for any function $\varphi : X\to \sou n$, 
$$F_{\pi^B} F_{j^B} (\varphi) = \pi^B j^B \varphi = (-1)^n \sum_{\emptyset \subseteq Y \subseteq [n]} (-1)^{|Y|} \rho_Y\varphi
=\varphi + \sum_{Y\neq[n]} (-1)^{n+|Y|} \rho_Y\varphi \mpoint$$
But $[n]\not\subseteq\rho_Y(\sou n)$ if $Y\neq[n]$, hence $[n]\not\subseteq\rho_Y\varphi(X)$.
In other words, $\rho_Y\varphi\in H_{\sou n}(X)$, so that $F_{\pi^B} F_{j^B} (\varphi) = \varphi \pmod{H_{\sou n}(X)}$.
Composing with the canonical map ${p:F_{\sou n}(X) \to F_{\sou n}(X)/H_{\sou n}(X)}$ and writing $p(\varphi)=\overline\varphi$, we obtain
$$p F_{\pi^B}\overline F_{j^B} (\overline\varphi) = \overline\varphi \mvirg$$
as was to be shown.\mpn

(b) This follows immediately from~(a).\mpn

(c) This follows immediately from~(a) and the obvious equality $\overline F_{j^B}p\,F_{\pi^B} = F_{j^B}F_{\pi^B}$.\mpn

(d) Since the functor $F_?:k\CL \to \CF_k$ is fully faithful by Theorem~\ref{fully faithful}, $j^B\pi^B$ must be an idempotent
because its image $F_{j^B}F_{\pi^B}$ is an idempotent. Alternatively, it is not difficult to compute directly
$$j^B\pi^Bj^B=j^B(-1)^n \sum_{\emptyset \subseteq Y \subseteq [n]} (-1)^{|Y|} \rho_Y
=j^B \Id_{\sou n} +(-1)^n \sum_{ Y \neq [n]} (-1)^{|Y|} j^B\rho_Y=j^B \mvirg$$
because $[n]\not\subseteq\rho_Y(\sou n)$ if $Y\neq[n]$, hence $j^B\rho_Y=0$ by Proposition~\ref{vanish-on-H}.\mpn

(e) The required equality has already been proved above.
\endpf

Our aim is to show that the idempotents $j^B\pi^B$ are orthogonal.
In order to understand the product of two idempotents $j^C\pi^C$ and $j^B\pi^B$ we need to have more information about $\pi^C j^B$.
This is the purpose of the next two propositions.

\result{Proposition} \label{compose-nonzero}
Let $T$ be a finite lattice and let
$$\begin{array}{lcl}
B=(b_0,b_1,\ldots,b_{n-1}) &\quad\text{ with } &\widehat 0 \leq b_0 < b_1 < \ldots < b_{n-1} < b_n=\widehat 1 \mvirg \\
C=(c_0,c_1,\ldots,c_{m-1}) &\quad\text{ with } &\widehat 0 \leq c_0 < c_1 < \ldots < c_{m-1} < c_m=\widehat 1 \mpoint
\end{array}$$
If $\pi^C j^B\neq 0$, then the restriction of $\pi^C$ to the subset $\{b_0,b_1,\ldots,b_{n-1}, b_n\}$ is injective.
In particular, $n\leq m$.
\fresult

\pf We have
$$\pi^Cj^B=(-1)^n\sum_A\mu(B,A) \pi^Cj_A^B
=(-1)^n \sum_{\psi:\sou n\to\sou m} \Big(\sum_{\substack{A\; \\ \pi^Cj_A^B=\psi}}\mu(B,A)\Big) \psi \mpoint$$
Now fix some morphism $\psi:\sou n\to\sou m$ and, for every $h\in [n]$, define
$$U_h=\{ a\in [b_{h-1},b_h]_T \mid \pi^C(a)=\psi(h)\} \subseteq [b_{h-1},b_h]_T \mpoint$$
Then the condition $\pi^Cj_A^B=\psi$ is equivalent to $j_A^B(h)\in U_h$ for every $h\in [n]$, that is, $a_h\in U_h$ for every $h\in [n]$.
In particular $U_h\neq\emptyset$ in that case.
Since all elements of~$U_h$ have the same image under~$\pi^C$, so has their join and therefore $U_h$ has a supremum
$$u_h=\sup (U_h) \in U_h \mpoint$$
Define
$$V_h=\{a\in [b_{h-1},b_h]_T \mid \pi^C(a)\leq \psi(h)-1\} \mpoint$$
Then we clearly have $[b_{h-1},u_h]_T = V_h\sqcup U_h$ because any $a\leq u_h$ satisfies either $\pi^C(a)\leq \psi(h)-1$ or $\pi^C(a)= \psi(h)$.
If $V_h$ is nonempty (i.e. if $b_{h-1}\notin U_h$), then $V_h$ has a supremum $v_h=\bigvee_{t\in V_h}t$ and $V_h=[b_{h-1},v_h]_T$.\par

Now, in the expression above for $\pi^Cj^B$, the coefficient of~$\psi$ is $(-1)^n z_\psi$, where
$$z_\psi=\sum_{\pi^Cj_A^B=\psi} \mu(B,A)=\Big(\sum_{a_1\in U_1} \mu(b_0,a_1)\Big) \Big(\sum_{a_2\in U_2} \mu(b_1,a_2)\Big) \ldots \Big(\sum_{a_n\in U_n} \mu(b_{n-1},a_n)\Big)$$
and we examine one inner sum $\sum_{a\in U_h} \mu(b_{h-1},a)$.
We have already observed that $U_h\neq\emptyset$.
Now there are two cases~:

\begin{itemize}
\item[(a)] If $V_h=\emptyset$, then $U_h=[b_{h-1},u_h]_T$ and $\sum_{a\in U_h} \mu(b_{h-1},a)=0$, unless $b_{h-1}=u_h$.
\item[(b)] If $V_h\neq \emptyset$, then we claim that $\sum_{a\in U_h} \mu(b_{h-1},a)=0$, unless $b_{h-1}=v_h$.
\end{itemize}
To prove the claim, suppose that $b_{h-1}<v_h$ (hence $b_{h-1}<u_h$ because $v_h<u_h$). Then we have $[b_{h-1},u_h]_T=[b_{h-1},v_h]_T \sqcup U_h$ and therefore
$$\sum_{a\in U_h} \mu(b_{h-1},a)=\sum_{a\in [b_{h-1},u_h]_T} \mu(b_{h-1},a) - \sum_{a\in [b_{h-1},v_h]_T} \mu(b_{h-1},a) =0 \mvirg$$
because $b_{h-1}<u_h$ and $b_{h-1}<v_h$.\par

It follows that the sum $\sum_{a\in U_h} \mu(b_{h-1},a)$ is nonzero in only two cases~:
\begin{itemize}
\item[(a)] $U_h=\{b_{h-1}\}$. Then $\sum_{a\in U_h} \mu(b_{h-1},a)=\mu(b_{h-1},b_{h-1})=1$.
\item[(b)] $V_h=\{b_{h-1}\}$. Then $U_h=]b_{h-1},u_h]_T$ and $\sum_{a\in U_h} \mu(b_{h-1},a)=-1$ because
$$\begin{array}{rcl}
0=\sum_{a\in [b_{h-1},u_h]_T} \mu(b_{h-1},a)&=&\mu(b_{h-1},b_{h-1})+\sum_{a\in ]b_{h-1},u_h]_T} \mu(b_{h-1},a) \\
&=&1+\sum_{a\in U_h} \mu(b_{h-1},a) \mpoint
\end{array}$$
\end{itemize}

If the coefficient $z_\psi$ is nonzero, then this sum must be nonzero for every $h\in[n]$ and we obtain $z_\psi=(-1)^p$ where $p$ is the number of times the integer $h$ satisfies case~(b) (and then the coefficient of~$\psi$ in the expression for $\pi^Cj^B$ is $(-1)^{n+p}$).\par

Now suppose that the restriction of $\pi^{C}$ to $\{b_0,b_1,\ldots,b_{n-1}, b_n\}$ is not injective.
We want to prove that the coefficient $z_\psi$ is zero (for our fixed morphism $\psi:\sou n\to\sou m$).
Since $\pi^C$ is order-preserving and not injective on the subset above, there exists $h\in[n]$ such that $\pi^C(b_{h-1})=\pi^C(b_h)$.
We know that $z_\psi=0$ if $U_h$ is empty. Thus we can assume that $U_h$ contains an element~$a$.
Since $b_{h-1}\leq a\leq b_h$, we have $\pi^C(b_{h-1})\leq \pi^C(a)\leq \pi^C(b_h)$, hence equality $\pi^C(b_{h-1})= \pi^C(a)= \pi^C(b_h)$.
This means that both $b_{h-1}$ and $b_h$ belong to~$U_h$ (because $\pi^C(a)=\psi(h)$ by definition of~$U_h$).
Now in both cases (a) and~(b) above, the set $U_h$ does not contain simultaneously $b_{h-1}$ and $b_h$ (it is either $\{b_{h-1}\}$ or~$]b_{h-1},u_h]_T$).
This shows that the sum $\sum_{a\in U_h} \mu(b_{h-1},a)$ cannot be nonzero. Therefore $z_\psi=0$.\par

This argument holds for every morphism $\psi:\sou n\to\sou m$ and so every coefficient $z_\psi$ is zero. It follows that $\pi^Cj^B=0$, as was to be proved.\endpf

\result{Proposition} \label{compute-composite}
Let $T$ be a finite lattice and let
$$\begin{array}{lcl}
B=(b_0,b_1,\ldots,b_{n-1}) &\quad\text{ with } &\widehat 0 \leq b_0 < b_1 < \ldots < b_{n-1} < b_n=\widehat 1 \mvirg \\
C=(c_0,c_1,\ldots,c_{m-1}) &\quad\text{ with } &\widehat 0 \leq c_0 < c_1 < \ldots < c_{m-1} < c_m=\widehat 1 \mvirg \\
D=(d_0,d_1,\ldots,d_{m-1}) &\quad\text{ with } &\widehat 0 \leq d_0 < d_1 < \ldots < d_{m-1} < d_m=\widehat 1 \mpoint
\end{array}$$
Then
$$j^D\pi^Cj^B=\begin{cases}{j^D} & { \text{ if } \, C=B \,,} \\
{0} & { \text{ if } \, C\neq B \,.}\end{cases}$$
\fresult

\pf
We assume that $j^D\pi^Cj^B\neq 0$, and in particular $\pi^Cj^B\neq 0$.
Write first
$$\pi^Cj^B=\sum_A \mu(B,A) \pi^Cj_A^B$$
where $A=(a_1,a_2,\ldots,a_n)$ and $a_h\in [b_{h-1},b_h]$ for every $h\in[n]$.
Let $\psi:\sou n\to \sou m$ be a map appearing with a nonzero coefficient in the expression of~$\pi^Cj^B\neq 0$ as a linear combination of morphisms in the category~$\CL$.
Let $A$ be such that $\pi^Cj_A^B=\psi$. 
Since $j^D\pi^Cj^B\neq 0$, we can also assume that $\psi$ is such that $j^D\psi\neq0$.
Proposition~\ref{vanish-on-H} implies that the function $\psi:\sou n\to \sou m$ must satisfy $[m]\subseteq \psi(\sou n)$.
Since $\psi$ is a morphism of lattices, we also have $\psi(0)=0$.
Therefore $\psi:\sou n\to \sou m$ must be surjective. In particular $n\geq m$.\par

By Proposition~\ref{compose-nonzero}, $\pi^Cj^B\neq 0$ implies that $n\leq m$. Therefore $n=m$.
Since the map $\pi^Cj_A^B:\sou n\to\sou n$ is order-preserving and surjective, it must be the identity map.
This shows that whenever $A$ is such that $j^D\pi^Cj_A^B\neq 0$, then $\pi^Cj_A^B=\Id_{\sou n}$.
Therefore, the functions $\psi$ which appear with a nonzero coefficient in the expression of~$\pi^Cj^B$
are $\Id_{\sou n}$ and maps in the kernel of~$j^D$.
\par

Now Proposition~\ref{compose-nonzero} also asserts that the restriction of $\pi^C$
$$\pi^C:\{b_0,b_1,\ldots,b_{n-1}, b_n\}\longrightarrow \sou m$$
is injective.
Since $n=m$, this must be a bijection, and since $\pi^C$ is order-preserving, we must have $\pi^C(b_h)=h$ for every $h\in \sou n$.
As in the proof of Proposition~\ref{compose-nonzero}, associated with the map $\Id_{\sou n}:{\sou n} \to {\sou n}$, there is a corresponding subset
$$U_h=\{a\in [b_{h-1},b_h]_T \mid \pi^C(a)=\Id_{\sou n}(h)\}=\{a\in [b_{h-1},b_h]_T \mid \pi^C(a)=h\} \mpoint$$
Thus we have $b_h\in U_h$ for every $h\in \sou n$.\par

Since $\Id_{\sou n}$ appears with a nonzero coefficient $(-1)^n z_{\Id}$ in the expression of~$\pi^Cj^B$, the proof of Proposition~\ref{compose-nonzero} shows that $\sum_{a\in U_h} \mu(b_{h-1},a)$ is nonzero. Therefore we must be in one of the two cases (a) or (b) of that proof.
But the first case cannot hold because $b_h\notin U_h=\{b_{h-1}\}$.
So we are in case~(b) and we have $b_h\in U_h=]b_{h-1},u_h]_T\subseteq ]b_{h-1},b_h]_T$. This forces $u_h=b_h$ and $U_h=]b_{h-1},b_h]_T$.\par

The equality $U_h=]b_{h-1},b_h]_T$ means that $\pi^C$ coincides with $\pi^B$ on this interval. This holds for every $h\in[n]$.
We now show, by descending induction, that $c_h=b_h$ for every $h\in[n]$. First $c_n=\widehat 1=b_n$.
Now since $\pi^C(b_h)=h$, we have $b_h\leq c_h$, but $b_h\not\leq c_{h-1}$, for every $h\in\sou n$.
Assume that $c_h=b_h$ for every $h\geq i$ and suppose for contradiction that $b_{i-1}<c_{i-1}$.
Then $b_{i-1}<c_{i-1}<c_i=b_i$, hence $c_{i-1}\in \, ]b_{i-1},b_i]$, which implies that $\pi^B(c_{i-1})=i$.
But $\pi^B$ coincides with $\pi^C$ on the interval, so $\pi^C(c_{i-1})=i$, which contradicts the definition of~$\pi^C$.\par

We have now proved that $n=m$ and $C=B$ whenever $j^D\pi^Cj^B\neq 0$.
Moreover, in that case, we have seen that $\pi^Bj^B$ is the sum of a multiple of~$\Id_{\sou n}$ and morphisms in the kernel of~$j^D$.
The proof of Proposition~\ref{compose-nonzero} shows that the coefficient of~$\Id_{\sou n}$ is equal to $(-1)^{n+p}$, where $p$ is the number of times the integer $h$ satisfies case~(b). But we have noticed that case~(a) cannot occur, so $p=n$ and the coefficient is~1.
Thus $\pi^Bj^B=\Id_{\sou n} \pmod{\Ker(j^D)}$ and it follows that $j^D\pi^Bj^B=j^D$.\par

Alternatively, the equality $\pi^Bj^B=\Id_{\sou n} \pmod{\Ker(j^D)}$ also follows from Proposition~\ref{vanish-on-H} and Proposition~\ref{ker-j},
because $\pi^B j^B=(-1)^n \sum_{\emptyset \subseteq Y \subseteq [n]} (-1)^{|Y|} \rho_Y$ and $j^D\rho_Y=0$ whenever $Y\neq[n]$
since $[n]\not\subseteq \rho_Y(\sou n)$.
\endpf

Given a finite lattice $T$, let $\CP_{T,n}$ be the set of all strictly increasing $n$-tuples $(d_0,d_1,\ldots,d_{n-1})$ in~$T-\{\widehat1\}$, that is, such that
$\widehat 0 \leq d_0 < d_1 < \ldots < d_{n-1} < d_n=\widehat 1$.
The set $\CP_{T,n}$ corresponds bijectively to the set of all surjective morphisms $T\to \sou n$.
If $D,C\in \CP_{T,n}$, we define
$$f_{D,C}=j^D\pi^C:T\longrightarrow T \mpoint$$
In particular, $f_{B,B}=j^B\pi^B$ is the idempotent of Proposition~\ref{ker-j}.

\result{Corollary} \label{matrix-coeff}
Let $T$ be a finite lattice.
\begin{enumerate}
\item Let $A,B \in \CP_{T,n}$ and $D,C\in \CP_{T,m}$, where $n,m\geq0$ are two integers. Then
$$f_{D,C}f_{B,A}=\begin{cases}{f_{D,A}} & { \text{ if } \, C=B \; (\text{hence } \, n=m) \,,} \\
{0} & { \text{ if } \, C\neq B \,.}\end{cases}$$
\item When $n\geq0$ varies and $B\in \CP_{T,n}$ varies, the idempotents $f_{ B,B}$ are pairwise orthogonal.
\end{enumerate}
\fresult

\pf 
(a) By Proposition~\ref{compute-composite}, we have
$$f_{D,C}f_{B,A}=j^D\pi^Cj^B\pi^A= j^D\pi^A$$
if $C=B$, and zero otherwise.\mpn

(b) This follows from (a).
\endpf

Let $M_{|\CP_{T,n}|}(k)$ denote the matrix algebra of size~$|\CP_{T,n}|$, with rows and columns indexed by the set~$\CP_{T,n}$.
If $D,C\in \CP_{T,n}$, we let $m_{D,C}$ denote the matrix with coefficient~1 in position~$(D,C)$ and 0 elsewhere.

\result{Theorem} \label{matrix-algebra}
Let $T$ be a finite lattice, let $\CP_{T,n}$ denote the set of all strictly increasing $n$-tuples in~$T-\{\widehat1\}$,
and let $N$ be the maximal length of a strictly increasing sequence in~$T-\{\widehat1\}$.
\begin{enumerate}
\item The map
$$\CI_T: \; \bigoplus_{n=0}^N M_{|\CP_{T,n}|}(k) \longrightarrow \End_{k\CL}(T) \,, \qquad m_{D,C} \mapsto f_{D,C}$$
is an algebra homomorphism  (without unit elements).
\item $\CI_T$ is injective.
\item The image of $\CI_T$ is equal to the subalgebra~$\CE_T$ (without unit element) of $\End_{k\CL}(T)$ having a $k$-basis consisting of all maps whose image is totally ordered (and which are morphisms in the category~$\CL$).
\end{enumerate}
\fresult

\pf
(a) If $|C|\neq|B|$, then $m_{D,C}$ and $m_{B,A}$ are not in the same block, so their product is~0, while the product $f_{D,C}f_{B,A}$ is also zero.
If $|C|=|B|$, then the relations of Corollary~\ref{matrix-coeff} are the standard relations within a matrix algebra of size~$|\CP_{T,n}|$.\mpn

(b) Since the elements $m_{D,C}$ form a $k$-basis of $\bigoplus_{n=0}^N M_{|\CP_{T,n}|}(k)$, it suffices to prove that their images $f_{D,C}$ are $k$-linearly independent. Suppose that
$$\sum_{n=0}^N \;\sum_{D,C\in\CP_{T,n}} \lambda_{D,C} f_{D,C} =0 \mvirg$$
where $\lambda_{D,C}\in k$.
Multiply on the left by the idempotent $f_{A,A}$ and on the right by the idempotent $f_{B,B}$.
By Corollary~\ref{matrix-coeff}, we obtain $\lambda_{A,B} f_{A,B}=0$.
Since $f_{A,B}$ is a linear combination of distinct maps $T\to T$, one of them appearing with coefficient $\pm1$, we must have $\lambda_{A,B}=0$.\mpn

(c) It is clear that $\CE_T$ is a subalgebra. Moreover, every map $j_A^D\pi^C$ has an image which is totally ordered, by construction.
Therefore $f_{D,C}=j^D\pi^C\in\CE_T$ and hence $\Im(\CI_T)\subseteq\CE_T$.\par

Now $\CE_T$ has a $k$-basis consisting of all maps $\varphi_{U,V}:T\to T$ described as follows.
First $U=(u_0,u_1,\ldots,u_{n-1})\in\CP_{T,n}$, while $V\in \CY_{T,n}$, where $\CY_{T,n}$ denotes the set of all strictly increasing $n$-tuples
$V=(v_1,v_2,\ldots, v_n)$ in~$T-\{\widehat0\}$ (that is $\widehat0=v_0 < v_1<v_2<\ldots<v_n\leq \widehat1$).
Define $\lambda^V:\sou n\to T$ by $\lambda^V(i)=v_i$ for every~$i\in\sou n$ and then set
$\varphi_{U,V}=\lambda^V \pi^U$. Then
$$\CB_n=\{\varphi_{U,V} \mid U\in \CP_{T,n}\,, \, V\in \CY_{T,n} \}$$
is a $k$-basis of the submodule $\CE_{T,n}$ generated by all endomorphisms whose image is isomorphic to~$\sou n$,
while $\CB=\bigcup_{n=0}^N \CB_n$ is a $k$-basis of~$\CE_T=\bigoplus_{n=0}^N \CE_{T,n}$.
We have seen that $\CB'=\bigcup_{n=0}^N \{f_{D,C} \mid D,C\in \CP_{T,n} \}$ is a $k$-basis of~$\Im(\CI_T)$.
It is an easy exercise to show that there is a bijection between $\CY_{T,n}$ and $\CP_{T,n}$.
Therefore $\CB$ and $\CB'$ have the same cardinality.
In other words $\Im(\CI_T)$ and $\CE_T$ are free $k$-modules of the same rank.\par

Now we allow the base ring $k$ to vary and we write an index~$k$ to emphasize the dependence on~$k$.
Thus we have an injective algebra homomorphism
$$\CI_{T,k}: \; \bigoplus_{n=0}^N M_{|\CP_{T,n}|}(k) \longrightarrow \CE_{T,k} \subseteq \End_{k\CL}(T)$$
and we let $\CQ_k:=\CE_{T,k}/\Im(\CI_{T,k})$, so that we have a short exact sequence
$$\xymatrix{
0 \ar[r] & \Im(\CI_{T,k}) \ar[r]^-{j_k} & \CE_{T,k} \ar[r]^-{p_k} & \CQ_k \ar[r] & 0 \mvirg
}$$
where $j_k$ is the inclusion map and $p_k$ the canonical surjection.
In the case of the ring of integers $\Z$, we see that $\CQ_\Z$ is a finite abelian group, because $\Im(\CI_{T,\Z})$ and $\CE_{T,\Z}$ are free $\Z$-modules of the same rank.
Tensoring with the prime field~$\F_p$, where $p$ is a prime, we obtain
$$\xymatrix{
0 \ar[r] & \Tor(\F_p,\CQ_\Z) \ar[r] & \F_p\otimes \Im(\CI_{T,\Z}) \ar[r]^-{1\otimes j_\Z} & \F_p\otimes \CE_{T,\Z} \ar[r]^-{1\otimes p_\Z} & \F_p\otimes \CQ_\Z \ar[r] & 0 \mpoint
}$$
Using the canonical bases of $\Im(\CI_{T,k})$ and $\CE_{T,k}$ respectively, we see that $k\otimes \Im(\CI_{T,\Z}) \cong \Im(\CI_{T,k})$,
and similarly $k\otimes \CE_{T,\Z} \cong \CE_{T,k}$. Moreover the inclusion map $1\otimes j_k$ corresponds, under these isomorphisms,
to the inclusion map~$j_k$.
In particular, we obtain
$$\xymatrix{
0 \ar[r] & \Tor(\F_p,\CQ_\Z) \ar[r] & \Im(\CI_{\F_p}) \ar[r]^-{j_{\F_p}} & \CE_{\F_p} \ar[r]^-{1\otimes p_\Z} & \F_p\otimes \CQ_\Z \ar[r] & 0 \mpoint
}$$
By (b), we know that $j_{\F_p}$ is injective, hence $\Tor(\F_p,\CQ_\Z)=\{0\}$.
But since this holds for every prime number~$p$ and $\CQ_\Z$ is finite, we must have $\CQ_\Z=\{0\}$.
It follows that the inclusion map $j_\Z:\Im(\CI_{T,\Z}) \to \CE_{T,\Z}$ is an isomorphism.
Tensoring with~$k$, it follows that the inclusion map $j_k:\Im(\CI_{T,k}) \to \CE_{T,k}$ is an isomorphism as well.
In other words, $\Im(\CI_{T,k})=\CE_{T,k}$, as was to be shown.
\endpf

\begin{rem}{Remark} \label{change-basis} Let $\CB$ be the canonical basis of $\CE_T$ described in the proof.
The change of basis from $\CB$ to the basis $\{f_{D,C}\}$ of~$\Im(\CI_T)$ is not obvious.
By construction, every map $j_A^D\pi^C$ belongs to~$\CB$, but beware of the fact that if $C$ and $D$ are $n$-tuples,
then $j_A^D\pi^C$ may be a composite $T\to \sou m \to T$ for some $m<n$, because the $n$-tuple $A= (a_1,a_2,\ldots,a_n)$ is increasing but not necessarily strictly increasing, hence may consist of $m$ distinct elements for some $m<n$.
\end{rem}

The image under~$\CI_T$ of the identity element of $\bigoplus_{n=0}^N M_{|\CP_{T,n}|}(k)$ is an idempotent~$e_T$ of~$\End_{k\CL}(T)$ and $e_T$ is an identity element of~$\CE_T$. We now prove that the we actually get central idempotents.

\result{Theorem} \label{central-e}
For every finite lattice~$T$, let $\CE_T=\Im(\CI_T)$ be the subalgebra of~$\End_{k\CL}(T)$ appearing in Theorem~\ref{matrix-algebra},
and let $e_T$ be the identity element of~$\CE_T$.
\begin{enumerate}
\item $e_T=\displaystyle\sum_{n=0}^N \sum_{B\in\CP_{T,n}} f_{B,B}$.
\item For any finite lattice~$T'$ and any morphism $\theta\in\Hom_{k\CL}(T,T')$, we have $\theta e_T = e_{T'} \theta$.
In other words, the family of idempotents~$e_T$, for $T\in\CL$, is a natural transformation of the identity functor~$\Id_{k\CL}$.
\item $e_T$ is a central idempotent of~$\End_{k\CL}(T)$.
\item The subalgebra $\CE_T$ is a direct product factor of~$\End_{k\CL}(T)$, that is, there exists a subalgebra $\CD$ such that $\End_{k\CL}(T)=\CE_T\times\CD$ (where $\CE_T$ is identified with $\CE_T\times\{0\}$ and $\CD$ with $\{0\}\times\CD$, as usual).
\end{enumerate}
\fresult

\pf (a) The identity element of $\displaystyle\bigoplus_{n=0}^N M_{|\CP_{T,n}|}(k)$ is equal to
$\displaystyle\sum_{n=0}^N \sum_{B\in\CP_{T,n}} m_{B,B}$.
Taking its image under~$\CI_T$ yields the required formula.\mpn

(b) We have seen in the proof of Theorem~\ref{matrix-algebra} that every element of the canonical basis of~$\CE_T$ has the form $\lambda^V \pi^U$,
where $U=(u_0,u_1,\ldots,u_{n-1})\in\CP_{T,n}$ is a strictly increasing $n$-tuple in~$T-\{\widehat1\}$, while $V=(v_1,v_2,\ldots, v_n)$ belongs to the set $\CY_{T,n}$ of all strictly increasing $n$-tuples in~$T-\{\widehat0\}$.
We now compute their opposite, as in Lemma~\ref{Lop}.\par

From the surjective morphism $\pi^U: T \to \sou n$, we obtain $(\pi^U)\op: \sou n\op \to T\op$ defined as follows~: for every $h\in \sou n$,
$$(\pi^U)\op(h)=\bigvee_{\pi^U(t)\leq h}t=\sup\big((\pi^U)^{-1}(h)\big)=u_h \mvirg$$
in view of the way $U$ is associated to~$\pi^U$ (see the beginning of Section~\ref{Section-embed-total}).
But $u_h=\lambda^{U\op}(h)$, hence $(\pi^U)\op=\lambda^{U\op}$.\par

Applying this to the case $U=V\op$, where $V\in \CY_{T,n}$, we obtain $(\pi^{V\op})\op=\lambda^{{(V\op})\op}=\lambda^V$.
Taking opposites yields $({\lambda^V})\op=(\pi^{V\op})$.\par

It follows that the opposite of the canonical basis element $\lambda^V \pi^U$ of~$\CE_T$ is the canonical basis element $\lambda^{U\op} \pi^{V\op}$ of~$\CE_{T\op}$. Therefore, the opposite of the identity element~$e_T$ of~$\CE_T$ must belong to~$\CE_{T\op}$.
Moreover, it must be the identity element of~$\CE_{T\op}$, because taking opposites behaves well with respect to composition, by Lemma~\ref{Lop}.
Therefore $(e_T)\op=e_{T\op}$.\par

Now if $\theta:T\to T'$ is a morphism in~$\CL$, then the image of a totally ordered subset of~$T$ is a totally ordered subset of~$T'$.
It follows that composition with~$\theta$ maps $\CE_T$ to a linear combination of maps with a totally ordered image, hence invariant under the idempotent element~$e_{T'}$. In other words, we have
$$\theta \,e_T=e_{T'} \,\theta \,e_T \mpoint$$
Applying this equation to ${T'}\op$, $T\op$, and the morphism $\theta\op:{T'}\op \to T\op$, we obtain $\theta\op \,e_{{T'}\op}=e_{T\op} \,\theta\op \,e_{{T'}\op}$.
Passing to opposites and using the above equality $(e_T)\op=e_{T\op}$, we get
$$e_{T'} \, \theta \, e_T = e_{T'} \, \theta \mpoint$$
The two displayed equations then yield $\theta e_T = e_{T'} \theta$.
This holds as well if $\theta$ is replaced by a linear combination of morphisms, as was to be shown.\mpn

(c) This is a special case of~(b).\mpn

(d) This follows immediately from~(c).
\endpf


\section{The case of a total order} \label{Section-total}

\bigskip
\noindent
In this section, we consider the case of a totally ordered lattice~$\sou n$, where $n\in\N$.
We determine completely the ring of endomorphisms of $\sou n$ in the category~$k\CL$ and we deduce a direct sum decomposition of~$F_{\sou n}$.
For the lattice~$\sou n$, the set $\CP_{\sou n,r}$ of all strictly increasing $r$-tuples $(d_0,d_1,\ldots,d_{r-1})$ in~$\sou n-\{n\}$
is just the set of all subsets of~$\{0,1,\ldots,n-1\}$ of size~$r$, because any such subset is totally ordered.
Throughout this section, we use subsets of size~$r$ instead of strictly increasing $r$-tuples.
In particular $|\CP_{\sou n,r}|={n\choose r}$.

\result{Theorem} \label{end total} Let $n\in \N$.
\begin{enumerate}
\item The homomorphism of $k$-algebras of Theorem~\ref{matrix-algebra}
$$\CI_{\sou n}:\bigoplus_{m=0}^n M_{|\CP_{\sou n,r}|}(k)\longrightarrow \End_{k\CL}(\sou{n}) \,, \qquad
m_{D,C}\mapsto f_{D,C} \mvirg$$
is an isomorphism.
\item In particular, if $k$ is a field, then $\End_{k\CL}(\sou{n})$ is semi-simple.
\end{enumerate}
\fresult

\pf
First note that (b) follows from~(a), because any matrix algebra $M_q(k)$ is simple, so that the direct sum is semi-simple.
Since any map $\varphi:\sou n\to \sou n$ has an image which is totally ordered,
the subalgebra $\CE_T$ of~$\End_{k\CL}(\sou{n})$ appearing in Theorem~\ref{matrix-algebra}
is the whole of~$\End_{k\CL}(\sou{n})$.
By Theorem~\ref{matrix-algebra}, the morphism $\CI_{\sou n}$ is surjective and injective, hence an isomorphism.
\endpf

\begin{rem}{Remark} \label{change-basis-2} Since $\CE_T$ is the whole of~$\End_{k\CL}(\sou{n})$, the canonical basis~$\CB$ of~$\CE_T$ is a basis of~$\End_{k\CL}(\sou{n})$, which also has another basis consisting of the morphisms~$f_{D,C}$.
As mentioned in Remark~\ref{change-basis}, the change of basis is not straightforward, but it can be made more explicit in the case of a totally ordered lattice considered here.\par

Every element of the canonical basis of~$\CE_T$ has the form $\lambda^V \pi^U$,
where $U=(u_0,u_1,\ldots,u_{n-1})\in\CP_{T,n}$ is a strictly increasing $n$-tuple in~$T-\{\widehat1\}$, while $V=(v_1,v_2,\ldots, v_n)\in \CY_{T,n}$ is a strictly increasing $n$-tuple in~$T-\{\widehat0\}$.
On the other hand, the morphisms $f_{D,C}$ are parametrized by pairs $D,C\in \CP_{T,n}$.
But we have an obvious bijection from $\CP_{T,n}$ to $\CY_{T,n}$, mapping an $n$-tuple $(u_0,u_1,\ldots,u_{n-1})\in\CP_{T,n}$ to $(u_0+1,u_1+1,\ldots,u_{n-1}+1)\in\CY_{T,n}$. Thus we can parametrize both bases by the same set $\CP_{T,n}\times \CP_{T,n}$. Then it is not hard to see that the matrix of the change of basis is unitriangular.
Actually, this provides another proof of the fact that the map $\CI_{\sou n}$ is an isomorphism.
\end{rem}

\begin{rem}{Remark} Theorem~\ref{end total} is similar to the result proved in~\cite{FHH} about the planar rook algebra. Over the field~$\C$ of complex numbers, this algebra is actually isomorphic to $\End_{\C\CL}(\sou{n})$. However, the planar rook monoid is not isomorphic to the monoid of endomorphisms of $\sou{n}$ in~$\CL$, because it turns out that they do not have the same number of idempotents, even when $n=2$. Only the corresponding monoid algebras become isomorphic (over~$\C$).
\end{rem}

Now we consider the central idempotents of $\End_{k\CL}(\sou{n})$ corresponding to the above decomposition into matrix algebras.

\result{Notation} \label{beta}
For an integer $m$ with $0\leq m\leq n$, set 
$$\beta_{n,m}:=\sum_{B\in\CP_{\sou n,m}} f_{B,B}\mpoint$$
In particular, for $m=n$ and $B=\{0,1,\ldots,n-1\}$, we have $\pi^B=\Id_{\sou n}$ and we define
$$\varepsilon_n:=\beta_{n,n}=f_{B,B}=j^B\pi^B=j^B=(-1)^n \sum_{\emptyset \subseteq Y \subseteq [n]} (-1)^{|Y|} \rho_Y\mvirg$$
where, as before, $\rho_Y \in\End_\CL(\sou{n})$ is defined by $\rho_Y(h)=h$ if $h\in Y$ and $\rho_Y(h)=h-1$ otherwise.
\fresult

\result{Proposition} \label{central-idempotents} 
The elements~$\beta_{n,m}$, for $0\leq m\leq n$, are orthogonal central idempotents of $\End_{k\CL}(\sou{n})$, and their sum is equal to the identity.
In particular, the central idempotent $\varepsilon_n$ satisfies
$$\varepsilon_n\End_{k\CL}(\sou{n})=k\varepsilon_n\mpoint$$
\fresult

\pf
For $B\in\CP_{\sou n,m}$, the inverse image of $f_{B,B}$ under the algebra isomorphism~$\CI$
of Theorem~\ref{end total} is the matrix $m_{B,B}$ of the component $M_{|\CP_{\sou n,m}|}(k)$ indexed by~$m$. 
Summing over all $B\in\CP_{\sou n,m}$, it follows that the inverse image of~$\beta_{n,m}$ under~$\CI$ is the identity element of $M_{|\CP_{\sou n,m}|}(k)$.
The first statement follows.\par

In the case $m=n$, the set $\CP_{\sou n,m}$ consists of the singleton~$B=\{0,1,\ldots,n-1\}$ and the corresponding matrix algebra has size~1.
We see that the inverse image of $\varepsilon_n$ under~$\CI$
is the identity element $m_{B,B}$ of the component $M_1(k)=k$.
Clearly $\varepsilon_n\End_{k\CL}(\sou{n})\cong M_1(k)$ is a free $k$-module of rank~1,
hence $\varepsilon_n\End_{k\CL}(\sou{n})=k\varepsilon_n$.
\endpf

We want to use the functor $F_?: k\CL \to \CF_k$ of Section~\ref{Section-FT} to deduce information on the correspondence functor~$F_{\sou{n}}$.
By Theorem~\ref{projective-distributive}, we already know that $F_{\sou{n}}$ is projective,
because the total order $\sou{n}$ is a distributive lattice.
If $n=0$ (in which case $[n]=\emptyset$) and if $n=1$ (in which case $[n]=\{1\}$), we recover the cases already considered in Section~5 of~\cite{BT2}.
Our purpose is to treat now the general case.\par

We apply the functor $F_?: k\CL \to \CF_k$ to the map $j^B\in \Hom_{k\CL}(\sou{m},\sou{n})$ defined in Section~\ref{Section-embed-total}, where $B\in\CP_{\sou n,m}$.
By Proposition~\ref{vanish-on-H} we obtain a morphism  
$$F_{j^B}:F_{\sou{m}} \longrightarrow F_{\sou{n}}$$
which vanishes on~$H_{\sou m}$. By Proposition~\ref{ker-j}, this induces an injective morphism
$$\overline F_{j^B}:F_{\sou{m}}/H_{\sou{m}} \longrightarrow F_{\sou{n}}$$
which embeds $F_{\sou{m}}/H_{\sou{m}}$ as a direct summand of~$F_{\sou{n}}$, corresponding to the idempotent $f_{B,B}=j^B\pi^B$.
In particular, for $m=n$, we have $f_{B,B}=j^B=\varepsilon_n$ and we obtain an idempotent endomorphism
$F_{\varepsilon_n}$ of~$F_{\sou{n}}$ with kernel~$H_{\sou{n}}$.

\result{Theorem} \label{functor-Sn}
Let $n\in\N$ and let $\S_n=F_{\sou{n}}/H_{\sou{n}}$.
There are isomorphisms of correspondence functors
\begin{eqnarray*}
F_{\varepsilon_n}F_{\sou{n}}&\cong& \S_n \mvirg \\
F_{\beta_{n,m}}F_{\sou{n}}&\cong& \S_m^{\oplus {n\choose m}}\mvirg\; \text{ for each }\; 0\leq m\leq n\mvirg\\
F_{\sou{n}}&\cong&\bigoplus_{\substack{ 0\leq m\leq n \\ B\in\CP_{\sou n,m} }} \S_{|B|}\mpoint
\end{eqnarray*}
\fresult

\pf  By Theorem~\ref{fully faithful}, the functor $F_?$ induces an isomorphism of $k$-algebras 
$$\End_{k\CL}(\sou{n})\cong\End_{\CF_k}(F_{\sou{n}})\mpoint$$
Now the idempotents $f_{B,B}$ of $\End_{k\CL}(\sou{n})$, for $B\in\CP_{\sou n,m}$ and $0\leq m\leq n$,
are orthogonal and their sum is equal to the identity, by Theorem~\ref{end total}.
It follows that the endomorphisms $F_{f_{B,B}}$ of $F_{\sou{n}}$ are orthogonal idempotents, and their sum is the identity.
Hence we obtain a decomposition of correspondence functors
$$F_{\sou{n}}=\bigoplus_{\substack{ 0\leq m\leq n \\ B\in\CP_{\sou n,m} }} F_{f_{B,B}}\big(F_{\sou{n}}\big)\mpoint$$
By surjectivity of $\pi^B:\sou n \to \sou m$, the image of $F_{f_{B,B}}=F_{j^B}F_{\pi^B}:F_{\sou n} \to F_{\sou n}$ is equal to the image of~$F_{j^B}:F_{\sou m} \to F_{\sou n}$.
Therefore $F_{f_{B,B}}\big(F_{\sou{n}}\big)=F_{j^B}\big(F_{\sou{m}}\big)$.
By Proposition~\ref{ker-j}, the image $F_{j^B}\big(F_{\sou{m}}\big)$ is isomorphic to~$\S_m=F_{\sou{m}}/H_{\sou{m}}$ and it follows that
$$F_{f_{B,B}}\big(F_{\sou{n}}\big) \cong \S_m \mpoint$$
Taking $m=n$ and $f_{B,B}=j^B=\varepsilon_n$, we obtain the first isomorphism $F_{\varepsilon_n}F_{\sou{n}}\cong \S_n$.
Summing over all $B\in\CP_{\sou n,m}$ for a fixed~$m$, we obtain the second isomorphism, because $|\CP_{\sou n,m}|={n\choose m}$.
Finally, summing over all $0\leq m\leq n$ and all $B\in\CP_{\sou n,m}$, we obtain the third isomorphism.
\endpf

\result{Corollary} \label{Hom-zero} Let $m, n\in\N$. Then
$$\Hom_{\CF_k}(\S_n,\S_m)
=\left\{\begin{array}{ll}0&\hbox{if}\;n\neq m \mvirg\\
k\cdot\Id_{\S_n}&\hbox{if}\;n=m \mpoint \end{array}\right.$$
\fresult

\pf Since $\S_n\cong F_{\varepsilon_n}F_{\sou{n}}$, the case $n=m$ follows from Proposition~\ref{central-idempotents}. Now for integers $l,m\in\{0,\ldots,n\}$, we use the central idempotents $\beta_{n,l}$ of~Proposition~\ref{central-idempotents} and we obtain
$$\Hom_{\CF_k}(F_{\beta_{n,l}}F_{\sou{n}},F_{\beta_{n,m}}F_{\sou{n}})\cong \Hom_{\CF_k}(\S_l,\S_m)^{\oplus {n\choose l}{n\choose m}}\mpoint$$
Since $F_{\beta_{n,l}}$ and $F_{\beta_{n,m}}$ are central idempotents of $\End_{\CF_k}(F_{\sou{n}})$, and since they are orthogonal if $l\neq m$, it follows that $\Hom_{\CF_k}(F_{\beta_{n,l}}F_{\sou{n}},F_{\beta_{n,m}}F_{\sou{n}})=0$ if $l\neq m$, hence $\Hom_{\CF_k}(\S_l,\S_m)=0$.
\endpf

Now we prove that the functor $\S_n$ is actually isomorphic to a fundamental functor and we compute the ranks of all its evaluations.

\result{Theorem} \label{S_n-fundamental} Let $\S_n=F_{\sou{n}}/H_{\sou{n}}$.
\begin{enumerate}
\item $\S_n$ is isomorphic to the fundamental functor $\S_{[n],\tot}$, where $\tot$ denotes the total order on~$[n]$.
\item For any finite set $X$, the $k$-module $\S_n(X)$ is free of rank 
$${\rm rank}\big(\S_n(X)\big)=\sum_{i=0}^n(-1)^{n-i}{n\choose i}(i+1)^{|X|}\mpoint$$
\end{enumerate}
\fresult

\pf 
(a) We are going to use the results of Section~\ref{Section-duality-lattices} applied to the lattice $T= \sou{n}\op$.
The set of its irreducible elements is
$$E=\{0,1,\ldots,n{-}1\} \mvirg$$
with a total order $\tot\op$ being the opposite of the usual order. Now we have
$$F_{T\op}=F_{(\sou{n}\op)\op}=F_{\sou{n}}$$
and its evaluation at~$E$ contains an element
$$\gamma_T=\gamma_{\sou{n}\op}=\sum_{A\subseteq E} (-1)^{|A|} \eta_A^\circ \mpoint$$
Recall from Notation~\ref{def-eta-gamma} that $\eta_A^\circ:E\to T\op=\sou{n}$ denotes the same map as
$\eta_A:E\to T=\sou{n}\op$ and that $\eta_A$ is defined by
$$\forall e\in E,\;\;\eta_A(e)=\left\{\begin{array}{ll}r(e)&\hbox{if}\;e\in A\mvirg \\
e&\hbox{if}\;e\notin A \mvirg \end{array}\right.
\quad\text{ that is, }\quad
\eta_A(e)=\left\{\begin{array}{ll}e+1&\hbox{if}\;e\in A \mvirg\\
e&\hbox{if}\;e\notin A \mvirg \end{array}\right.$$
because $r(e)=e+1$ in the lattice $\sou{n}\op$.\par

Now we define $\omega:E\to \sou{n}$ by $\omega(e)=e+1$.
Then $\omega\in F_{\sou{n}}(E)$ and when we apply the idempotent $F_{\varepsilon_n}$ we claim that we obtain
\begin{equation}\label{omega-gamma}
F_{\varepsilon_n} (\omega)=(-1)^n \gamma_{\sou{n}\op} \mpoint
\end{equation}
The definition of $\varepsilon_n$ yields
$$F_{\varepsilon_n} (\omega)=(-1)^n\sum_{Y\subseteq[n]} (-1)^{|C|} \rho_Y \, \omega$$
and the definition of $\rho_Y$ gives
$$(\rho_Y \, \omega)(e)=\rho_Y(e+1)
=\left\{\begin{array}{ll}e+1&\hbox{if}\;e+1\in Y \mvirg \\
e&\hbox{if}\;e+1\notin Y \mpoint \end{array}\right.$$
Setting $Y=A+1$ for each $A\subseteq E$, we see that $\rho_Y \, \omega = \eta_A^\circ$ and it follows that
$$\begin{array}{rcl}
F_{\varepsilon_n} (\omega)&=&
(-1)^n\displaystyle\sum_{Y\subseteq[n]} (-1)^{|Y|} \rho_Y \, \omega \\
&=&(-1)^n\displaystyle\sum_{A\subseteq E} (-1)^{|A|} \eta_A^\circ \\
&=&(-1)^n\gamma_{\sou{n}\op} \mpoint
\end{array}$$
This proves Claim~\ref{omega-gamma} above.\par

Now $F_{\sou{n}}$ is generated by $\omega\in F_{\sou{n}}(E)$, because it is generated by $\iota\in F_{\sou{n}}([n])$
(where $\iota:[n]\to \sou{n}$ is the inclusion), hence also by any injection from a set~$E$ of cardinality~$n$ to~$\sou{n}$
(by composing with a bijection between $E$ and~$[n]$).
Since $F_{\varepsilon_n}$ is an idempotent endomorphism of the correspondence functor $F_{\sou{n}}$,
we see that $F_{\varepsilon_n}F_{\sou{n}}$ is generated by $F_{\varepsilon_n} (\omega)$.
In other words, in view of Claim~\ref{omega-gamma} above, $F_{\varepsilon_n}F_{\sou{n}}$ is generated by~$\gamma_{\sou{n}\op}
\in F_{\sou{n}}(E)$.
Now Theorem~\ref{SERgamma} asserts that the subfunctor of~$F_{\sou n}$ generated by
$\gamma_T=\gamma_{\sou{n}\op}$ is isomorphic to~$\S_{E,\tot\op}$.
But $(E,\tot\op)\cong ([n],\tot)$ via the map $e\mapsto n-e$.
Therefore, using the isomorphism of Theorem~\ref{functor-Sn},  we obtain
$$\S_n\cong F_{\varepsilon_n}F_{\sou{n}}=\langle \gamma_{\sou{n}\op}\rangle
\cong \S_{E,\tot\op} \cong \S_{[n],\tot} \mpoint$$

(b) The canonical $k$-basis of the $k$-module $\S_n(X)=F_{\sou n}(X)/H_{\sou n}(X)$ is the set $Z_n(X)$ of all maps $\varphi:X\to \sou{n}$ such that $[n]\subseteq \varphi(X) \subseteq \sou{n}$.
Therefore $\S_n(X)$ is free of rank $|Z_n(X)|$.
The number of maps in~$Z_n(X)$ has been computed in Lemma~8.1 of~\cite{BT2} and the formula is actually well-known.
The formula shows that this rank is equal to
$$|Z_n(X)|=\sum_{i=0}^n(-1)^i{n\choose i}(n+1-i)^{|X|}=\sum_{j=0}^n(-1)^{n-j}{n\choose j}(j+1)^{|X|}$$
as required.
\endpf

\begin{rem}{Remark}
We shall see in~\cite{BT3} that a similar formula holds for the rank of the evaluation of any fundamental functor,
but the proof in the general case is much more elaborate.
Also, Corollary~\ref{Hom-zero} holds more generally for fundamental functors and the general case will be proved in~\cite{BT3}.
\end{rem}

\bigskip
When $k$ is a field, we get even more.

\result{Corollary} \label{simple-Sn}
Let $k$ be a field.
\begin{enumerate}
\item The functor $\S_n$ is simple, isomorphic to $S_{[n],\tot,k}$, where $k$ is the trivial module for the trivial group $\Aut([n],\tot)=\{\Id\}$. 
\item $\S_n$ is simple, projective, and injective in $\CF_k$.
\end{enumerate}
\fresult

\pf 
(a) It is clear that $\Aut([n],\tot)$ is the trivial group, with a single simple module~$k$.
Recall from Section~\ref{Section-functors} that the fundamental functor $\S_{E,R}$ and the simple functor $S_{E,R,V}$ are given by
$$\S_{E,R}=L_{E,\CP_E f_R}/J_{E,\CP_E f_R} \,, \qquad S_{E,R,V}=L_{E,T_{R,V}}/J_{E,T_{R,V}} \mvirg$$
where $\CP_E f_R$ is the fundamental module corresponding to $(E,R)$
and $T_{R,V}=\CP_E f_R\otimes_{k\Aut(E,R)}V$.
Here $E=[n]$, $R=\tot$, $k\Aut(E,R)=k$, and $V=k$, hence
$$T_{\tot,k}=\CP_{[n]} f_{\tot}\otimes_kk\cong \CP_{[n]} f_{\tot} \mpoint$$
Therefore $\S_{[n],\tot} \cong S_{[n],\tot,k}$.

(b) Since $\sou{n}$ is a distributive lattice, $F_{\sou{n}}$ is projective and injective by Corollary~\ref{self-injective}. 
Therefore so is its direct summand~$\S_n$. It follows that $\S_n$ is simple, projective, and injective.
\endpf

Our last purpose in this section is to find, for any finite lattice~$T$, all the direct summands of $F_T$ isomorphic to a functor~$\S_n$ corresponding to a total order.
Recall that $e_T$ denotes the central idempotent of~$\End_{k\CL}(T)$ which is an identity element for the subalgebra $\CE_T$ (see Theorem~\ref{central-e}).

\result{Theorem} \label{facteurs directs totaux}
Let $T$ be a finite lattice and let $N$ be the maximal length of a strictly increasing sequence in~$T$.
For every finite set $X$, let $F_T^{\rm tot}(X)$ be the $k$-submodule of~$F_T(X)$ generated by all the maps $\varphi:X\to T$ such that $\varphi(X)$ is a totally ordered subset of~$T$.
\begin{enumerate}
\item $F_T^{\rm tot}$ is a subfunctor of $F_T$, equal to $F_{e_T}(F_T)$.
\item $F_T^{\rm tot}$ is a direct summand of~$F_T$, isomorphic~to
$$F_T^{\rm tot}\cong \bigoplus_{\substack{ 0\leq m\leq N \\ B\in\CP_{T,m} }} \S_{|B|}\mpoint$$
\item The image of any morphism $F_{\sou n} \to F_T$ is contained in~$F_T^{\rm tot}$.
In particular, any subfunctor of $F_T$ isomorphic to a functor~$\S_n$ is contained in~$F_T^{\rm tot}$.
\item $\Hom_{\CF_k}\big(F_T^{\rm tot}, F_{\Id-e_T}(F_T)\big)=\{0\}$ and
$\Hom_{\CF_k}\big(F_{\Id-e_T}(F_T), F_T^{\rm tot}\big)=\{0\}$.
\item The splitting of the surjection $F_T\to F_T^{\rm tot}$ is natural in~$T$.
\end{enumerate}
\fresult

\pf
(a) Let $\varphi\in F_T^{\rm tot}(X)$. Writing $\Im(\varphi)=\{t_1,t_2,\ldots,t_n\}$ in increasing order, we can write $\varphi=j\psi$, where $\psi:X\to\sou n$ is the map defined by $\psi(x)=i$ if $\varphi(x)=t_i$, and where $j:\sou n\to T$ is the
map defined by $j(i)=t_i$ for $1\leq i\leq n$ and $j(0)=\widehat0$. Clearly $j$ is a morphism in the category~$\CL$.
By Theorem~\ref{central-e}, we have $e_T j=j e_{\sou n}=j$, because $e_{\sou n}$ is the identity element of~$\End_{k\CL}(\sou n)$ by Theorem~\ref{end total}. Therefore
$$\varphi=j\psi=e_T j\psi = e_T\varphi =F_{e_T}(\varphi) \mvirg$$
proving that $\varphi\in F_{e_T}(F_T(X))$.\par

Conversely, if $\varphi\in F_{e_T}(F_T(X))$, then we can write $\varphi= F_{e_T}(\psi)=e_T\psi$ for some map $\psi:X\to T$.
Since $e_T$ is, by construction, a linear combination of maps with a totally ordered image, so is $e_T\psi$,
proving that $\varphi\in F_T^{\rm tot}(X)$.\par

This shows that $F_T^{\rm tot}=F_{e_T}(F_T)$ and the latter is a subfunctor of~$F_T$.\mpn

(b) As in the proof of Theorem~\ref{functor-Sn}, we apply the fully faithful functor $k\CL\to\CF_k$ defined by $T\mapsto F_T$.
There is direct sum decomposition of functors
$$F_T=F_{e_T}(F_T) \oplus F_{\Id-e_T}(F_T) = F_T^{\rm tot} \oplus F_{\Id-e_T}(F_T) \mpoint$$
The idempotent $e_T$ is the sum of the orthogonal idempotents $f_{B,B}$ of~$\End_{k\CL}(T)$, for $B\in\CP_{T,m}$ and $0\leq m\leq N$.
It follows that the endomorphisms $F_{f_{B,B}}$ of $F_T$ are orthogonal idempotents with sum~$F_{e_T}$.
Hence we obtain a direct sum decomposition of correspondence functors
$$F_T^{\rm tot}=F_{e_T}(F_T)
=\bigoplus_{\substack{ 0\leq m\leq n \\ B\in\CP_{T,m} }} F_{f_{B,B}}\big(F_T\big) \mpoint$$
By Proposition~\ref{ker-j}, the image of~$F_{f_{B,B}}$ is isomorphic to~$F_{\sou m}/H_{\sou m}=\S_m$, where $m=|B|$, proving the result.\mpn

(c) Let $\alpha: F_{\sou n} \to F_T$ be a morphism of correspondence functors.
Since the functor $T\mapsto F_T$ is full, $\alpha$ is the image of a morphism $\sou n\to T$ in~$k\CL$, which is in turn a linear combination of order-preserving maps $f:\sou n\to T$.
For such a map~$f$ and for any function $\varphi:X\to \sou n$, the image of $f\varphi$ is a totally ordered subset of~$T$.
It follows that the image of the map $F_f:F_{\sou n}(X) \to F_T(X)$ is contained in $F_T^{\rm tot}(X)$.\par

The special case follows from the fact that $S_{\sou n}$ is a subfunctor of~$F_{\sou n}$, by Theorem~\ref{functor-Sn}.\mpn

(d) The first statement is a consequence of (b) and~(c), while the second one follows from a dual argument. Details are left to the reader.\mpn

(e) By Theorem~\ref{central-e}, the family of idempotents~$e_T$, for $T\in\CL$, is a natural transformation of the identity functor~$\Id_{k\CL}$. Therefore the family of idempotents~$F_{e_T}$, for $T\in\CL$, is a natural transformation of the identity functor~$\Id_{\CF_k}$.
\endpf


\bigskip
\noindent
Serge Bouc, CNRS-LAMFA, Universit\'e de Picardie - Jules Verne,\\
33, rue St Leu, F-80039 Amiens Cedex~1, France.\\
{\tt serge.bouc@u-picardie.fr}

\medskip
\noindent
Jacques Th\'evenaz, Section de math\'ematiques, EPFL, \\
Station~8, CH-1015 Lausanne, Switzerland.\\
{\tt Jacques.Thevenaz@epfl.ch}

\end{document}